   \documentclass[12pt]{article}

\usepackage{amsmath,amsthm,amssymb}

\begin{document}

\theoremstyle{plain}

\newtheorem{Thm}{Theorem}[section]
\newtheorem{Prop}{Proposition}[section]
\newtheorem{Defi}{Definition}[section]
\newtheorem{Cor}{Corollary}[section]
\newtheorem{Lem}{Lemma}[section]
\newtheorem{Rem}{Remark}[section]

\makeatletter\renewcommand{\theequation}{%
\thesection.\arabic{equation}}
\@addtoreset{equation}{section}\makeatother

\setlength{\baselineskip}{15pt}
\newcommand{\vlimsup}{\mathop{\overline{\lim}}}
\newcommand{\vliminf}{\mathop{\underline{\lim}}}
\newcommand{\spec}{{\rm spec}}
\def\textmc{\rm}
\def\({(\!(}
\def\){)\!)}
\def\R{{\bf R}}
\def\Z{{\bf Z}}
\def\N{{\bf N}}
\def\C{{\bf C}}
\def\T{{\bf T}}
\def\E{{\bf E}}
\def\H{{\bf H}}
\def\Prob{{\bf P}}
\def\M{{\cal M}}     
\def\F{{\cal F}}
\def\G{{\cal G}}
\def\D{{\cal D}}
\def\X{{\cal X}}
\def\A{{\cal A}}
\def\B{{\cal B}}
\def\L{{\cal L}}
\def\a{\alpha}
\def\b{\beta}
\def\e{\varepsilon}
\def\de{\delta}
\def\ga{\gamma}
\def\k{\kappa}
\def\la{\lambda}
\def\fa{\varphi}
\def\th{\theta}
\def\si{\sigma}
\def\t{\tau}
\def\om{\omega}
\def\De{\Delta}
\def\Ga{\Gamma}
\def\La{\Lambda}
\def\Om{\Omega}
\def\Th{\Theta}
\def\lan{\langle}
\def\ran{\rangle}
\def\lbr{\left(}
\def\rbr{\right)}
\def\const{\;\operatorname{const}} 
\def\dist{\operatorname{dist}} 
\def\Tr{\operatorname{Tr}}
\def\quadd{\qquad\qquad}
\def\n{\noindent}
\def\beq{\begin{eqnarray*}}
\def\eeq{\end{eqnarray*}}
\def\supp{\mbox{supp}}
\def\beqn{\begin{equation}}
\def\eeqn{\end{equation}}
\def\bp{{\bf p}}
\def\sg{{\rm sgn\,}}
\def\1{{\bf 1}}
\def\pf{{\it Proof.}}
\def\v2{\vskip2mm}
\def\n{\noindent}
\def\z{{\bf z}}
\def\x{{\bf x}}
\def\y{{\bf y}}
\def\bv{{\bf v}}
\def\bu{{\bf u}}
\def\be{{\bf e}}
\def\0{{\bf 0}}
\def\pr{{\rm pr}}
\def\cp{{\rm Cap}}
\def\tst12{{\textstyle\frac12}}

\begin{center}
{\bf  The Brownian hitting distributions in space-time 
 of  bounded sets 
  and \\ the expected volume of the Wiener sausage  for a Brownian bridge } \\
\vskip4mm
{K\^ohei UCHIYAMA} \\
\vskip2mm
{Department of Mathematics, Tokyo Institute of Technology} \\
{Oh-okayama, Meguro Tokyo 152-8551\\
e-mail: \,uchiyama@math.titech.ac.jp}
\end{center}

\vskip8mm
\n
{\it running head}:  Brownian hitting distributions in space-time

\vskip2mm
\n
{\it key words}: Brownian hitting time, caloric measure,  parabolic measure, Wiener sausage, Brownian bridge. 
\vskip2mm
\n

\vskip6mm

\begin{abstract}
The   space-time distribution, $Q_A(\x, dtd\xi)$  say, of  Brownian hitting of a bounded Borel set $A$  of $\R^d$ is studied.  We derive the  asymptotic form  of the leading term of the time-derivative  $Q_A(\x, dtd\xi)/dt$ for each $d =2, 3,\ldots$, valid uniformly with respect to the starting point  $\x$ of the Brownian motion, which result significantly  extends  the classical ones for $Q_A(\x, dtd\xi)$ itself by Hunt ($d=2)$, Joffe and Spitzer ($d\geq 3$).  The results obtained  are applied to find the asymptotic form of the expected volume of   Wiener sausage  for the Brownian bridge joining the origin to  a  distant point.
 \end{abstract}
\vskip6mm

\section { Introduction and summary of main results}

In  this paper we primarily focus on the  space-time distribution of  Brownian hitting of a bounded Borel set $A$   of $\R^d$ expressed as 
\beqn\label{Eq1}
Q_A(\x,dt d\xi) = P_{\x}[ B_{\sigma(A)}\in d\xi, \sigma_A\in dt]\quad (t>0, d\xi\subset \partial A).
\eeqn
 Here $ P_\x$ denotes the law of  a standard Brownian motion $B_t$  started at $\x$, $\sigma_A $ (or $\sigma(A)$)  
the first hitting time of $A$ by $B_t$, namely  $\sigma_A= \inf\{t>0: B_t\in A\}$,  and $\partial A$ the Euclidian boundary of $A$. When $A$ is a disc ($d=2$) or a ball ($d\geq 3$), the distribution $P_\x[\sigma_A <t]$ or its density $P_\x[\sigma_A \in dt]/dt$ are investigated by several recent works  \cite{HM},   \cite{BMR}, \cite{Ubh},\cite{Ubes}  seeking   the asymptotic behavior of them  for $\x \notin A$ as $t\to \infty$.
For  general $A$  the asymptotic  form of the distribution $P_\x[\sigma_A <t]$ is given by Hunt \cite{Ht} for $d=2$ and by Joffe \cite{J} and  Spitzer \cite{S2} for $d\geq 3$. 
Their classical results  may read as follows: if $A$ is compact and $\R^d\setminus  A$ is connected, then  for each $\x\in \R^d\setminus  A$,  as $t\to \infty$
\beqn\label{Hunt}
P_\x[t< \sigma_A] \sim \frac{2 e_A(\x)}{\lg t} \quad\quad  d=2,
\eeqn
\beqn\label{J-S}
P_\x[t<\sigma_A< \infty] \sim  \cp(A) P_\x[\sigma_A =\infty] \frac{t^{ -d/2+1}}{(2\pi)^ {d/2}(d/2-1) }\quad\quad d\geq 3.
\eeqn
  where  $e_A$ is a  Green function for  the open set $\R^2\setminus  A$ with a pole at infinity  and $\cp(\cdot)$ a Newtonian capacity
(see (\ref{i}) and  (\ref{ii}) below).   The latter result  is uniform for $\x$ in each compact set of  $\R^d \setminus  A$.   Of the former one there is given an elementary proof in  \cite{CMM},   which ensures the uniformity, while Hunt's proof, resting on a  Tauberian theorem, does not.  In \cite{Ht},  \cite{J} and  \cite{S2}  $A$ is assumed compact as in above but the extension to bounded Borel  sets is immediate with a  knowledge about  measurability question of the hitting time \cite{Ht2}, \cite{BG}, \cite{B}.  M.  van den Berg \cite{vB} recently improves  the result in the case $d\geq 3$  by obtaining sharp remainder estimates. 
In this paper we derive explicit  asymptotic forms  of  the time derivative $Q_A(\x,dt d\xi)/dt$ as $t\to\infty$ valid uniformly for $|\x|=o(t)$ or for  $t=o(|\x|)$.  
By easy integration  the  asymptotic forms of  $P_\x[t<\sigma_A<\infty]$ and $P_\x[\sigma_A<t|\sigma_A <\infty]$  can be computed   from our results given shortly in this introduction (see Appendix A.2).

 The measure kernel $Q_A(\x,dtd\xi)$ plays a significant role in the theory of  heat operator. If  $A$ is compact  and $\Om_A$ denotes the unbounded component of $\R^d \setminus  A$, then $Q_A(\x,dtd\xi)$ is identified with   the lateral part of the caloric measure (or parabolic measure) for the heat operator $\frac12\De -\partial_t$ in  the space-time domain
$D=\{(\x,t)\in \R^d\times (0,\infty) :  \x\in \Om_A\}$, the exterior of a cylinder (see Appendix   A.1).  The other part of it  is nothing but the measure whose density is given by the heat kernel for $\Om_A$ with Dirichlet zero boundary condition and its (uniform)  asymptotic estimate for large time  is recently  obtained by \cite{CMM} for  the space variables restricted to any compact set and by \cite{Udsty} without  restriction. 
The present work is partly motivated and steered by  a  study  of   Wiener sausage  swept by the set $A$ attached to a $d$-dimensional  Brownian motion started at the origin. 
Our interest  is in finding a  correct  asymptotic form of the expected volume of  the sausage of length $t$ as $t\to\infty$ under  the conditional law given that the Brownian motion at time $t$  is  at a given site  $\x$ which is  outside a parabolic region so that $|\x|^2> \e t$ for  any positive  $\e$. To this end  it is  needed in our approach to  estimate 
the density  $Q_A(\x,dt d\xi)/dt$;  we shall derive  asymptotic forms of the expected volume by applying the results  on $Q_A$.  Fine estimates are obtained in the case when the process is pinned at the origin  by
McGillivray \cite{M} ($d\geq 3$) and \cite{M2} ($d=2$) (cf. also \cite{BB}) and in the case when $d=2$, the value at time $t$ is   pinned  within 
a parabolic region and  $A$ is a disc  by  \cite{Usaus}. 
There are many works for the sausage in the unconditional  case (see, e.g.,  the references of \cite{BB}).  
    
\v2
  Let  $U(a)\subset \R^d$ denote  the open ball  about the origin of radius $a>0$.  Let $A$ be a bounded Borel set as above,  denote  by $A^r$ the set of all regular points of $A$, i.e., the set  of $\y\in \bar A$ such that
$P_\y[\sigma_A=0]=1$ and put 
$$\Om_A = \{\x\in \R^d : P_\x[\sigma_{\partial U(R)}< \sigma_A] >0 \} $$
with   any $R$ such that $\bar A\subset U(R)$.
 (The over bar designates  the Euclidean closure: $\bar A = A \cup \partial A$, where $\partial A$ denotes the  Euclidean boundary of $A$.)  
The set  $A^r$  is Borel 
 and   $A\setminus A^r$ is polar, so that  $P_\x[\sigma_A =\sigma_{A^r}]=1$ and  $P_\x[B_{\sigma(A)}\in A^r|\sigma_A<\infty]=1$ for all $\x$ (see e.g. \cite{PS}, \cite{B}); it is natural from the view point of intrinsic topology 
to consider  $\Om_A$, which agrees with  the unbounded  fine component of $\R^d\setminus {A^r}$ (see \cite{C}, p.370 for what   \lq fine component' means).
We always suppose that $A^r$ has no \lq cavities' isolated from $\Om_A$: 
\beqn\label{cavity}
{\bf R}^d \setminus \Om_A = A^r,
\eeqn
 which will require no  change of the intrinsic content of the paper,
  Brownian paths  started at a point of $\Om_A$ being kept out of the fine interior of ${\bf R}^d \setminus \Om_A $ before hitting  $A$. 

Define for $\x\in \Om_A$,
$$q_A(\x,t) =\frac{d}{dt}P_\x[\sigma_A \leq t] $$
and
\beqn\label{H_A}
H_A(\x,t; d\xi) = \frac {Q_A(\x,dt d\xi)}{dt} =\frac{P_{\x}[ B_{\sigma(A)}\in d\xi, \sigma_A\in dt]}{dt}\quad\quad (d\xi \subset \partial A, t>0).
\eeqn
(Although both $q_A$ and $H_A$ depend on $d$, we do not designate this dependence in the notation.)
Our purpose of this  paper is to find the asymptotic form of  $H_A(\x,t; \cdot)$ as $t\to\infty$. To this end 
it is often convenient to factor it   into a  product by conditioning on $\sigma _A$ as follows:
\beqn\label{factor}
H_A(\x, t; d\xi)= P_\x[ B_{\sigma(A)}\in d\xi|\sigma _A=t]q_A(\x,t).
\eeqn
Put 
$$ p_t^{(d)}(x)   = (2\pi t)^{-d/2} e^{-x^2/2t}\quad\quad (x=|\x|)$$  
and for $d\geq 3$
$$ G^{(d)}(x) = \int_0^\infty p_t^{(d)}(x) dt = \frac{\Ga(\nu)}{2\pi^{\nu+1}x^{2\nu}}\quad (\nu=\textstyle{\frac12} d -1).$$

We usually write 
$x$ for $|\x|$, $\x\in \R^d$  (as above)  and  sometimes $p^{(d)}_t(\x)$ for $p^{(d)}_t(x)$  and $\sigma(A)$ for $\sigma_A$ unless doing these causes any confusion. Denote   by nbd$_\e(A)$ the open    $\e$-neighborhood of  $A$ in $\R^d$. We may write $\x\notin  {\rm nbd}_\e(A^r)$ instead of $\x\in \Om_A\setminus {\rm nbd}_\e(A^r)$  in view of (\ref{cavity}). 
We write $f(t) \sim g(t)$ if
 $f(t)/g(t)\to 1$ in any  process of taking limit like $\lq t\to\infty$\rq. 

The results of this paper are summarized in the following propositions  (i) through (viii), where   $A$ is   bounded Borel and non-polar, i.e.,  $P_\x[\sigma_A<\infty]>0$ for some $\x$. We   indicate in  square brackets at the head of each statement the theorem given in a succeeding  section  which the result presented below is taken from. 

\v2\n
{\sc  Summary of Main Results.} \, 
\v2

I.  Case $x/t\to 0$.

 {\it For each $\e>0$, the following asymptotic formulae}  (i) {\it and}  (ii) {\it hold  uniformly for }$\x \notin$ nbd$_{\e}(A^r)$:
\v2

(i)  [Theorem \ref{thm3.1}]\,  {\it   if $d\ge 3$,  
as  $x/t\to 0$  and  $t\to\infty$   
\beqn\label{i}
q_A(\x,t) =\cp(A) P_\x[\sigma_A=\infty] p^{(d)}_t(x) (1+o(1)),
\eeqn
where $\cp(A)$ denotes the Newtonian capacity 
of $A$  (see Section 3.1.1)  normalized so that }
$$\cp(U(a))={a^{d-2}}/{G^{(d)}(1)};$$

(ii)  [Theorem \ref{thm3.2}]  \, {\it  if $d=2$,   as $x/t\to 0$  and  $t\to\infty$ 
\begin{eqnarray}\label{ii}
\qquad q_A(\x,t) = p^{(2)}_t(x) 
\times
\left\{ \begin{array} {ll}  
 {\displaystyle \frac{4\pi e_A(\x)\,}{(\lg t)^2}   \bigg(1+ O\Big(\frac1{\lg t}\Big)\bigg) } \quad& (x  \leq \sqrt {t}\,),\\ [5mm]
 {\displaystyle  \frac{\pi } {\lg (t/x) }    \bigg(1+ O\Big(\frac1{\lg (t/x)}\Big)\bigg) } \quad&(  x   > \sqrt t\,),
  \end{array} \right.   
\end{eqnarray}
where $e_A(\x)$ denotes  the Green function for $\Om_A$ with a pole at infinity  normalized so that 
$$e_{U(a)}(\x) = \lg (x/a)$$
 (cf. \cite{Kll}, p.369; see also  Section 6.1);  } 

\v2
(iii)  [Theorems \ref{thm3.3} and \ref{thm3.40}] \, {\it for  $d\geq 2$,  as  $x/t\to 0$  and  $t\to\infty$ 
  \beqn\label{iii}
\Big|P_\x[B_{\sigma(A)}\in \cdot\,\,|\, \sigma_A = t]  - H^\infty_A \Big|_{\rm t.var} \longrightarrow \, 0.
\eeqn
 Here,   $H^\infty_A$ stands for the harmonic measure for the Brownian motion started at infinity  conditioned on 
$\sigma_A<\infty${\rm:}  $H^\infty_A(d\xi) =\lim_{x\to\infty} P_\x[B_{\sigma(A)} \in d\xi\,|\,
\sigma_A<\infty]${\rm ;} and $|\cdot|_{\rm t.var}$ designates the total variation.}

\v2
Given a   compact set $K$ of $\R^d$, let $S_K(t)$ denote the {\it Wiener sausage} of length $t$ swept by  $K$ attached to a Brownian motion $B_t$ :
$$S_K(t) = \{\z\in \R^d: \z- B_s  \in  K \quad\mbox{for some} \quad s\in [0,t] \}.$$
The $d$-dimensional volume of  $A\subset \R^d$ is denoted  by vol$_d(A)$.  Suppose $K$ is non-polar. Then

\v2
 (iv) [Theorem \ref{thm2}]\,  {\it if   $d\geq3$, as \, $x/t\to 0$  and  $t\to\infty$ }
  \beqn\label{EQ01}
E_{{\bf 0}} \,[\, {\rm vol}_d(S_K(t)) \,|\, B_t=\x]  \,\sim\, \cp(K)t;
\eeqn

(v) [Theorem \ref{thm3}]\, {\it  if $d=2$,   $x/t\to 0$  and  $t\to\infty$  }
   \begin{eqnarray}\label{EQ02} \qquad
E_{{\bf 0}} \,[\, {\rm vol}_2(S_K(t))  \,|\, B_t=\x]  \,=\,  \left\{ \begin{array}{ll} 
{\displaystyle
\frac{2 \pi t }{\lg t} (1+o(1)) }\quad\quad&\mbox{if}   \quad x\leq  \sqrt t, \\[5mm]
 {\displaystyle
\frac{\pi t}{\lg (t/x)} (1+o(1)) }\quad\quad&\mbox{if}   \quad x> \sqrt t.
\end{array} \right.
\end{eqnarray}

 II. Case $x/t\to v$.   
 
 In the case when  $x/t$ is bounded away from zero and infinity the results are not explicit  as above. 
Put $R_A = \sup \{|\y|:y  \in A^r\}$ and  
   define the  measure kernel  $\la_A({\bf v}; d\xi)$ on $\R^d\times \partial A$ by 
$$\la_A({\bf v}; \Ga) = \int_{\partial U(1)} E_{R_A\xi}\Big[e^{-\frac12 v^2\sigma_A}; B_{\sigma(A)}\in \Ga \Big] g_{R_Av}(\th_{\xi,{\bf v}}) m_{1}(d\xi) \quad ({\bf v} \in \R^d, \Ga \subset \partial A),$$
if $v :=|{\bf v}|>0$,   where $\th=\th_{\xi,{\bf v}}\in [0,\pi]$ denotes the angle that  $\xi\in \partial U(1)$ forms with  ${\bf v}$ so that $\cos \th= \xi\cdot {\bf v}/v$; and  $g_\a(\th)$, $\a\geq 0$ (given  in (\ref{h_v}) of Section 2)  is a probability density  relative to the uniform probability measure  $m_1(d\xi)$ on $\partial U(1)$; if $v =0$, replace $g_{R_Av}(\th_{\xi,{\bf v}})$ by unity (which amounts to passing to the limit as $v\downarrow 0$).
 It is shown that  $\la_A({\bf v}; \partial A)$ is positive,  continuous   and bounded   (in ${\bf v}$). 
 We can assert the following 
\v2

(vi)\,  {\it  Let $d\geq 2$.  As  $x/t \to v >0$ and $t\to\infty$,
$$\bigg|\frac{H_A(\x,t;\cdot)}{R_A^{2\nu}\La_\nu(R_Ax/t)p_t^{(d)}(x)} - \la_A(\x/t;\cdot)\bigg|_{\rm t.var} \longrightarrow  \, 0$$
and for $K$ compact and non-polar,
$$ E_{{\bf 0}} \,[\, {\rm vol}_d(S_K(t)) \,|\, B_t=\x]  \,\sim\,  t R_K^{2\nu}\La_\nu(R_Kx/t)\int_{\partial K} e^{-\xi\cdot\x/t} \la_K(\x/t;d\xi),$$
with $\La_\nu(y)$  a positive function of $y$ given in (\ref{Lammda}) of the next section.
Here,  the convergence is  uniform in $v\leq M$ in both the formulae for each $M>1$ and they become continuously linked at $v=0$  to the formulae of {\rm (iii)} through {\rm (v)}. }  

{\it For the balls or discs the last  formula  may be given somewhat more explicitly:}
\beqn\label{EQ24}
E_{{\bf 0}} \, [ {\rm vol}_d(S_{U(a)}(t)) \,|\, B_t=\x] \,\sim\, a^{2\nu}\,t \La_\nu(av)\int_{\partial U(1)} e^{- a\xi\cdot\x/t} g_{av}(\th_{\xi,\x})m_1(d\xi).
\eeqn

\v2 
 III. Case $x/t\to \infty$.

Let $\be$ be  a unit vector of $\R^d$.   Denote by  $\De_{\be}$   the hyper-plane perpendicular to $\be$  passing through the origin and   ${\rm pr}_{\be}A$  the orthogonal projection of a set  $A$ on  $\De_\be$. Let $K$ be a compact set. Define a function $h=h_{\be, K}$ on $\De_\be$ by
$$h(\z)= \sup\{s\in \R: \z+ s\be \in  K\},  \quad\quad \z\in \De_{\be}$$
with $\sup \emptyset = -\infty$  and let $\mbox{dis-ct}_{\be}(K)$ be the set of discontinuity points of $h$ (see (\ref{mmm})).
Suppose that 
$$
{\rm vol}_{d-1}\Big(\overline{\mbox{dis-ct}_{\be}(K)}\Big)=0,$$
where  ${\rm vol}_{d-1}(\cdot)$  denotes the $(d-1)$-dimensional volume of a set of $\De_\be$.  Then:
\v2

(vii)  [Theorem \ref{thm3.4}]\, {\it  for  $d\geq 2$, as $v:= x/t\to \infty$ and  $t\to\infty$  
\v2

  \beqn\label{iv}
\frac{H_K(x\be, t;  d\xi)}{vp_t^{(d)}(x)e^{v\be\cdot\xi} } \,\Longrightarrow \,  m_{K,\be}(d\xi),
\eeqn
where \lq $\,\Rightarrow\!$'  designates the weak convergence of measures and $m_{K,\be}$ stands for the Borel measure on $\partial K$ induced from the $(d-1)$-dimensional Lebesgue measure on $\pr_\be K\subset \De_\be$  by  the mapping $\z\in \pr_\be K \mapsto \z+h(\z)\be \in \partial K$}  (see  Section 3.2 for more details); 
\v2

(viii) [Theorem \ref{thm4.3}]\,    {\it   for $d\geq 2$, as  $x/t\to\infty$ and $t\to\infty$, }
$$E_{{\bf 0}} \,[ \, {\rm vol}_d(S_K(t)) \, | \, B_t= x\be] = {\rm vol}_{d-1}({\rm pr}_{\be}K) x +o(x).$$

\v2
For convenience of later citation here we record  the following scaling properties:
\begin{eqnarray}\label{scp}
\left.
\begin{array}{c}
H_A(\x,t;  d\xi) = R^{-2}H_{R^{-1}A}(\x/R, t/R^2; R^{-1}d\xi) \quad (d\xi\subset \partial A),\\ [2mm]
q_A(\x,t) = R^{-2}q_{R^{-1}A}(\x/R, t/R^2),\\[2mm]
 e_A(\x)= e_{R^{-1}A}(\x/R) \quad (d=2) \quad \mbox{and} \\[2mm]
 \cp(A)= R^{d-2}\cp(R^{-1}A) \quad  (d\geq 3) 
 \end{array}\right.
 \end{eqnarray}
 ($R>0$; the heading factor $R^{-2}$ on the right-hand sides of the first two identities comes simply from the differential $d (t/R^2)$). Note that  the formula (\ref{i}) is consistent to these  relations.  

The rest of the paper is organized into the following sections. 
 \v2
  
  \S 2.  The hitting  distribution for a disc/ball.
  
  \S  3. The hitting  distribution for a bounded Borel set.
 
 \quad 3.1. Case $x/t\to 0$.\quad 3.2. Case $x/\sqrt t\to\infty$.
 
  \S  4.   The Wiener sausage for a Brownian bridge.
  
  \quad 4.1. Case $x/t\to 0$.\quad 4.2.  Case $x/\sqrt t\to\infty$.
  
 \S  5.  Brownian motion with a constant drift. 
   
  \S 6. Miscellaneous estimates concerning  $\sigma_A$.
  
 \quad   6.1.  Uniform   estimates for  \;  $e_A(\x)/P_\x[\sigma_{\partial U(r)} < \sigma_A] .$ 
 
 \quad  6.2.  An upper bound of $q_A$ $(d\geq 3)$.  \quad   6.3.  Some   \quad  upper  bounds of $q_A$ $(d = 2).$  
 
\quad  6.4.  A lower bound of  $P_\x[\sigma_A<t]$ in case   $x/t>1$.
 
   \S  Appendix.
   
  \quad A.1.  Harmonic measure of heat operator.\, \quad  A.2.   Asymptotics of the distribution of  $\sigma_A$.

\v2  
  In \S 2  we summarize the results of \cite{Ubes} and \cite{Ucal_b} that are particularly relevant to the present subject. 
The results of both \S 3 and \S 4 heavily depend on those from \S  2. 

The  subjects of \S 3 and \S 4 are inter-related. Because of this  the item (vi) is proved in Section 4.2.  The asymptotic behavior of $H_A(\x, t;  d\xi)$ is closely related to that of the expected volume of Wiener sausage swept by a compact set  $K$ attached to Brownian bridge joining ${\bf 0}$ and $\x$. If $x/t \to 0$, the results on $H_A$ entails those of the sausage. In the case $x/t \to\infty$, the situation is not so simple; our proof of (vii)  relies on a result on the upper estimate of the expected volume of the sausage, and the  lower estimate of it is obtained by using  (vii).

Most  of the statements advanced above may  translate into the ones corresponding  to Brownian motion with a constant drift and some  of them  will be  presented in \S 5.  In \S 6  we prove miscellaneous results that we need to use in the proofs of main results; some asymptotic evaluation  of $P_\x[\sigma_{\partial U(r)}<\sigma_A]$ as $r\to\infty$  in terms of  $e_A(\x)$; the upper bounds of $q_A$ that are  fundamental  to the proofs of the main results (i) through (viii); some estimates of $P_\x[\sigma_A<t]$ used in \S 4.2.3.
In Appendix we provide a brief  exposition of the  well-known relation of the hitting distribution to caloric measure and 
some asymptotic forms  of $P_\x[\sigma_A \leq t]$ which are derived  by elementary computation  from those of the density exhibited above.

\section{The hitting   distribution for a disc/ball}
 Here we consider the case when  $A=U(a)$,  the open ball  centered at the origin of radius $a$, and state   the  asymptotic estimates   of  $H_{U(a)}(\x,t;d\xi)$ obtained in \cite{Ubes} and \cite{Ucal_b}.   The following notation is used throughout the paper.  
 \begin{eqnarray}\label{Lammda}
&&\nu =\frac{d}{2}-1\quad  (d=1,2,\ldots); \nonumber
\\
&& \La_\nu(y)= \frac{(2\pi)^{\, \nu+1}}{2y^{\nu}K_{\nu}(y)}\quad(y>0);\quad\quad
 \La_\nu(0) = \lim_{y\downarrow 0}\La_\nu(y) \quad\quad (\nu>0).
 \end{eqnarray}
Here $K_\nu$ is the usual modified Bessel function (of the second kind) of order $\nu$.  
 We write 
$$q(x,t;a)$$
for $q_{U(a)}(\x,t)$.  
The definition of   $q(x,t;a)$    may be naturally extended to  the Bessel process of order $\nu$ and the results concerning it given below may be applied to such extension if $\nu\geq 0$.
The following result from \cite{Ubes} provides a precise asymptotic form of the hitting time density for a ball ($d\geq 3$) (or  disc ($d=2$)).

\begin{Thm}\label{thm01}\,  Uniformly for $x > a$, as $t\to\infty$, 
\beqn\label{R2}
q(x,t;a)  \, =\, a^{2\nu}\La_\nu\bigg(\frac{a x}{t}\bigg) p^{(d)}_t(x) \bigg[1-\bigg(\frac{a}{x}\bigg)^{2\nu}\bigg](1+o(1))\quad\quad \mbox{if}\quad  d\geq 3,
\eeqn
and
\beqn\label{R21}
q(x,t;a) = p^{(2)}_t(x) 
\times
 \left\{ \begin{array} {ll}  
  {\displaystyle \frac{4\pi\lg(x/a)\,}{(\lg (t/a^2))^2}\Big(1+o(1)\Big)    } \quad& (x  \leq \sqrt t\,)\\ [5mm]
 {\displaystyle  \La_0\bigg(\frac{a x}{t}\bigg) \Big(1+o(1)\Big)  }\quad&(  x   > \sqrt t\,)   \end{array} \right.
 \quad \mbox{if}\quad  d =2.
\eeqn
If the  right-hand sides are multiplied by $e^{-a^2/2t}$, both the formulae (\ref{R2}) and (\ref{R21}) so modified  hold true  also as $x\to\infty$ uniformly for $t>0$ (see also (\ref{result6}) below).
\end{Thm} 
\v2
 From  known results  on $K_\nu(z)$ it follows that 
\beqn\label{K01}
 \La_\nu(y)= \frac{2\pi}{\int_0^\infty  \exp(-\frac1{4\pi u}{y^2 })e^{-\pi u } u^{\nu-1}du}~~~~~~~~~  (y>0, \nu\geq 0);
\eeqn
 \[
\La_\nu(y) =(2\pi)^{\nu+1/2} y^{-\nu+1/2}\, e^{y} ( 1+O(1/y)) \quad \mbox{as}\quad y\to\infty;
\]
\[
\La_\nu(0)  = 1/G^{(d)}(1)= {2\pi^{\nu+1}}/ {\Ga(\nu)} = (d-2)\pi^{d/2}/\Ga(d/2)
\quad\mbox{for}\quad\nu >0; \,\mbox{and}
\]
\[\La_\nu(y) = \left\{ \begin{array}{ll}  {\displaystyle \frac{\pi}{-\lg (e^\ga y/2)} }(1 +O(y^2)) & (\nu=0)\\[4mm]
\La_\nu(0) + O(y^2)  \quad &(\nu>0, \nu\neq 1)   \\
\La_\nu(0) + O(y^2\lg y)  \quad &(\nu=1)    
\end{array} \right. \quad \mbox{as}\quad y \downarrow 0.
\]
Here $\ga =\int_0^1(1-e^{-t} -e^{-1/t})t^{-1}dt $ (Euler's constant).
\v2
\v2

The following  two results are also valid for Bessel processes of order $\nu\geq 0$. The first one  is easily deduced from Theorem \ref{thm01} by elementary computation (see  the last section  of \cite{Ubes}).
The second one (a reduced version of \cite[Lemma 4]{BMR})  is also easily derived from the one-dimensional result with the help of a drift-transformation formula (\cite[(12)]{Ubes}). 
\begin{Thm}\label{thm2.2}\, Uniformly for $x>a$ if $\nu>0$ and  uniformly for  $x > \sqrt{t/\lg t}$ if $\nu=0$, as $t\to \infty$, 
$$\bigg(\frac{a}{x}\bigg)^{-2\nu}\int_0^t q(x,s;a)ds = \La_{\nu}\bigg(\frac{ax}{t}\bigg)
\frac{2^\nu }{(2\pi )^{\nu+1}}\int_{x^2/2t}^\infty e^{-y}y^{\nu-1}dy(1+o(1)).$$
(Note that $(a/x)^{2\nu} = \int_0^\infty q(x,s;a)ds $, so that the left side represents a conditional probability.)
\end{Thm}

\begin{Lem}\label{thm5} \,  For each $\nu\geq 0$ it holds that  uniformly
 for  all $0< t<a^2$ and $x>a$, 
\beqn\label{result6}
 q(x,t;a) =  \frac{x-a}{\sqrt{2\pi}\, t^{3/2}} e^{-(x-a)^2/2t}\bigg(\frac{a}{x}\bigg)^{(d-1)/2}\bigg[1+ O\bigg(  \frac{t}{x}\bigg)\bigg].
 \eeqn
 \end{Lem}

A trite computation shows  that in case  $x/t\to\infty$ the function form of the leading term  for $q(x,t;a)$ as $t\to \infty$   given above coincides  with the one  for  $t\to\infty$ given in Theorem \ref{thm01} so that (\ref{result6}) also holds in this case if
$O(t/x)$ is replaced by $o(1)$.
 
The case  $x^2/t\to 0$  with $\lim (\lg x)/\lg t =1/2$ of $\nu=0$ (i.e. $d=2$), not included in Theorem \ref{thm2.2},  is somewhat delicate.    The following result is a reduced form of Theorem 3 of \cite{Ubh}.
\vskip2mm
 \begin{Thm}\label{thm2.20}\, Let $\nu=0$ and $\k= 2e^{-2\ga}$. Uniformly for $\sqrt t> x>a$,  as $t\to\infty$
$$
\int_0^t q(x,s;a)ds
=  \frac{1}{\lg (\k t/a^2)}\bigg[1-\frac{\ga}{\lg (\k t/a^2)}\bigg]\int_{x^2/2t}^\infty\frac{e^{-y}}{y}dy +  O\bigg(\frac1{(\lg t)^2}\bigg).
$$
\end{Thm}
\v2

We sometimes need only the upper  or lower bounds  that follow immediately  from Theorem \ref{thm01}     and for convenience sake we write   down them. In the following corollary (obtained also  in \cite{BMR}) we 
include  the bounds for the case  $0< t<a^2$  and $d\geq 3$ that follows from Lemma \ref{thm5}.
We write $x\vee y$ and $x\wedge y$ for the maximum and  minimum of real numbers $x, y$, respectively. For positive functions $f$ and $g$ defined on a set  $\La$, the expression  $f(\la) \asymp g(\la)$ signifies that  there exist  constants $c_1$ and   $c_2$ such that  $c_1 g(\la) \leq  f(\la)\leq c_2 g(\la)$ for  $\la\in \La$.

\v2
\begin{Cor} \label{upper_bd} For all $a>0, t> 0, x>a$,
$$
q(x,t;a) \asymp  \Big(1-\frac{a}{x}\Big)a^{2\nu}p^{(d)}_t(x-a)\bigg(1\wedge \frac{t}{ax}\bigg)^{\nu-\frac12} \quad\quad\quad \mbox{if}\quad  d\geq 3,
$$
and
$$q(x,t;a) \asymp \left\{\begin{array}{ll} {\displaystyle  \frac{\lg x/a}{t(\lg t/a^2)^2} }&\quad   \quad  (x<\sqrt t) \\[4mm]
{\displaystyle \frac{1}{1+ \lg (t/ax)} \, p^{(2)}_t(x) }&\quad   (\sqrt t \leq x<t/a) \\[4mm]
{\displaystyle \sqrt{ \frac{ax}{t}} \, e^{ax/t}p^{(2)}_t(x) }&\quad  (x\geq t/a >2a)
\end{array}\right. \quad  \qquad \mbox{if}\quad  d =2,$$
where all the  constants involved in the symbol  $\asymp$ depend only on $d$. 
\end{Cor}
\v2
   
We also deduce the following corollary  of Theorem \ref{thm2.2}   (cf.  \cite[Theorem 12]{Ubes}) in the case $d\geq 3$. (See (\ref{6.9}) and (\ref{6.10})  for the corresponding results for $d=2$.)
\begin{Cor} \label{upper_bd2}  \, Let $d\geq 3$ and $a>0$.    Then
 \begin{eqnarray}\label{6.9_d_Int}
 P_\x[\sigma_{U(a)}< t] \asymp  \left\{ \begin{array} {ll}  (a/x)^{2\nu} \quadd \qquad\,\, &\mbox{if}\quad t>x^2>a^2,  \\[2mm]
{\displaystyle \frac{a^{2\nu}t^2}{x^2} \La_\nu\bigg(\frac{ax}{t}\bigg) p_t^{(d)}(x)} \quad &\mbox{if}\quad x^2\geq t>a^2.
 \end{array}\right.
 \end{eqnarray}
 Here all the  constants involved in the symbol  $\asymp$ depend only on $d$. 
  \end{Cor}

\v2
Next we present  a result from \cite{Ucal_b}  on the conditional distribution of the hitting site $B_{\sigma(U(a))}$ given  $\sigma_{U(a)}=t$.  For $\a\geq 0$ let $g_\a(\th)$ be a function of  $\th\in [-\pi,\pi]$ given by
\beqn\label{h_v}
g_\a(\th) = \sum_{n=0}^\infty \frac{K_0(\a)}{K_n(\a)} H_n(\th)
\eeqn
if  $\a>0$ and  $g_0(\th) \equiv1$,
where  
$$H_n(\th) = \left\{\begin{array} {ll}  \cos n \th   \quad &\mbox{if} \quad d=2,\\
\k_{n,\nu} C_n^\nu(\cos \th) \quad 
  &\mbox{if} \quad d \geq 3.
 \end{array}\right.$$ 
Here  $\k_{n,\nu} = (n+\nu)\Ga(\nu)/\sqrt \pi \, \Ga(\nu+ \frac12)$ and  $C_n^\nu(z)$ is  the Gegenbauer polynomial of order $n$ associated with $\nu$ (cf. \cite{W}). It follows that $H_0 \equiv 1$ and  $g_\a(\th)$ is jointly continuous in $(\a,\th)$. Let $\th_{\xi,\x}\in [0,\pi]$ denote the colatitude of a point $\xi\in \partial U(a)$ with the  vector $a\x/x$ taken to be the north pole, namely  $\cos \th_{\xi,\x} = \xi\cdot \x/ax$.
Let  $m_a(d\xi)$ denote the uniform probability measure on  $\partial U(a)$. It follows that $g_\a(\th_{\xi,\x})m_a(d\xi)$ is a probability measure on $\partial U(a)$.

\begin{Thm}\label{thm2.3}  Let   $d\geq 2$. The function  $g_\a(\th)$ is positive on  $[0,\pi]$, and  for  $v\geq 0$,
as  $x/t\to v$ and $t\to\infty$,  
\beqn\label{s-dst}
\frac{P_{\x}[B_{\sigma(U(a))} \in a d\xi\,|\, \sigma_{U(a)}=t]}{m_1(d\xi)}  =
\left\{ \begin{array}  {ll} 1 + O(\ell_d(x,t)x/t)\quad &\mbox{if} \quad v=0,\\[2mm]
 g_{av}(\th_{\xi,\x})(1+o(1))    \quad&\mbox{if}  \quad v>0
 \end{array}\right.
\eeqn
 uniformly for $(\xi, v) \in  \partial U(1) \times [0, M]$ for each $M>1$.  Here   $\ell_d (x,t) \equiv 1$   for  $d\geq 3$ and
 $$ \ell_2(x,t) = (\lg t)^2/{\lg(2+ x)} \quad (x< \sqrt t); \quad =  \lg(t/x) \quad ( x\geq \sqrt t).$$ 
  \end{Thm}
\v2

Theorem \ref{thm2.3}  asserts  that  the limit distribution is uniform on the sphere $\partial U(a)$ if $\x/t\to 0$, while it is distributed  with a  positive and continuous density function   $g_{av}(\th_{\xi,\x})$  if $x/t\to v>0$.  
In the case  when $x/t$ is unbounded, the situation becomes different: the weak limit concentrates at $a\x/x$ as one can infer from a result stated next.

 The following theorem, applied in Section 4.2.3, follows from Corollary 2.1 and Lemma 5.6 of \cite{Ucal_b} ((i) is immediate from the former; as for   (ii)  use  the latter in addition).   
\begin{Thm}\label{cor4.1}   Let $d\geq 2$ and $v:=x/t$.

 {\rm (i)}\,     Under the constraint   $ \cos \th_{\xi,\x} > v^{-1/3}$,  uniformly for $\x$, $\xi \in \partial U(a)$ and $t>a^2$,  as $v\to\infty$
$$
\frac{H_{U(a)}(\x, t; d\xi)}{\om_{d-1}a^{2\nu}m_a(d\xi)}  
= \frac{\x\cdot \xi}{t} p^{(d)}_t(|\x-\xi|)\Bigg[1+O\bigg( \frac1{v\cos ^3\th_{\xi,\x}}\bigg)\Bigg],
$$
where $\om_{n}$
 stands for the area of $n$-dimensional unit sphere. 

{\rm (ii)}\, For each  $\e>0$ there exists $M_\e>1$ such that  if $E(\e;\x):=\{\xi \in \partial U(a): \cos \th_{\xi,\x} \leq \e \} $,  $v>M_\e$ and  $t>a^2$, then
$$ H_{U(a)}(\x, t; E(\e;\x)) \leq \k_d \e a^{2\nu}e^{a\e v}p_t^{(d)}(x),$$
where  $\k_d$ is a constant depending only on $d$. 
\end{Thm}

Taking account of the last statement of Theorem 2.1, Theorem 2.5(ii) may be paraphrased as
\beqn\label{xx}
P_\x[ B_{\sigma(U(a))} \in E(\e;\x)\,|\, \sigma_{U(a)} =t]  
\leq \k'_d \e e^{-(1-\e)av} (va)^{\nu-1/2} \quad (v>M_\e, t>a^2).
\eeqn

\section{The  hitting distribution for a bounded Borel set}

  Here we seek   an exact    asymptotic form, as $t\to\infty$,  of   $H_A(\x,t;d\xi)dt$ for  $\x \in  \Om_A$, where $A$ is a bounded Borel set of  $\R^d$. 
  
\subsection{Case $x/t\to 0$} 
We deal with cases $d\geq 3$ and $d=2$ separately. In the case $d\geq 3$ the capacity of $A$ is involved in the leading term, while for  $d=2$, the logarithmic capacity   appears only  in the next order term that we shall not identify.

Our basic strategy is to use the Huygens property of $H_A$  in the form  
\begin{eqnarray} \label{-70}
&&H_A(\x, t; E)  =\int_0^t ds\int_{\partial U(R)} H_{U(R)}(\x, t-s; d\xi)H_A(\xi, s; E) \\
&&\quadd \quadd \quadd \quadd \qquad \mbox{for}\quad  E\subset \bar A, \,\x \notin U(R) \nonumber
\end{eqnarray}
valid if  $A\subset  U(R)$.
 For the Brownian motion started at a  point very distant from  $A$   to hit $A$ for the first time at   time $t$ it must hit $U(R)$  (with $R/x << 1$)  in a relatively small  time  interval just before $t$ with high probability, so that the outer integral of (\ref{-70}) must concentrate on  such an  
interval and  the exact asymptotic forms for balls $U(R)$  described in the preceding section will yield the results for general $A$.  In the case when the starting point is not distant from $A$ the process must make a big  excursion  in the interval $[0,\sigma_A]$ (if $\sigma_A=t$ is very large) and exit 
from  a large sphere within a relatively short time interval so that the problem is reduced to the  case of distant starting points. 
In the actual process of verification there arise  various error terms  that must be properly estimated and 
 to  obtain  certain upper bounds of  $q_A(\x,t)$, being crucial for that purpose, constitutes a substantial part of the proofs.   However,  we  deal with them  in  Section 6  and here we
  focus on the problem of obtaining an exact asymptotic form of $H_A$ by taking  for granted the results verified in Section 6. 

We put 
 $R_A = \sup \{|\y|: \y  \in A^r\}$ as in Section 1 and 
 designate by $c_R, c_R'$ etc. constants depending only on $R_A$ and  $d$ whose exact values are not significant for the present purpose and  may vary  at different occurrences of them.

\v2

{{\bf 3.1.1.} \sc Density of hitting time distribution $(d\geq 3)$.}
Let $d=3,4,\ldots$ and 
 Cap$(A)$ denote  the Newtonian capacity of $A$. For the present purpose it is convenient to  define  it by 
$$\cp(A) =\frac1{G^{(d)}(1)} \lim_{x\to\infty} x^{2\nu}P_{\x}[\sigma_A<\infty].$$ 
  The identity 
 $$ x^{2\nu}P_{\x}[\sigma_A<\infty] = R^{2\nu} E_{\x}\Big[P_{B_{\sigma(U(R))}}[\sigma_A<\infty]\,\Big|\, \sigma_{U(R)} <\infty\Big]$$
   holds true  whenever   $x>R\geq R_A$ and  shows the existence of the limit defining $\cp(A)$ as well as the  formula
\beqn\label{cap1}
\cp(A) = \cp(U(R))P_{m_R}[ \sigma_A<\infty],
\eeqn
where $P_{m_R}$ denotes the law of  $B_t$ started with the uniform probability measure on the sphere $\partial U(R)$. (See also Remark 1 (d) below.) 

\begin{Thm}\label{thm3.1} \,  Let $d\geq 3$.   Uniformly for $\x\in \Om_A$, as $t\to\infty$ and $x/t\to 0$ 
\beqn\label{-52}
q_A(\x,t) = \cp(A)  p_t^{(d)}(x) \Big(P_\x[\sigma_A=\infty](1+o(1))  + {\rm err}(\x,t)\Big),
\eeqn
where   if $x\geq 2R_A$, ${\rm err}(\x,t) =0$ and  if $x < 2R_A$, for any decreasing function $\de(t)$ that tends to zero
\beqn\label{err}
|{\rm err}(\x,t)| \leq CP_\x[\sigma_{A\cup \partial U(2R_A)} > t\de(t)]  \leq C'e^{- \la t\de(t)/R_A^2} 
\eeqn
with  some universal constants $C$ and $\la>0$;  the constant involved  in   $o(1)$  may depend only  on  $d$ and the choice of $\de(t)$ (apart from the dependence on $A$). 
\end{Thm}

 To be precise the term  $o(1)$ in (\ref{-52})  must depend also on  $A$ but the constants involved can be taken independently of  $A$ if   the   limit is  taken   under  $t/R_A\to \infty$ and $R_Ax/t \to 0$ in place of  $t\to \infty$ and $x/t\to 0$. This or similar points  will  not be  stated explicitly in what follows  as in the theorem above. 
\v2
{\sc Remark 1}. (a)\, In view of our definition of Cap$(A)$  and Theorem \ref{thm01} the formula (\ref{-52})  entails   that  as $x\to \infty$ and  $x/t\to 0$
$$q_A(\x,t) =  P_\x[ \sigma_A<\infty] x^{2\nu}q(x,t;1) (1+o(1)),$$
which,  on noting that $1/x^{2\nu} = P_\x[\sigma_{U(1)} <\infty]$,  may be  paraphrased    in terms of conditional probability as follows:  as $x\to \infty$ and $x/t\to 0$
$$P_\x[\sigma_A\in dt \,|\,\sigma_A<\infty] /dt= P_\x[\sigma_{U(1)} \in dt \,|\,\sigma_{U(1)}<\infty] /dt (1+o(1)).$$
This formula does  not hold in general  if $x/t$ is bounded away from zero.   

(b) \,   As  $\x$ gets close to $A$ and $t$ large $q_A(\x,t)$  itself approaches zero very fast, but it is  not so simple a matter  to express  in general how fast it does. In any case,  in (\ref{-52}) the error term  err$(\x,t)$, though very small  for $t$ large  (see (\ref{err})),  cannot be absorbed into  $o(1)$ that precedes it  since for each  $t>1$,  $P_\x[\sigma_{A\cup \partial U(2R_A)} >t]$ may be much larger than  $P_\x[\sigma_A= \infty]$ for some $\x \in \Om_A$. An example is easily constructed by considering a ball  pockmarked by infinitely
many cave-like holes, each one containing a relatively spacial chamber connected by a narrow tunnel to the outside and   among them  there  being one such that   the ratio of diameter of the  tunnel  to   that  of the chamber is smaller than any prescribed number.   

For each $\e>0$,  for reasons of continuity  $\inf_{\x\notin {\rm nbd}_{\e}(A^r)}  P_\x[\sigma_A=\infty]>0$,  and hence  the asymptotic formula  \beqn\label{-55}
q_A(\x,t) = \cp(A)P_\x[\sigma_A=\infty] p_t^{(d)}(x) (1+o(1))
\eeqn
holds  uniformly for $\x \notin$ nbd$_{\e}(A^r)$.
 With this remark one may readily verify that  (\ref{-55})  holds uniformly for $\x\in \Om_A$ if $A$ satisfies some regularity condition such as smoothness of boundary. 

(c) \, For each $\e>0$, we can find  a constant $\eta>0$ such that on  $\{\x\notin {\rm nbd}_{\e}(A^r):  x< \sqrt{2t \lg t} \}$, the factor $o(1)$ may be replaced by $ O(t^{-\eta})$ in (\ref{-52}).  A sketch of  proof is given at the end of this subsection.  

(d) \,  There exists a finite measure $\mu_A$  supported by the closure $\bar A$, called the equilibrium measure,   such that    $P_\x[\sigma_A<\infty]= \int G^{(d)}(|\y-\x|)\mu_A(d\y)$,  $\x\in \R^d$ (see for a concise proof  the arguments given in  \cite[pp.248, 249]{IM}) and  the capacity of $A$ is usually  defined as the total charge of  $\mu_A$. Our definition  conforms to it as is well known: if  $R> R_A$,
$$\cp(A) = \cp(U(R))\int m_R(d\xi) \int G^{(d)}(|\y-\xi|)\mu_A(d\y) = \mu_A(\bar A)$$
where the identities $\cp(U(R))=1/G^{(d)}(R)$ and $ \int G^{(d)}(|\y-\xi|)m_R(d\xi)=G^{(d)}(R)$ ($|\y| < R$) are
used for the  last equality.  (Cf. e.g. \cite[Sections 5.1 and 5.2]{C}; \cite[Proposition (2.5.8)]{B}.)

\v2

A substantial part  of Theorem \ref{thm3.1} is contained in  the next lemma.
\begin{Lem}\label{cor-1}  \,  Let $d\geq 3$. As $x\to\infty$  and  $x/t\to 0$ 
\beqn\label{-5}
q_A(\x,t) =\frac{ \cp(A)}{\cp(U(1))}q(x,t;1) (1+o(1)). 
\eeqn
Here the constant involved  in   $o(1)$  may depend only  on  $d$ and the choice of $\de(t)$
 (apart from the dependence on $R_A$) . 
\end{Lem}

\v2\n
\pf \,  Put  $R= 2R_A$. 
Suppose $x\to\infty$. By strong Markov property 
\beqn\label{3.6a}
q_A(\x,t)  = \int_0^t ds\int_{\partial U(R) } H_{U(R)}(  \x, t-s; d\xi)  q_A(\xi, s).
\eeqn
Take a (large) number  $M$ from the interval $(1, \sqrt t)$ and  split the range of outer integral at $s=M$ and $s=t/2$. We write 
$$I_{[0,M]},\quad  I_{[M,t/2]}   \quad\mbox{and}\quad  I_{[t/2,t]}$$   
for the corresponding integrals over  the Cartesian products  $[0,M]\times \partial U(R)$, etc.
 On employing  Theorem \ref{thm2.3} (with $v=0$),  
\begin{eqnarray}\label{l3.1} \qquad
 I_{[0,M]}
=\int_0^{M} q(x, t-s;R) ds \int_{\partial U(R)}q_A(\xi, s) m_R(d\xi)\bigg(1+ O\Big(\frac{x}{t}\Big)\bigg).
 \end{eqnarray}
 Noting  $p_{ t-s}^{(d)}(x)/p_t^{(d)}(x) = \exp\{- \frac{x^2s}{2t(t-s)}\}(1+o(1))$ ($s<M)$  we apply  Theorem \ref{thm01} to see
 \beqn\label{q/q}
 \frac{q(x,t-s;R)}{q(x,t;R)}=\exp \Big\{- \frac{x^2s}{2t(t-s)} \Big\}(1+ o(1))
 \eeqn
as $s/t\to 0$,  whereas,  since $\sup_{\xi \in \partial U(R)} P_\xi[ M< \sigma_A<\infty] \leq CM^{-\nu} \cp(A) $
 in view of (\ref{J-S}) (cf. also Lemma \ref{lem6.1} or (\ref{-00})  below), the identity   (\ref{cap1}) yields
$$\int_0^{M} ds \int_{\partial U(R)}q_A(\xi, s) m_R(d\xi) =P_{m_R}[\sigma_A<M] = \frac{\cp(A)}{\cp(U(1))}\Big( R^{-2\nu}+ O(M^{-\nu})\Big).$$
 Now suppose  $x/t \to 0$, and make substitution from   these relations in (\ref{l3.1}) and let $M\to \infty$ under  $(1\vee x^2)M/t^2\to 0$.  Then,  observing  
 $$R^{-2\nu}q(x,t;R) = q(x,t;1)(1+o(1))$$
   we find  that 
 $ I_{[0,M]}$ is asymptotic to the right-hand side of (\ref{-5}). 
  
  The integral over $[M,t/2]\times \partial U(R)$  may be  disposed of by means of the inequality
 $$I_{[M,t/2]}  \leq \bigg(\sup_{M<s<t/2}\frac{q(x,t-s;R)}{q(x,t;R)}\sup_{\xi \in \partial U(R)} P_\xi[ M< \sigma_A<\infty]\bigg) q(x,t;R)\bigg( 1+ O\Big(\frac{x}{t}\Big)\bigg).$$
 The first supremum is bounded (owing to Theorem \ref{thm01}), while the  second one is  $\cp(A)\times O(M^{-\nu})$ as noted right after (\ref{q/q}). Thus $I_{[M,\, t/2]}$ is  negligible.
 
It remains to ascertain that  $I_{[t/2,t]}$  is also    negligible.  We prove  that  
  \beqn\label{-02}
  I_{[t/2,t ]} \leq \k_d \cp (A)  q(x,t;R) \times  (R/x)^{2\nu}  e^{- x^2/3t}\qquad (R< x<t/R_A).
  \eeqn
The proof  rests on the bound
\beqn\label{-00}
\sup_{s \geq  t/2} q_A(\y ,s) \leq \k_d\,{\rm Cap(A)} t^{-d/2} 
\eeqn
  valid for all  $t>R^2$ and  $\y \in \Om_A$, where the constant $\k_d$  depends on $d$  only.  The  proof  of this bound  is somewhat involved and   postponed  to the last section (see  Lemma \ref{lem6.3}). (The bound (\ref{-00})  plays a key role also in the proof of Theorem \ref{thm3.1},  but only  for estimation of error terms;  if one is content with not sticking to perfection   the error estimate  of $q_A$ having the heading factor   $R_A^{2\nu}$ in place of $\cp(A)$, he may simply apply  Theorem \ref{thm01} instead of (\ref{-00}).)

Clearly we have
 $$I_{[t/2,t ]} \leq P_{\x}[\sigma_{U(R)}<t/2] \times \sup_{t/2\leq s \leq t}\, \sup_{|\xi|=R} q_A(\xi,s).  $$
  On taking (\ref{-00})  for granted and employing  Theorem \ref{thm2.2}  the right-hand side of this  inequality is dominated by
\[ \k_d \frac{\,{\rm Cap(A)}}{t^{d/2}} \times  \bigg(\frac{R}{x}\bigg)^{2\nu} \int_{x^2/t}^\infty e^{-u}u^{\nu-1}du 
 \leq \k'_d \frac{\,{\rm Cap(A)}}{t^{d/2}} \times  \bigg(\frac{R}{x}\bigg)^{2\nu}  \Big([e^{-x^2/t}(x^2/t)^{\nu-1}] \wedge 1 \Big),
 \]
hence by the right-hand side of   (\ref{-02}) as is  ensured by a crude estimation  using  the asymptotic form  of  $q(x,t;R)$ given in Theorem \ref{thm01}.  The proof of Lemma \ref{cor-1} is complete.  \qed  
\v2

 In comparison with the one-dimensional result  we have 
 \beqn\label{spect}
 \sup_{\y\in U(r)}P_\y[\sigma_{ \partial U(r)} >t]  \leq C e^{-\la t/r^2}
 \eeqn
with some universal constants $C$ and  $\la>0$.

 \begin{Lem}\label{crucial_d} 
Let $r >2 R_A$,  $\x\in U(2R_A)\cap \Om_A$ and $T>r^2$. Then,  with  universal constants $C>0$ and $\la>0$,
$$
P_\x[\sigma_{A\cup \partial U(r)} >T] \leq Ce^{-\la T/ 2r^2} \Big(P_\x[\sigma_A=\infty] + P_\x[\sigma_{A\cup \partial U(2R_A)} >\tst12 T]\Big).
$$
 \end{Lem} 
 \n
 \pf\,   Let $R=2R_A$. 
We break the event $\sigma_{A\cup \partial U(r)} >T$ according as 
 $\sigma_{ \partial U(R)}$ is less than $T/2$ or not  and  infer that  
   \beq
&&P_\x[\sigma_{A\cup \partial U(r)} >T]  \\
 && = \iint_{[0,\frac12 T]\times \partial U(R)}P_\xi[ \sigma_{A\cup \partial U(r)} > T -s ]P_\x[\sigma_{\partial U(R)}<\sigma_A,  \sigma_{\partial U(R)}\in ds, B_{\sigma_{\partial U(R)}}\in d\xi]\\
&& \,\,+\,  \int_{U(R)\cap \Om_A}P_\y[ \sigma_{A\cup \partial U(r)} > \tst12 T ]P_\x[B_{\tst12 T}\in d\y, \sigma_{A\cup\partial U(R)}>\tst12 T]\\
  &&  \leq\Big( P_\x[\sigma_{\partial U(R)} < \sigma_A] + P_\x[\sigma_{A\cup \partial U(R)} >\tst 12 T] \Big)\sup_{\y\in U(R)}P_\y[\sigma_{\partial U(r)} >\tst12 T].
 \eeq
Thus  the assertion of the lemma  follows from (\ref{spect}).  \qed

 \v2
\n
{\it Proof of Theorem \ref{thm3.1}.}  \, In view of Lemma \ref{cor-1} (together with $\cp(U(1))=\La_\nu(0)$),  Theorem \ref{thm01} and the fact that  $\lim_{x\to\infty}P_\x[\sigma_A=\infty]=1$
we may suppose that   $x$ remains in a bounded set as $t\to\infty$.   Given $T<t$ and $r>x$, we decompose
 \begin{eqnarray}\label{-51}
 q_A(\x,t)  = \int_0^{T} \int_{\partial U(r) } P_\x[\sigma_{\partial U(r)} \in ds, B_{\sigma(\partial U(r))}\in d\xi, \sigma_A>s]  q_A(  \xi, t-s)  + \e(\x,t),
 \end{eqnarray}
 where  $\e(\x,t) ={\textstyle P_\x[\sigma_{\partial U(r)} > T, \sigma_A\in dt]/dt}$.  Both $T$ and $r$ may depend on $t$;  we  take  $T=2 t\de(t)$ with  $\de(t) > 1/(1\vee \lg t)$  as well as     $\lim_{t\to\infty} \de(t) = 0$; also,  e.g. $r(t) = t^{1/3}$,  
 so that 
 $T/ r^2> t^{1/4},  $ entailing   
 \beqn\label{T/r}
 p^{(d)}_{t-s}(\xi) \sim  p^{(d)}_{t}(x) \sim p^{(d)}_{t}(0) \quad \mbox{ uniformly for}  \quad s<T,  \,\, \xi\in \partial U(r).
 \eeqn

 The term  $\e(\x,t)$ contributes  only to the error terms. Indeed,
  \begin{eqnarray}\label{e0}  
\e(\x,t)   & =&
  \int_{ U(r)} {\textstyle P_\x \Big[ B_{T}\in d\y, \sigma_{A\cup \partial U(r)} > T \Big]q_A(\y , t-T)}  \nonumber \\
  &\leq& \k_d \cp(A) t^{-d/2}P_\x[ \sigma_{A\cup \partial U(r)} > T], 
   \end{eqnarray}
where we have applied  (\ref{-00}) for the inequality.  Hence,  $\e(\x,t)$  is  absorbed into the error terms 
represented by $o(1)$ or  err$(\x,t)$ in (\ref{-52})  if  $x\leq R$ owing to Lemma \ref{crucial_d}, while
it is negligible in the case $x>2R_A$ when $P_\x[\sigma_A=\infty]\geq 1-2^{-2\nu}$, for  by (\ref{spect}) 
 $
P_\x[ \sigma_{A\cup \partial U(r)} > T] \leq P_\x[ \sigma_{\partial U(r)} >  T] \leq C e^{-\la T/r^{2}}\leq Ce^{-\la t^{1/4}}$.

The first  term on the right-hand side of (\ref{-51})   contains the principal  part. We need some care to obtain the exact asymptotic form.
Since $r\to \infty$,  Lemma \ref{cor-1} applies to $q_A(\xi,t-s)$ and  noting (\ref{T/r}) we  thus obtain   that uniformly for $\xi\in \partial U(r)$ and $s\leq T$, 
 \beqn\label{R200}
q_A(  \xi, t-s) =\cp(A)p^{(d)}_t(x) (1+o(1)).
\eeqn
On the other hand, defining    $\eta(\x,t, T)$ (for $T>1$)  via 
\begin{eqnarray}\label{R201}
\int_0^{T} P_\x[\sigma_{\partial U(r)} \in ds, \sigma_A> s]  &=& P_\x[\sigma_{\partial U(r)} <T \wedge \sigma_A] \nonumber\\
& =& P_\x[\sigma_A=\infty] - \eta(\x,t, T),
\end{eqnarray}
we observe that   
$\eta(\x,t, T) = P_\x[ \sigma_{\partial U(r)} \geq T, \sigma_A=\infty]-  P_\x[\sigma_{\partial U(r)} <T \wedge \sigma_A, \sigma_A<\infty]$,  hence 
 \beqn\label{-8}
-  P_\x[  \sigma_{\partial U(r)} <\sigma_A <\infty]    \leq \eta(\x,t, T)  \leq P_\x[\sigma_{\partial U(r)} \geq T, \sigma_A=\infty]. 
 \eeqn
 Note that $P_\xi[\sigma_A<\infty] \leq (R_A/r)^{2\nu} < 1/2$ for $\xi\in \partial U(r), r >R$  and use  the strong Markov property to  deduce first $P_x[\sigma_A=\infty] \geq  \frac12 P_\x[\sigma_{\partial U(r)} <\sigma_A] $ and then 
  \beqn\label{R202}
P_\x[\sigma_{\partial U(r)} <\sigma_A < \infty]  \leq  2(R_A/r)^{2\nu}P_\x[\sigma_A=\infty],
   \eeqn
provided  $x\vee R <r$. 
 Applying Lemma \ref{crucial_d}  to the right-most  member  in (\ref{-8})  therefore yields
   $$|\eta(x,t)|  \leq C\Big[ (R_A/r)^{2\nu}P_\x[ \sigma_A=\infty] + P_\x[ \sigma_{A\cup \partial U(R)}  >\tst12 T]\Big] \quad \mbox{for}\quad  x<R, 
 $$ 
 while  for  $x\geq R$, $\eta(\x,t)$   plainly makes only a negligible contribution in view of (\ref{spect}).
   Putting  (\ref{R200}) (\ref{R201}), (\ref{R202})  and this bound of $\eta$  together  we  find that the repeated integral in (\ref{-51}) agrees with the asserted asymptotic  formula of the theorem.  By what is remarked on $\e(\x,t)$ right after  (\ref{e0}) this completes the proof. \qed

\v2\n
{\it Sketch of the proof of  Remark 1 (c)}.  Let  $d\geq 3$. Theorem 3 of \cite{Ubes} provides an error estimate for the asymptotic form of  
 $q(x,t;a)$ given in Theorem \ref{thm01} and  according to  it   the error term $o(1)$ in (\ref{R2}) can be replaced by $O(t^{-1}\lg t)$ for $x<\sqrt {2t\lg t}$   (the asymptotic form being exact if $d=3$).  On examining the proof of Lemma \ref{cor-1}  this shows that 
under the same constraint on $x$, (\ref{-5}) can be refined to 
 \beqn\label{refine}
  q_A(\x,t) = \cp (A)p^{(d)}_t(x) (1 + O(x^{-2\nu}\vee t^{-\eta}))   \quad\,\,  (x<\sqrt {2t\lg t}) 
  \eeqn
 hence  $o(1)$  in   (\ref{R200}) to  $O(r^{-2\nu}\vee t^{-\eta})$,  for some $\eta>0$.
The rest is done by noting that $P_\x[\sigma_A=\infty] = 1+ O(x^{-2\nu})$, the argument in the proof of Theorem \ref{thm3.1} is valid also on  e.g.  $\{\x\in \Om_A: R< x< t^{1/5}, \x\notin {\rm nbd}_{\e}(A^r)\}$, and   the terms  appearing in the form $e^{-\la T/r^2}$  therein are all negligible. 
     
 \v2
 {{\bf 3.1.2.} \sc Density of hitting time distribution $(d= 2)$.} \,
Let $A$ be a bounded Borel set of $\R^2$ that is non-polar. Define 
\beqn\label{3.1.2}
e_A(\x) = \pi \lim_{|\y|\to\infty} g_{\Om_A}(\x,\y),
 \eeqn
 where $g_{\Om_A}(\x,\y)$  denotes the Green function for $\Om_A$, which is continuous in the interior of  $\Om_A\times \Om_A$.  
  If $A$ contains a curve connecting the origin with $\partial U(R_A)$,  then   $ 0\leq e_A(\x)\leq C$ if  $|\x|\leq R_A$,  as   being assured  by the scaling property of $e_A$; if $A$ does not, $\sup_{\x\in U(R_A)} e_A(\x)$ could be arbitrarily large (depending on  $A$). It may be well known that $e_A(\x) \sim \lg x$ as $x\to \infty$; in fact, we can show  (cf. (\ref{00*}))  that for $x\geq R_A$,
 \beqn\label{e2}
 0\leq e_A(x)- \lg (x/R_A) \leq C\b_A, \quad \b_A:= m_{2R_A}(e_A),  
 \eeqn
  where $C$ is a universal constant and  $m_R(e_A) = \int_{\partial U(R)} e_A dm_R$.
(Cf.   Section 6.1 for  further properties of $e_A$.)
 
 \begin{Thm}\label{thm3.2}  Let $d=2$ and suppose $A$ to be  non-polar. Then,  as $t\to\infty$  and $x/t\to 0$ 
\begin{eqnarray}\label{-53} 
q_A(\x,t) = p^{(2)}_t(x) 
\times
\left\{ \begin{array} {ll}  
 {\displaystyle  \frac{ 4\pi }{(\lg t)^2}\bigg[e_A(\x)\bigg(1+O\Big(\frac1{\lg t}\Big)\bigg) +{\rm err}(\x,t)\bigg] }\quad& (x  \leq \sqrt t\,),\\ [5mm]
 {\displaystyle  \frac{\pi } {\lg (t/x) } \bigg( 1+ O\Big(\frac1{\lg (t/x)}\Big)\bigg)}  \quad&(  x   > \sqrt t\,),
  \end{array} \right.   
\end{eqnarray}
where   if $x\geq 2R_A$, ${\rm err}(\x,t) =0$ and  if $x < 2R_A$, 
\beqn\label{err2}
|{\rm err}(\x,t)| \leq CP_\x[\sigma_{A\cup \partial U(2R_A)} > t/\lg t]  \leq C'e^{- \la t/R_A^2\lg t} 
\eeqn
with  some universal constants $C$ and $\la>0$. 
\end{Thm}

\v2
{\sc Remark 2.}  \,(a)  The   $O$ terms in (\ref{-53}) depend on $A$. If   we write
them  in the form 
$\b_A \times O(\cdot)$, then the  $O(\cdot)$'s in  this form  are independent of $A$ (except via scaling of $t$ by $R^2_A$  and $x$ by  $R_A$).
 One may realize this by observing that   the dependence on  $A$ other than  $R_A$ comes in via  (\ref{sigma-0}) given below.
It is noted that $ \sup_{\xi\in \partial U(2R_A)} e_A(\xi) < C\b_A$  owing to  Harnack's  inequality, and that
since $-\lg R + m_R(e_A)$ is independent of $R\geq R_A$, $\b_A= \lg 2 + m_{R_A}(e_A)$.  (Cf. (\ref{lim_thm}) and (\ref{*0*})).  

(b) $t/\lg t$ in (\ref{err2}) may be replaced by $\de(t)t$ with $\de(t)\to 0$ as in (\ref{err}).

\v2\n
\pf\,  Let $R=2R_A$. 
The case $x\to \infty$    is dealt with in almost  the same way as in  the proof of Lemma \ref{cor-1} (apart from  $O(\cdot)$ terms in (\ref{-53})). It is  only pointed out that for the estimation of  $I_{[0, M]}$ and $I_{[M,t/2]}$,  we apply
\beqn\label{A_M}
P_\x[\sigma_A > t]  \leq \frac{Ce_A(\x)}{\lg (t/R_A^2)}  \quad\mbox{ for}\quad  x\geq 2R_A \,\,\mbox{and}  \,\, t>R_A^2
\eeqn
(with $\x$ and $ t$ replaced by $\xi \in \partial U(2R_A)$ and  $M$, respectively) and, instead of the bound  (\ref{-02}),  we have
 \beqn\label{-03}
 I_{[t/2,t]}   \leq  \frac{C\b_Aq(x,t;R)}{\lg (x/R)}   e^{- x^2/3t}\quad \mbox{for}\quad R< x<t.
  \eeqn  
These relations  are deduced from Proposition \ref{prop6.4}  of Section 6.3: the former one immediately; while the latter as in the proof of Lemma \ref{cor-1} by observing    
$$q_A(\y, s) \leq C\b_A q(x, t; R_A)e^{x^2/2t}/\lg(x/R)\qquad (s>t/2,\y\in \partial U(R))$$
with the help of  Theorems \ref{thm01}  (a better bound is provided by   Lemma \ref{lem6.10} whose proof  however requires  an additional work). 

  In order to obtain  the error terms $O(\cdot)$ in (\ref{-53}) 
  we need to use the following refinement of (\ref{R21}): uniformly for $x>a$, as $t\to\infty$ 
 \begin{eqnarray}\label{imp} 
  q(x,t;a) = \left\{ \begin{array} {ll} {\displaystyle  \frac{4\pi \lg( x/a)}{(\lg t)^2}\,p_t^{(2)}(x)\Big(1+ O\Big(\frac1{\lg t}\Big)\Big)\quad} \,\,&\mbox{if}\quad  x<\sqrt{t}, \\[4mm]
{\displaystyle  \frac{\pi}{\lg (t/x)}\, p_t^{(2)}(x)\Big(1+ O\Big(\frac1{\lg (t/x)}\Big)\Big)  }
  \quad & \mbox{if}\quad \sqrt{t} \leq x <t, 
  \end{array}\right.
 \end{eqnarray}
  which is a part of  Corollary 4 of \cite{Ubes}.   This allows us to replace $o(1)$ by  \,  $O\Big(1/{\lg (t/x)}\Big)$ \,   in   (\ref{q/q}) if $s/t\to 0$, so that   taking $M = t/(x\vee \sqrt t) =\frac{t}{x} \wedge \sqrt t$ in the proof of Lemma \ref{cor-1} and employing  (\ref{A_M}) and (\ref{-03})  we deduce 
  \beqn\label{imp1}
  q_A(\x,t) = q(x,t; R)\Big(1+ O\Big(\frac1{\lg (t/x)}\vee \frac1{\lg x}\Big)\Big)\quad \mbox{for}\quad  2\vee R< x<t/2,
\eeqn
  which implies the second half of (\ref{-53}).  Since $e_A(\x) = {\lg x} + O(1)$ (cf.  (\ref{e2})), it also  disposes of
the case $t^{1/4} \leq x <\sqrt{ t}$ of  (\ref{-53})   as well.
  
  Finally  consider  the case $x< t^{1/4}$. A crucial issue arising in this case is dealt with    by using a result from Section 6.1  (in addition to   (\ref{A_M})),  which we state  as  Lemma \ref{lem3.01}  after the present proof.  
 Taking  $r=  t^{1/3}$ and $T=t^{3/4}$ 
 we make use  of 
  the same expression  as (\ref{-51}):
  \begin{eqnarray}\label{-51_2}
 q_A(\x,t)  = \int_0^{T} \int_{\partial U(r) } P_\x[\sigma_{\partial U(r)} \in ds, B_{\sigma(\partial U(r))}\in d\xi, \sigma_A> \sigma_{\partial U(r)}]  q_A(  \xi, t-s)  + \e(\x,t),
 \end{eqnarray}
 where  $\e(\x,t) = P_\x[\sigma_{\partial U(r)} > T, \sigma_A\in dt]/dt$.
 Observe that  for $s<T$,
  $q_A(\xi, t-s)$ $(\xi\in \partial U(r))$  may be replaced by
 \beqn\label{ooo}
 \frac{4\pi \lg r}{(\lg t)^2} p^{(2)}_t(x)\Big(1+ O\Big(\frac1{\lg t}\Big)\Big)
 \eeqn
  in view of (\ref{imp}) and  (\ref{imp1}) and   the repeated integral  on the right-hand side of (\ref{ooo}) is accordingly   factored into the product of this quantity and $\int_0^TP_\x[\sigma_{\partial U(r)} \in ds, \sigma_A>s]$,   of which  the latter is reduced to $P_\x[\sigma_{\partial U(r)} < T\wedge \sigma_A]$.
We make the decomposition
\[
P_\x[\sigma_{\partial U(r)} < T\wedge \sigma_A] = P_\x[\sigma_{\partial U(r)} <  \sigma_A] -   
P_\x[T\leq \sigma_{\partial U(r)} <  \sigma_A].
\]
Multiply the first probability on the right   by the quantity   (\ref{ooo})  and applying  Lemma \ref{lem3.01} below, we  find  the first formula on the right-hand side of (\ref{-53}) with err$(\x,t)$ discarded.   
  On the other hand,  analogously to Lemma \ref{crucial_d} we have for $\x\in U(R)\cap\Om_A$
 \beqn\label{crucial}
P_\x[\sigma_{A\cup \partial U(r)} \geq T] \leq Ce^{-\la T/ 2r^2} \Big(e_A(\x) + P_\x[\sigma_{A\cup \partial U(R)} \geq \tst12 T]\Big),
\eeqn
which together with (\ref{A_M})  shows that the second one, namely  
$P_\x[T\leq \sigma_{\partial U(r)} <\sigma_A]$,  contributes only to the  error terms in (\ref{-53}), and accordingly  that  the repeated  integral in (\ref{-51_2})  agrees with  the first  formula on the right-hand side  of   (\ref{-53}).  
The proof of (\ref{crucial}) is omitted, being the same  as that of  Lemma \ref{crucial_d} except for the use  of Lemma \ref{lem3.01}.

Plainly $\e(\x,t)$ is dominated by 
$$  P_\x[\sigma_{A\cup\partial U(r)} \geq T]\sup_{\y\in U(r) \cap\Om_A}q_A(\y, t-T)$$
and  is    absorbed into the error terms  as readily shown by using  (\ref{crucial}) as well as  the fact that  the supremum above is bounded by a universal constant (see  Lemma \ref{lem0} of Section 6.3).  Proof of Theorem \ref{thm3.2} is complete.   \qed

 \begin{Lem}\label{lem3.01}  \, 
 Let $d =2$.  For  $r\geq 2R_A$ and  $\x\in \Om_A\cap U(r)$,  
\beqn\label{sigma-0}
e_A(\x) \bigg( 1 -\frac{3m_{R_A}(e_A)}{\lg (r/R_A)}\bigg)\leq  P_\x[\sigma_{\partial U(r)} < \sigma_A] \lg\Big(\frac{r}{R_A}\Big) \leq e_A(\x).  
\eeqn
 \end{Lem}
 \v2
Proof of this lemma  is given in  Section 6.1 (see Proposition \ref{prop6.1}).

\v2

{{\bf 3.1.3.} \sc Asymptotic form of $H_A$ as $x/t\to 0$.} \,
Remember that for a bounded Borel set $A$,  $H_A(\x,t; d\xi)$ is defined in (\ref{H_A}) and $H^\infty_A(d\xi)$   in (iii) of Section 1.  
 \begin{Thm}\label{thm3.3}  \, Let $d\geq 3$.  Uniformly for  Borel subsets $E$ of $\partial A$,   as $t\to\infty$ and  $x/t\to 0$ 
  \beqn\label{-6}
H_A(\x, t;  E)= q_A(\x,t) \Big[H^\infty_A(E)+ \e_{x,t}(E) \Big]  + \cp(A) \, {\rm err}_{\x,t}(E)
\eeqn
where  $ \e_{x,t} $ and  ${\rm err}_{\x,t}$ are  signed-measures on $\partial A$ such that $|\e_{\x,t}|_{{\rm t.var}}= o(1)$ and that 
if $x\geq 2R_A$, ${\rm err}_{\x,t} \equiv 0$ and  if $x < 2R_A$, $|\rm{ err}_{\x,t}|_{{\rm t.var}}$ admits  the same bound as err$(\x,t)$ in Theorem \ref{thm3.1}, namely  $|{\rm err}_{\x,t}|_{{\rm t.var}} \leq CP_\x[\sigma_{A\cup \partial U(2R_A)} > t\de(t)]$, with  some universal constant $C$,  for any decreasing function $\de(t)$ that tends to zero. 
\end{Thm}

\v2\n
\pf\,   Suppose   $x \to \infty $ and $x/t \to 0$.  We proceed as in  the  proof of Lemma \ref{cor-1} with    (\ref{3.6a}) replaced by the Huygens relation (\ref{-70}). As before,  let   $R= 2 R_A$ 
and  restrict  the outer integral of (\ref{-70})  to $[0,M]$ with   $M$ to be made indefinitely large under $(1\vee x^2)M/ t^2 \to 0$, $M\leq \sqrt t$.
 Observe  $\cp(A) = \La_\nu(0)R^{2\nu}P_{m_R}[\sigma_A<\infty]$ (see   (\ref{cap1}))  and plug it in the formula of  Theorem \ref{thm3.1}.   Then using  Theorem \ref{thm01} we  find 
\beqn\label{-80}
q_A(\x,t) = P_\x[\sigma_A=\infty] P_{m_R}[\sigma_A<\infty]q(x,t;R)(1+o(1)),
\eeqn
in which  the factor $ P_\x[\sigma_A=\infty]$, tending to unity,  may be dropped from the right-hand side. 
Now  apply   Theorem \ref{thm2.3}, the relation $ q(x,t-s;R) = q(x,t;R)(1+ o(1))$ valid for $s<M$  and (\ref{-80}) in turn.
The   integral in (\ref{-70})   restricted on $[0,M]$ then   may be  written  as
\begin{eqnarray}\label{-9}
&&\int_0^M q(x,t-s;R)ds\int_{\partial U(R)}H_A(\xi, s; E) m_R(d\xi)(1+ O(x/t)) \\
 &&= q(x,t;R)P_{m_R}[ B_{\sigma(A)}\in E,\sigma_A<M] (1+o(1)) \nonumber\\
 &&= q_A(\x, t)P_{m_R}[B_{\sigma(A)} \in E\,|\,\sigma_A<\infty](1+ o(1)), \nonumber
\end{eqnarray}
of which  the last expression may be  further  rewritten as  the right-hand side of (\ref{-6}) since $P_{m_R}[B_{\sigma(A)} \in E\,|\,\sigma_A<\infty]= H^\infty_A(E)$ (that follows from $ P_{m_R}[B_{\sigma(A)} \in E, \sigma_A<\infty] = \lim_{x\to\infty}  P_{\x}[B_{\sigma(A)} \in E, \sigma_A<\infty \,|\,\sigma_{U(R)}<\infty]$). For the remaining integral we may consider only the case  $E=\partial A$ and the required bound of it  obviously  follows from those obtained in the proof of   Lemma  \ref{cor-1}.

The case when  $x$ remains in a bounded set can be dealt with in a similar way by tracing the corresponding part of the proof of Theorem \ref{thm3.1}. Here we only mention a way of how the signed measure err$_{\x,t}(d\xi)$ is determined. In the proof of Theorem \ref{thm3.1} the term err$(\x,t)$  comes from $\e(\x,t)$ and $\eta(\x,t)$ and estimated by using Lemma \ref{crucial_d}. These are adapted for construction as well as  estimation of   err$_{\x,t}(d\xi)$ as follows. The decomposition     (\ref{3.6a}) is replaced by  
\beq
H_A(\x, t; E) &=&\int_0^{T} \int_{\partial U(r)} P_\x[\sigma_{\partial U(r)}\in ds, B_{\sigma(\partial U(r))}\in d\xi,\sigma_A> \sigma_{\partial U(r)}]H_A(\xi, t-s; E) \\
&&+\,\, \e(\x,t;E).
\eeq
with $\e(\x,t;E) = P_x[ \sigma_{\partial U(r)}> T, \sigma_A\in dt, B_{\sigma(A)} \in E]/dt$; for $x<R$,  the part of  err$_{\x,t}(E)$ that comes from $\e(\x,t;E)$ is accordingly given by
$$\frac1{\cp(A)} \int_{U(r)} P_\x[\sigma_{A\cup \partial U(R)} >\tst12 T, \sigma_{A\cup \partial U(r)} > T, B_{T}\in d\xi,\sigma_A> \sigma_{\partial U(r)}]H_A(\xi, t-T; E).$$
We may analogously define $\eta(\x,t;E)$  and  determine the corresponding part for it. The estimation of $|\rm{ err}_{\x,t}|_{{\rm t.var}}$ is trivially made  by using  Lemma \ref{crucial_d}.
  The details are omitted. \qed 
 \v2
 {\sc Remark 3.} \,   Under the constraint $x<\sqrt{2t\lg t}$  the $o(1)$ term  in (\ref{-6})   can be replaced by $O(t^{-\eta})$ with some constant  $\eta>0$.  Verification is  the same as that of Remark 1 (c).
 
 \v2
Here we record the formula: for $d\geq 3$,
\beqn\label{H/H}
\int_0^\infty dt\int_{\partial U(R)} H_A(\xi,t; E)m_R(d\xi) = G^{(d)}(R)\cp(A) H^\infty_A(E),  \quad R\geq R_A.
\eeqn
The proof may be performed by  recalling our definition of $\cp(A)$ to see that $H_A^\infty(E)$ equals  $ \lim_{x\to\infty} x^{2\nu} P_\x[B_{\sigma(A)} \in E, \sigma_A<\infty]/[G^{(d)}(R)\cp(A)]$, and 
  applying the Huygens property to the probability under the limit sign. 
 
 \begin{Thm}\label{thm3.40}  \, Let $d =2$.  Uniformly for  Borel subsets $E$ of $\partial A$,   as $t\to\infty$ and  $x/t\to 0$ 
 \beqn\label{eq3.4}
H_A(\x, t;  E)= q_A(\x,t)\Big[ H^\infty_A(E) + \e_{\x,t}(E)\Big] + {\rm err}_{\x,t}(E),
\eeqn
where  $ \e_{\x,t} $ and  ${\rm err}_{\x,t}$ are  signed-measures on $\partial A$ such that the total variation 
$ |\e_{\x,t}|_{{\rm t. var}} \leq  C\b_A/ \lg [t/(1\vee x)] ;$ and that 
if $x\geq 2R_A$, ${\rm err}_{\x,t} \equiv 0$ and  if $x < 2R_A$, $|\rm{ err}_{\x,t}|_{{\rm t.var}}$ admits  the same bound as err$(\x,t)$ in Theorem \ref{thm3.2}, namely  $|{\rm err}_{\x,t}|_{{\rm t.var}} \leq CP_\x[\sigma_{A\cup \partial U(2R_A)} > t /\lg t]$ ($C$ is   a universal constant in both places). 
\end{Thm}

\v2\n
\pf\,  The  proof is substantially the same as that of Theorem \ref{thm3.3}  apart from the  estimate $ |\e_{\x,t}|_{{\rm t. var}} \leq  C\b_A/ \lg [t/(1\vee x)]$.  We consider only the case  $ R< x <  t^{1/2}$ ($R=2R_A$), the other case being similar.  For  $t^{1/4}< x <  t^{1/2}$, we apply Theorem \ref{thm2.3} and  (\ref{imp})   to see  
$ q_A(\x,t) = 
q(x,t; R)(1+ O(1/\lg t))$
(in place of (\ref{-80})) and observe that  both $O(x/t)$ and $o(1)$ in (\ref{-9}) can be replaced  by $O(1/\lg t)$, which suffices for the required error estimate,  $O(\cdot)$ term being  absorbed into $\e_{\x,t}(E)$.
  For  $x\leq t^{1/4}$, the argument given in the proof of Theorem \ref{thm3.2} leads to the  error estimate
  asserted in the theorem. \qed

\subsection{Case $x/t\to\infty$} 
\v2
In this subsection we suppose  $K$ to be a compact set. Let $\be$ be a unit vector and ${\rm pr}_{\be}A$  the orthogonal projection of a set  $A$ on $\De_{\be}$, the hyper-plane   perpendicular to $\be$  passing through the origin. We often write $\pr_\x$ for  $\pr_\be$ if $\be = \x/x$. We  bring in   the subset $\lan K\ran_{\be} $ of $\partial K$  given by
\beqn\label{eq_h}
\lan K\ran_{\be} = \{\xi \in \partial K : \xi+t\be \notin K \,\,\mbox{ for}  \,\,\, t>0\},
\eeqn
the mapping  $h= h_{\be,K} : \De_{\be}   \mapsto \R$ 
by the relation
$$
 \z+ h(\z)\be \in \lan K\ran_{\be} \quad \mbox{if} \quad \z\in {\rm pr}_{\be}K;\quad h(\z) =-\infty\quad\mbox{otherwise}
$$
 and   the  Borel measure  $m_{K,\be}$ on $\partial K$ by
\beqn\label{mm}
m_{K,\be}(E) = {\rm vol}_{d-1}\Big(\pr_{\be}(\lan K\ran_{\be} \cap E)\Big)   \quad\quad\quad (E\subset \partial K\,).
\eeqn
 Here we regard  ${\rm pr}_\be K$  $(={\rm pr}_{\be} \lan K\ran_{\be}$) as a subset of $(d-1)$-dimensional Euclidean space.  Denote by $\partial_{d-1}({\rm pr}_\be K)$ the boundary of it as such. Note that $\lan K\ran_{\be}$ is a Borel subset of  $\partial K$ and $m_{K,\be}$ is supported by  $\lan K\ran_{\be}$. We  need to further bring in the set of discontinuity points of $h$ given by 
\beqn\label{mmm}
\mbox{dis-ct}_{\be}(K)
= 
\{\z\in {\rm int}({\rm pr}_{\be} K): h \mbox{ is discontinuous at } \z\} \cup \partial_{d-1}( \pr_\be K). 
\eeqn

\begin{Thm}\label{thm3.4}   Let  $\be\in  \partial U(1)$ and  suppose that {\rm dis-ct}$_{\be}(K)$ is of zero  Jordan measure, 
namely 
{\rm
\beqn\label{cond1}
{\rm vol}_{d-1}\Big(\overline{\mbox{dis-ct}_{\be}(K)}\, \Big)=0.
\eeqn 
}
Then, as $t\to\infty$ with $v:= x/t\to \infty$    
  \beqn\label{eq3.5}
H_K(x\be, t;  d\xi) \approx  vp^{(d)}_t(x)e^{v\be\cdot\xi} m_{K,\be}(d\xi).  
\eeqn
\end{Thm}
\v2

By (\ref{eq3.5}) we mean that  the Borel measure 
\beqn\label{300}
\mu_{t, \x}^K(d\xi):= \frac{e^{-v\be\cdot\xi} H_K(x\be, t;  d\xi)}{vp^{(d)}_t(x)}
\eeqn
 ($d\xi\subset \partial K$, $\x = x\be$)  asymptotically concentrates on $\lan K\ran_\be$ so that  
 $\lim \mu_{t,\x}^K(\partial K \setminus \lan K\ran_\be) =0$
and converges to $m_{K,\be}(d\xi)$ on $\lan K\ran_\be$  in the sense that (as $x/t\to\infty$)
 \beqn\label{K6}
 \mu_{t, \x}^K(E) \to m_{K,\be}(E)  \quad \mbox{for each Borel set}\,\, E\subset \lan K\ran_\be\,\,  \mbox{with}\,\,  {\rm vol}_{d-1}(\partial_{d-1}({\rm pr}_\be E))=0.
 \eeqn
 With  the condition (\ref{cond1}) being assumed we can readily verify that the relation (\ref{eq3.5}) implies the weak convergence of $\mu_{t, \x}^K$ as stated in (vii) (the coverse is also  true).

    If  $a=\max \{\xi\cdot \be: \xi\in K^r \} $ and  $m_{K,\be}(\{\xi\in \lan K\ran_\be: \xi\cdot \be >  a-\e\})>0$ for any $\e>0$, then it follows from (\ref{eq3.5}) that the hitting site distribution  $H_K(x\be,t; \cdot)/q_A(x,t)$ tends to concentrate  on $\{\xi\in K^r  : \xi\cdot \be =a\}$ as $x/t \to\infty$.    The formula is also  useful for a study of Wiener sausage 
 (see the end of Section 5)  as well as the hitting site distribution for the Brownian motion with large drift, while 
 for the part $\partial K \setminus \lan K\ran_{\be}$,  (\ref{eq3.5})  provides only a crude upper bound of $H_K(x\be;t,d\xi)$;  there is no direct way to derive from (\ref{eq3.5})   any exact asymptotic form of $q_A(\x,t)$. 

Formula  (\ref{eq3.5}) may be formally  inferred by looking at the space-time distribution of  the first arrival of $(B_t)$  on the plane passing through $\xi$ and perpendicular to $\be$. Indeed,  if $B_0=x\be$, its density is given by 
\beqn\label{hd_pln}
\frac{x-\xi\cdot \be}{t}p^{(1)}_t(x- \xi\cdot \be)p^{(d-1)}_t(|{\rm pr}_{\be}\xi|),
\eeqn
 which  is asymptotic to  the right-hand side of (\ref{eq3.5}) divided by  $m_{K,\be}(d\xi)$. The actual  proof  is postponed to the last part of Section 4.2,  since we need to use a result given there.

\v2
{\sc Remark 4.}   
(a)\, There is a large class of  compact sets  $K$ such that   ${\rm vol}_{d-1}(\pr_\be K)= 0$ for all  $\be$ and
 $K$ is  non-polar. This condition obtains  if $K$ has Hausdorff dimension larger than  $d-2$ and zero    Hausdorff measure of $(d-1)$-dimension (\cite{Fl} Theorem 6.4 ($d\geq 3$), \cite{B} Exercise 27 of p.373 ($d = 2$)). 
 
 (b)\, It is not clear whether   the condition (\ref{cond1}) needs to be   required  for (\ref{eq3.5}).     If $\lan  K\ran_{\be}$ is contained in the union of a finite number of  hyper-planes perpendicular to $\be$, it may be removed. For the upper bound the   condition (\ref{cond1}) is not needed.

\v2
 Here we state  and prove a  lower bound of $H_A(\x,t;E)$.    

\begin{Lem}\label{lem3.5}  There exists a  constant $\k_d>0$ depending only on  $d$  such that for $0<\eta<t/2$,  $x> R:=2R_A$, $t > R_A^2+\eta$,  and $E \subset \partial A$,  
$$H_A(\x,t;E) \geq  \k_d q(x,t-\eta; R)\inf_{\xi\in \partial U(R)} P_{\xi}[B_{\sigma(A)}\in E, \sigma_A \leq \eta].$$
\end{Lem}
\n
\pf\, First consider the case  $x>t$. In the Huygens representation (\ref{-70}) of $H_A(\x,t;E)$  we restrict the range of integration to $[0,\eta]\times \{\xi\in \partial U(R):  \x\cdot \xi> 0\}$ and apply Theorem \ref{cor4.1}   to see
\[
H_A(\x,t;E) 
\geq R^{2\nu} \int_0^\eta ds\int_{\xi\cdot \x > 0}\frac{\xi\cdot\x}{t-s}p_{t-s}^{(d)}(\xi-\x)H_A(\xi, s; E) m_R(d\xi)(1+o(1)),
\]
where $o(1)\to 0$ as  $x \to\infty$ (cf. \cite[Lemma 5.1]{Ucal_b}  if necessary).
 In the integrand we may replace $t-s$ by $t- \eta$ for the lower bound and    integration  of $H_A(\xi,s;  E)$ over $0< s<\eta$  yields $P_\xi[B_{\sigma(A)}\in E, \sigma_A \leq \eta]$.
We note  $p_{t-\eta}^{(d)}(\xi-\x)= e^{\xi\cdot \x/(t-\eta)}p_{t-\eta}^{(d)}(x)(1+ O(1/t))$, and deduce that
$$\int_{\xi\cdot \x >0} \frac{\xi\cdot\x}{t- \eta}\, e^{\xi\cdot \x/(t- \eta)}m_R(d\xi) \geq \k_d \int_0^{1} \frac{Rx}{t-\eta}e^{ \frac{Rx}{t-\eta}(1-\frac12  \th^2) }\th^{d-2}d\th \geq \k_d'\La_\nu\bigg(\frac{Rx}{t-\eta}\bigg).$$
Using  Theorem \ref{thm01}  we now find the lower bound of the lemma.

 In the case  $x\leq t$ we have
 $H_A(\x,t-s;d\xi)\geq \k_d q(x,t-\eta; R)m_R(d\xi),$ $s<\eta$ and 
the same argument as above immediately leads to the result.
\qed

\section{The  Wiener sausage for a Brownian bridge}

Given a compact  set $K$, let $S_K(t)$ be a Wiener sausage of length $t$ swept by  $K$ attached to a Brownian motion $B_t$: 
$$S_K(t) = \{\z\in \R^d: \z- B_s  \in  K \,\,\,\mbox{for some} \,\,\, s\in [0,t) \}.$$
The $d$-dimensional volume of a Borel set $A$ is denoted by  vol$_d(A)$. We sometimes write area$(A)$ for  vol$_2(A)$. In this section   $(d)$ is usually suppressed   from $p^{(d)}_t$, $R_K=\sup\{|\y|: \y\in K^r\}$ and $K$ always denotes a non-polar compact set of $\R^d$  satisfying  (\ref{cavity}) with  $K$ in replace of $A$.

\subsection {Case  $\x/t\to 0$} 

 \begin{Thm}\label{thm2}\,   If  $d\geq 3$, uniformly for   $\x$, 
as $t\to\infty$ and  $x/t\to 0$,  
 \beqn\label{EQ1}
E_{{\bf 0}} \,[\, {\rm vol}_d(S_K(t)) \,|\, B_t=\x]  \,\sim\, \cp(K) t.
\eeqn
\end{Thm}

 \begin{Thm}\label{thm3} \,   If  $d=2$, uniformly for   $\x$, 
as $t\to\infty$ and  $x/t\to 0$,  
   \begin{eqnarray}\label{EQ2} \qquad
E_{{\bf 0}} \,[ {\rm area}(S_K(t))  \,|\, B_t=\x]  \,=\,  \left\{ \begin{array}{ll} 
{\displaystyle
\frac{2 \pi t }{\lg t} (1+o(1)) }\quad\quad&\mbox{if}   \quad x\leq  \sqrt t, \\[5mm]
 {\displaystyle
\frac{\pi t}{\lg (t/x)} (1+o(1)) }\quad\quad&\mbox{if}   \quad x> \sqrt t.
\end{array} \right.
\end{eqnarray}
 \end{Thm}
 \v2

{\sc Remark 5}.  (a) \, In the case when  $B_t$ is pinned at $\x={\bf 0}$ an asymptotic expansion of $E_{{\bf 0}} \,[\, {\rm vol}_d(S_K(t)) \,|\, B_t= {\bf 0}] $ is obtained in \cite{M} ($d\geq 3$) and \cite{M2} ($d =2$) (see also \cite{BB}). 

(b)\, In 
 {\rm \cite{Usaus}},  it is shown that if $d=2$,
 for each $M>1$, uniformly for $|\x|\leq M\sqrt{t}$, 
  $$  E_{{\bf 0}}\, [ {\rm area} (S_{U(a)}(t))\,|\, B_t=\x] = 2\pi t N(\k t/a^2) + \frac{\pi x^2}{(\lg t)^2}\Big[\lg \frac{t}{x^2\vee 1} +O(1)\Big]  +O(1) $$
  as $t\to \infty$,
where    $\k=2e^{-2\ga}$ and $N(\la) = \int_0^\infty e^{-\la u}  [(\lg u)^2 +\pi^2)]^{-1} u^{-1}du$,which admits asymptotic expansion in powers of $1/\lg t$:  $N(\a t) = \frac1{\lg t} - \frac{\ga+\lg \a}{ (\lg t)^{2}} + \frac{(\ga+\lg \a)^2 +\pi^2/6 }{ (\lg t)^{3}}  + \cdots$ ($\a>0, t\to\infty$).

\v2
 Let $-K+\z$ denote the set $\{-\xi + \z: \xi\in K\}$.

 \begin{Lem}\label{lem4.1}   
 For  $0 \leq \a<\b\leq t$, $\x\in\R^d$ and a Borel set $W\subset \R^d$ such that $K\cap W=\emptyset$,
 \begin{eqnarray}\label{AreaK}
&& E_\0[{\rm vol}_d (\{\z\in W: \a\leq  \sigma_{- K+\z} < \b\})\,|\, B_t=\x] \nonumber\\
&&= \frac1{p_t(x) }\int_{\z\in W}|d\z| \int_{\a}^{\b} ds\int_{\xi \in \partial K } p_{t-s}(\z-\x-\xi) H_K(\z,s; d\xi),  
\end{eqnarray}
 where  $|\cdot|$ designates the Lebesgue measure on $\R^d$.
\end{Lem}
\v2\n
\pf\,  
The set $\{\z\in W: \a\leq  \sigma_{- K+\z} < \b\}$, depending on Brownian path $\om :=(B_s)_{0\leq s<t}$, is considered as  an $\om$-cross-section of the
set 
\beqn\label{mble}
 \{(\z, \om):  (\z, \sigma_{- K+\z})  \in W\times [\a, \b)\},
 \eeqn 
 which is measurable  w.r.t. the product $\sigma$-field ${\bf B(\R^d)} \times \sigma (B_s\!:s\leq t)$ and of which  the conditional probability of the $\z$-cross-section is given by
\beq
&& P_\0[ \sigma_{-K +\z} \in [\a, \b)\,|\, B_t = \x]\\
&& =\frac1{p_t(x)}\int_{ [\a,\b) \times \partial (-K+\z)} p_{t-s}(\x-\xi) P_\0[ \sigma_{-K +\z} \in ds, B(\sigma_{-K+\z})\in d\xi] \\
&& = \frac1{p_t(x)} \int_{\a}^{\b} ds \int_{\xi\in -\partial K} p_{t-s}(\x-\z-\xi)H_{-K}(-\z,s; d\xi).
\eeq
Now, applying the equality  $\int_{ -\partial K} \fa(\xi) H_{-K}(-\z,s; d\xi)= \int_{\partial K} \fa(-\xi)  H_{K}(\z,s; d\xi)$   as well as Fubini's theorem we conclude  the formula to be  shown. \qed
\v2

Put
\beq
F_K(t, \x) =\frac1{p_t(x) } \int_{\z  \in \Om_K}|d\z| \int_0^t ds\int_{\xi\in \partial K}p_{t-s}(\z-\x-\xi) H_K(\z,s; d\xi). 
\eeq
Since 
$S_K(t)\cap \Om_K  
= \{\z\in \Om_K:  \sigma_{- K+\z} < t\}$ valid under the condition (\ref{cavity}), i.e., $\R^d\setminus \Om_K = K^r$,
we obtain the following corollary of Lemma \ref{lem4.1}.
 \begin{Cor}\label{cor4.-2} \quad
$
E_\0[{\rm vol}_d (S_K(t))\,|\, B_t=\x]   = F_K(t, \x)+  {\rm vol}_d(K). 
$
\end{Cor}

We use the monotone class theorem to  extend the formula (\ref{AreaK})
 to  space-time Borel sets. The result is stated as the  following corollary,  which  we shall need   in the  proof (Step 2) of Lemma \ref{lem0a}.
 \begin{Cor}\label{cor4.-1}  If $D$ is a Borel set of $(\R^d\setminus K) \times [0, t)$, then
 $$E_\x\Big[{\rm vol}_d\Big(\!\{\z: (\z,\sigma_{- K+\z})  \in D\}\!\Big) \Big| B_t=\x\Big] = \frac1{p_t(x)} \int\!\!\!\int_{D}|d\z| ds\int_{\partial K}p_{t-s}(\z-\x-\xi) H_K(\z,s;d\xi). $$
 \end{Cor}

\v2
Put
$$F_K^*(t, \x) =  \frac1{p_t(x) }\int_{\z  \in \Om_K}|d\z| \int_0^t p_{t-s}(\z-\x) ds\int_{\xi\in \partial K}e^{-\x\cdot\xi/t-|\xi|^2/2t}H_K(\z,s;d\xi).$$
By the scaling property of Brownian motion we 
have
\beqn\label{scl_F}
F_K(t,\x) =r^d F_{K/r}(t/r^2,\x/r),
\eeqn
 and similarly for  $F_K^*$. The function $F_K^*$ is easier to deal with than  $F_K$, and the difference of the two  is negligible for the present purpose as one can read off from the following lemma.

\begin{Lem} \label{lem3.50} \, Let $M$ be a positive constant.   For $R_Kx/t <M$ and $t>R_K^2$,
$$|F_K(t,\x)  - F_K^*(t, \x)| \leq c_M\ga_K R_K^2  \Big [1+ \sqrt t  /R_K\Big],$$
where   $c_M$ is a constant depending only on $ M$ and $d$, and
$$\gamma_K =1\,\, \mbox{ or  \,\,$\cp(K)$\quad  according as \quad $d=2$\,\, or}\,\,\, d\geq 3.$$  
\end{Lem}
\v2\n
\pf\, The proof is similar to that of Lemma 4.3 of \cite{Usaus}. 
Let $h > 4 R_K^2$ be a  constant  that will be  suitably  chosen later depending on  $M$ (see the  part (c) below). Putting  $T:=s(t-s)/t$, we bring in the following sub-regions of $[0,t] \times \Om_K$
 \beq
&&D^h= \{(s,\z): t-h\leq s < t, \, \z\in \Om_K \}, \quad\,\\
&&D_h = \{(s,\z): 0 \leq s < a,\, \z\in \Om_K \}, \quad\,\\
&&D_> =\{(s,\z): h\leq s < t-h, \,   |\z-(s/t)\x| \geq \sqrt{8T \lg (T/R_K^2)}, \, \z\in \Om_K\},\\
&&D_< =\{(s,\z): h\leq s<t-h, \,  |\z-(s/t)\x| < \sqrt{8T \lg (T/R_K^2)}, \,\z\in \Om_K \},
\eeq
and restrict to them the   integrals that define  $F_K$ or $F^*_K$.
 Denote the ratios to $p_t(x)$ of the corresponding integrals for  $F_K$  by $J\{D^h\}, J\{D_h\}$,  etc.  and those for $F_K^*$   by $J^*\{D^h\}, J^*\{D_h\}$,  etc. 
 
Suppose $d\geq 3$.  Of the first three regions above we evaluate the corresponding  integrals $J$ and $J^*$ separately,
 and verify  that they are all  bounded from above by 
 \beqn\label{remark}
  c_M R_K^{2}\cp(K) 
  \eeqn
   in the paragraphs (a), (b) and (c) below. The actual  computations are given only for $F_K$, those for  $F_K^*$  being   similar  and much simpler.     
   From  Proposition \ref{prop6.3} we have
    \beqn\label{00}
 q_K(\z, t) \leq \k_{d} \cp(K)  p_{t}(\z)e^{2R_K|\z|/t} \quad\mbox{for}\quad t>R_K^2, \z  \in \Om_K.
 \eeqn
  For simplicity let $R_K=1$  in what follows, which gives no loss of generality because of the scaling property (\ref{scl_F}),  and    suppose  $t> 2h$  and $ x< Mt$.
  Denote by  $c_M, c_M'$ etc. unimportant  positive constants that depends only on $M$ and $d$.
\v2
(a) \,$D^h$:  Use the inequality 
 $$H_K(\z,s;d\xi) \leq \int_{\Om_K}  p_1(\y-\z) H_K(\y, s-1;d\xi) |d\y|\quad\quad\quad (s>1)$$
  and perform  integration by $\z$ first. Then, after changing variable of integration,   we have 
  $$J\{D^h\}\leq \frac1{p_t(x) }\int_0^ h ds \int\!\!  \int_{\Om_K\times \partial K} p_{s+1}(\y-\x-\xi)H_K(\y, t-s-1;d\xi)|d\y| .$$
 Using  the inequality $|\a-\xi|^2\geq \frac12 |\a|^2 - |\xi|^2$ ($\a\in \R^d$) we find  that $p_{s+1}(\y-\x-\xi)\leq 2^{d/2} e^{|\xi|^2/2(s+1)}p_{2s+2}(\y-\x)$. Thus  the inner double  integral over $\Om_K\times \partial K$ is dominated by
 $$2^{d/2}e\int_{\Om_K}p_{2s+2}(\y-\x)q_K(\y,t-s-1)|d\y|,$$ 
and owing to (\ref{00})   we  obtain
$J\{D^h\} \leq c'_{M}\cp(K) \int_0^ h p_{t+s+1}(x)ds/p_t(x) \leq c_{M}\cp(K)$
as desired,  where we  also have applied the assumption  $x<Mt$.

\v2
(b) \,$D_h$:\, Noting $\sup_{\xi\in K} p_{t-s}(\z- \x -\xi) <\k_d  p_{t+1}(\z-\x)$  ($s<h$) and $\sigma_K \geq \sigma_{U(1)}$ we see 
$$
J\{D_h\}  \leq  \frac{\k_d}{p_t(x) }\int_{\Om_K} p_{t+1}(\z-\x) P_\z[\sigma_K< h]  |d\z|  \leq \frac{\k_d}{p_t(x) }  P_\x[t+1\leq \sigma_{K} \leq t+1+h],$$
and   applying   (\ref{00})  again we conclude  $J\{D_h\} \leq c_{M}\cp(K)$.
\v2
For the rest  of this proof we  shall use the identity
\beqn\label{id}
p_s(\z)p_{t-s}(\z-\x) =p_t(\x)p_{T}\bigg( \z -\frac{s}{t}\x \bigg),
\eeqn
where
$$T=T(t,s) = \frac{s(t-s)}{t}.$$
(the second factor on the right-hand side of (\ref{id}) is the probability density of the Brownian bridge at time  $s\in  (0,t)$).

\v2
(c) \,$D_>$:\, Using the inequality  $p_s(x) e^{2x/s}\leq \k_d\, p_{s+1}(x)$ ($s>1$) as well as (\ref{00}) we have
 \beq 
 J\{D_>\} \leq \frac{c_M\cp(K)}{p_t(x)} \int_ h^{t-h} ds \int_{|\z- (s/t)\x|>\sqrt{8T\lg T}}  p_{s+1}(\z)p_{t-s+1}(\z-\x) |d\z|. 
  \eeq
On writing   $t'=t+2, s'=s+1$ and $T'=s'(t'-s')/t'$ the  integral on the right-hand side divided by  $p_{t+2}(x)$ equals
 $ \int_ h^{t-h} ds\int_{|\z- (s/t)\x|>\sqrt{8 T \lg T}}  p_{T'}(\z-\frac{s'}{t'}\x)|d\z|$. Noting  
   $|s'/t' -s/t| x \leq x/t <M$,   choose  $h=h_M>4$  so that $\sqrt{8T\lg T} -M \geq \sqrt{5T'\lg T'}$ for $h<s<t-h$, and  we can then replace the range of  the inner integral by $|\z-\frac{s'}{t'} \x|>\sqrt{5T'\lg T'}$ to deduce  that
 $$  \frac{J\{D_>\}}{\cp(K)} \leq c_M\int_h^{t-h} ds\int_{\sqrt{5\lg T'}}^\infty e^{-u^2/2}u^{d-1}du \leq c'_M\int_{h+1}^{t'-h}\frac{ds'}{[s'(t'-s')/t']^2}  \leq c''_M h^{-1}$$
   as required.
\v2

(d) \,$D_<$:  We continue to suppose  $R_K=1$.  Using (\ref{id}) we observe
\begin{eqnarray}\label{B1}
&&\frac{p_s(\z)}{p_t(\x)}\Big[p_{t-s}(\z- \x-\xi) - e^{-\x\cdot\xi/t-|\xi|^2/2t}p_{t-s}(\z-\x) \Big] \\
&&=\frac{p_t(\x+\xi)}{p_t(\x)}\left[  p_{T}\Big(\z- \frac{s}{t}(\x+\xi)\Big) -  p_{T}\Big(\z- \frac{s}{t}\x\Big)  \right]. \nonumber
\end{eqnarray}
For $(s,\z)\in D_<$ we have  $|\z-\frac{s}{t}\x|/T \leq \sqrt{8(\lg T)/T} < 3$ and, on using  $|e^{\a-p}-1| \leq e^{|\a|}(|\a|+p)$ ($p\geq 0$),
\begin{eqnarray}\label{B2}
&&\bigg|p_{T}\Big(\z- \frac{s}{t}(\x+\xi)\Big) -  p_{T}\Big(\z- \frac{s}{t}\x\Big)  \bigg| \\
&&=
p_{T}\bigg(\z-\frac{s}{t}\x\bigg)  \Bigg| \exp\Bigg\{\frac{s/t} { T}\bigg[\bigg(\z-\frac{s}{t}\x\bigg)\cdot \xi - \frac{s|\xi|^2}{2t}\bigg] \Bigg\}-1\Bigg|\nonumber \\
&&\leq e^3  p_{T}\bigg(\z-\frac{s}{t} \x\bigg)    \frac{1}{t-s}  \Bigg[|\xi|\bigg|\z- \frac{s}t \x\bigg| + \frac{s|\xi|^2}{2t}\Bigg]. \nonumber
\end{eqnarray}
Hence, if $x\leq Mt$,
 \begin{eqnarray*}
&&|J\{D_<\} - J^*\{D_<\}|\\
&&\leq e^3 e^{ M} \int_{D_<}    \frac{1}{p_s(\z)}\cdot p_{T}\bigg(\z-\frac{s}{t} \x\bigg)    \frac{1}{t-s}  \Bigg( \bigg|\z- \frac{s}t \x\bigg| + \frac{s}{2t}\Bigg) q_K(\z,s)ds  |d\z| .
\end{eqnarray*}
 On $D_<$ we have  $|\z| < (M+1)s$ and, noting $\int p_T(\y)|y||d\y| \leq \k_d \sqrt T$,  we   apply (\ref{00}) to  deduce that the integral above  is at most  a constant multiple of  
\begin{eqnarray}\label{A10}
\cp(K) \int_ h^{t-h}  \frac{1}{(t-s)} \bigg[ \sqrt T + \frac{s}{t}\bigg]ds 
  \leq C' \cp(K){\sqrt t}.
\end{eqnarray}
This completes the proof of Lemma \ref{lem3.50} if $d\geq 3$.
In the case $d=2$ one may go through  with the same proof except that he uses Proposition \ref{lem6.9} in place of Proposition \ref{prop6.3}.  \qed

 \v2
 
\n
{\it Proof of Theorems \ref{thm2} and \ref{thm3}.}\, 
 Owing to Lemma \ref{lem3.50} as well as Corollary \ref{cor4.-2}  it suffices to show that   as $x/t\to 0$ and $ t\to \infty$  the function $F_K^*(t, \x) $, which may obviously  be written as  
\beqn\label{*21}
  \frac1{p_t(x) } \int_0^t  ds\int_{\z\in \Om_K}p_{t-s}(\z-\x) q_K(\z,s)|d\z| ( 1+O(x/t)),
\eeqn
has the same asymptotic form as given on the right-hand side of (\ref{EQ1}) if $d\geq 3$ and  that of (\ref{EQ2}) if $d=2$.  Put $D_t= \Big\{(s,\z): 4 <s <t/2,  \, |\z-(s/t)\x| < \sqrt{2s\lg s }, \z\in \Om_K\Big\}.$

Let $d\geq 3$.  Considering  $D_t$ in place of  $D_<$
we argue as in the proof of   Lemma \ref{lem3.50}. Observe  $\int_{|\z-(s/t)\x| \geq  \sqrt{2s\lg s}}\,   p_{T}(\z-\frac{s}{t} \x)  |d\z| = O( 1/ s) $, $s>1$
(valid even  if $s, t$ are replaced by $s+1, t+1$) and 
use Theorem \ref{thm3.1} to see that the inner  integral  restricted to $D_t$ equals
\beq
&&\cp(K) p_t(\x)  \int_{|\z-(s/t)\x| <  \sqrt{2s \lg s}}\, P_\z[\sigma_K=\infty]  p_{T}\bigg(\z-\frac{s}{t} \x\bigg)  |d\z| (1+o(1))\\
&&= \cp(K)  p_t(\x)(1 +o(1)),
 \eeq
 where in the last line $o(1)\to 0$ as $s\to\infty$.
 Hence  the integral on the lower half interval $0<s<t/2$ gives half the asserted leading term in Theorem \ref{thm2}. The other half is dealt with in an analogous way and the details are omitted.
 
The case  $x<M\sqrt t$ of  $d=2$    follows from the result for $U(a)$ given in Remark 5 (b)  on noting  that $q_K$ has the same asymptotic form as $q_{U(1)}$ (for large space variable) owing to  Theorems \ref{thm01}  and \ref{thm3.2} and that the repeated integral (\ref{*21}) restricted on the outside of $D_t$ is negligible. 
As for the other case of $d=2$, of which the proof is somewhat involved if argued as above, we apply Proposition \ref{prop5} of the next subsection, which immediately implies the asserted result in view of Theorem \ref{thm3.2}.
 \qed

\vskip2mm


\subsection{Case $x/\sqrt t\to\infty$}
Here we consider the case $x/\sqrt t \to\infty$, mostly the case  $x/t\to \infty$.  The main results  are  presented  in the first part  of this subsection and   proofs of them will be given in the second  through fourth parts.  Throughout this subsection  we fix a unit vector $\be\in \partial U(1)$.
\v2

{\bf 4.2.1.} {\sc Statements of results.} \, 
Combined with Theorems \ref{thm3.1} and \ref{thm3.2} the next proposition  covers the case  $x/\sqrt t\to\infty$  of  Theorems  \ref{thm3} as noted at the end of the preceding  subsection. It also plays a substantial  role in the proof of 
Theorem \ref{thm3.4} and Theorem \ref{thm4.3} below. $K$ continues to denote a non-polar compact  set of  $\R^d$.
\begin{Prop} \label{prop5}\,  As $x^2/t\to\infty$ and $ t\to\infty$ 
\beq
E_{{\bf 0}} \,[ \, {\rm vol}_d(S_K(t)) \, | \, B_t= \x]  = \frac{t}{p^{(d)}_t(x)  }\cdot\frac{ E_\x[ e^{- B_{\sigma(K)}\cdot {\bf x}/t}; \sigma_K \in dt ] }{ dt}+ o(\zeta_d(x,t)),
\eeq
where $\zeta_d(x,t) = t\vee x$  if $d\geq 3$ and  $\zeta_2(x,t) = t/\lg(t/x)$ $(x<t/2)$; $= x$ $(x\geq t/2)$.
\end{Prop}
\v2
\begin{Thm} \label{thm4.3}\,  Let $d\geq 2$,  $\x=x\be$ for $\be\in \partial U(1)$ fixed. Suppose   the  condition (\ref{cond1})  to be  satisfied. Then as $x/t\to\infty$ and $ t\to\infty$, 
\beqn\label{Eq_thm4.3}
E_{{\bf 0}} \,[ \, {\rm vol}_d(S_K(t)) \, | \, B_t= \x] = {\rm vol}_{d-1}({\rm pr}_{\be}K) x +o(x).
\eeqn
\end{Thm}
\v2

The proofs of these two  results and that of Theorem \ref{thm3.4} are interrelated in a way. Our proof of Proposition \ref{prop5} is quite involved: the difficulty occurs when  $x/t$ becomes indefinitely large,       otherwise the proof being  easy.
The upper bound for the asymptotic relation (\ref{Eq_thm4.3}) is relatively easy and the result is useful  for  the proof of  Proposition \ref{prop5} in the case $x/t\to\infty$.  The lower bound of some cases of  $K$ of special shape is easy  and it together with the upper bound is used for the proof  of Theorem \ref{thm3.4}, which combined with Proposition \ref{prop5} immediately gives Theorem \ref{thm4.3}. 

In the rest of this part  we give a proof of the assertion  (vi) of Section 1,  which concerns the case when  $x/t$ approaches a positive constant $v$, and need  Proposition \ref{prop5}. 
We remind the reader that  for each $a >0$, the function
$g_{av}(\th)=\sum_{n=0}^\infty \frac{K_0(av)}{K_n(av)} H_n(\th)$
 is a probability density   on $\partial U(1)$ with respect to   $m_1(d\xi)$ on the understanding that    $\th=\th(\xi)$ is the colatitude of $\xi$ relative to the (arbitrarily chosen)   vector $\be$ which is   taken for the north pole (see (\ref{h_v})).
 
\v2\n
{\it Proof of  {\rm (vi)} of {\sc Summary of Main Results}.}  
\v2\n
 The second formula of  (vi) follows   immediately  from the first owing to Proposition \ref{prop5}.   As for the first one we recall the Huygens  decomposition (\ref{-70})  with $R=R_A$ (under the convention that
for $\y\in A^r$,  $H_A(\y,s; \Ga) = \de(s)\1(\y\in \Ga)$ where $\de(s), s\geq 0$ is the Dirac delta function)  
and
apply  Theorems \ref{thm01} and \ref{thm2.3}. 
We split  the outer integral at $\sqrt t$ and $t/2$ as in the proof of Lemma  \ref{cor-1} and denote the corresponding integrals by 
$I_{[0,\sqrt t]}(\Ga)$, $I_{[\sqrt t, t/2]}(\Ga)$  and  $I_{[t/2, t]}(\Ga)$. (Here we should take  $\sqrt{t}/R_A$ in place of $\sqrt t$, but we do not do that since the exponent $1/2$ is rather arbitrarily chosen from $(0,1)$.) Plainly, 
$$I_{[\sqrt t,t/2]}(\partial A) \leq \sup_{\sqrt t<  s\leq t/2} q(x,t-s; R_A) \int_{\sqrt t}^{t/2}\sup_{\xi\in \partial U( R_A)}q_A(\xi,s)ds, \quad\mbox{and}
$$
 $$I_{[t/2,t]}(\partial A) \leq  P_\x[ \sigma_{U( R_A)} <\tst12 t]\sup_{\xi\in \partial U( R_A), \,t/2\leq s\leq t}q_A(\xi,s).
$$
By Propositions \ref{prop6.4} (d=2) and \ref{prop6.3} ($d\geq 3$) and Theorem \ref{thm01}
\beqn \label{adbd}
q_A(\xi,s) \leq C_{A}q(2 R_A,s; R_A)\qquad (\xi\in \partial U( R_A), s >  R_A^2), 
\eeqn
where $C_{A} = C \b_A e_A(\xi)$ if $d=2$  and $= \k_d \cp(A)/R_A^{2\nu}$ if  $d\geq 3$ (note that $q_A(\xi, s)=0$ if $\xi \in A^r, s>0$).
Write $\tilde v = x/t$. Noting that $p_{t-s}(x)$ is decreasing in $s$ whenever  $0< t-s < d^{-1}x^2$,  that 
 $$p_{t-\sqrt t}(x)/ p_{t}(x) \leq (1-{\textstyle\frac1{\sqrt t}})^{-d/2}  e^{-x^2/2t\sqrt t} \leq 2e^{-\frac12 \tilde v^2\sqrt t},$$
 and that $\La_\nu$ is increasing,
 we then  infer  that 
 $$I_{[\sqrt t,t/2]}(\partial A) \leq 
2C_A R^{2\nu}\La_\nu(R\tilde v) p_t(x) e^{- \frac12 \tilde v^2\sqrt t} \int_{\sqrt t}^\infty q(2R_A,s; R_A)ds.$$ 
 In view of  Theorem \ref{thm2.2} and  the upper bound of $q_A(\x,t)$  given in  (\ref{adbd}), $I_{[t/2,t]}(\partial A)$  
 is much smaller than  the right-hand side above because $p_{t/2}(x) \asymp p_t(x) e^{-\tilde v^2t/2}$.
Putting these together we apply Theorem \ref{thm2.3} to conclude that as  $x/t\to v$
  \beq
&&H_A(\x,t;\Ga) /p_t(x)\\
&&= \, R_A^{2\nu}\La_{\nu}(R_A \tilde v)\Bigg[\int_0^{\sqrt t} \frac{p_{t-s}(x)ds}{p_t(x)}\int_{\partial U(1)}H_A(R_A\xi,s;\Ga)  g_{R_Av}(\th_{\xi,\x})m_1(d\xi) (1+o(1)) + \e_{t,x}(\Ga)\Bigg],
\eeq
where  $\e_{t,\x}$ is a measure  on $\partial A$ such that
\beqn\label {err_bd}
\e_{t,\x}(\partial A)\leq 2C_A  e^{-\frac12\tilde v^2\sqrt t}\int_{\sqrt t}^\infty q(2R_A,s; R_A)ds
\eeqn
and  $\th=\th_{\xi,\x}$ is the colatitude of $\xi$ relative to  $\be=\x/x$.  The measure kernel $ \la_A$ introduced  in  Main result  II of Section 1 may be written as
 \begin{eqnarray}\label{lamda}
 \la_A({\bf v}; \Ga)   
= \int_0^{\infty} e^{-v^2 s/2}ds\int_{\partial U(1)}H_A(R_A\xi,s;\Ga)g_{R_Av}(\th_{\xi,{\bf v}}) m_1(d\xi)   \nonumber
 \end{eqnarray}
for  $v= |{\bf v}|>0$.  Noting  $p_{t-s}(x)\,\sim\, p_t(x) e^{-\tilde v^2s/2}$, we then deduce from 
  the expression of $H_A(\x,t;\Ga) /p_t(x)$ given  above  that  for each  $\Ga\subset \partial A$ with $H^\infty_A(\Ga)>0$, as $\tilde v =x/t \to v>0$
$$H_A(\x,t;\Ga) \sim  p_t(x)R_A^{2\nu}\La_{\nu}(R_A x/t) \la_A(\x/t; \Ga).$$ 
 For $\Ga=\partial A$, this may be written as  $ \la_A({\bf v}; \partial A) \sim q_A(\x,t)/p_t(x)R_A^{2\nu}\La_{\nu}(R_Ax/t)$. Hence, using  Propositions \ref{prop6.3} and \ref{prop6.4} again  we obtain
 $$ \la_A({\bf v}; \partial A) \leq c_M \gamma_A/R_A^{2\nu},$$
 where $\gamma_A $   $=1$  if  $d=2$;  $= \cp(A)$ if  $d\geq 3$  as in Lemma \ref{lem3.50}.
By the identity  
 \begin{eqnarray}\label{lamda5}
 \la_A(\0; \Ga)&=&\int_{\partial U(1)} P_{R_A\xi}[B_{\sigma_A}\in\Ga, \sigma_A<\infty]m_1(d\xi) \nonumber\\
 &=&  \left\{\begin{array}{ll} H^\infty_A(\Ga)\quad & (d =2) \\[2mm]
G^{(d)}(R_A) \cp(A) H^\infty_A(\Ga)\quad &(d\geq 3)
\end{array} \right.
\end{eqnarray}
(see (\ref{H/H}) for the second equality) and  the continuity of  $\la_A({\bf v}; \cdot)$  at ${\bf v} =\0$ one   observes  that   the  two asymptotic formulae of (vi) are valid  uniformly  about  $v= 0$ and  in particular they recover the results in the regime  $x/t\to 0$ given   in (iii), (iv) and the second half of  (v).
 
 The asserted
 uniform convergence in total variation norm is now easily checked  in view of  (\ref{err_bd}). 
Formula (\ref{EQ24})  is an immediate consequence of Theorem \ref{thm2.3}
and Proposition \ref{prop5}. \qed

\v2
{\bf 4.2.2.} {\sc Upper bound for  Theorem \ref{thm4.3}.}\, 
The upper bound for  the asymptotic relation (\ref{Eq_thm4.3})  follows from the next lemma as we shall see after its proof.  The proof of the lower bound is postponed to   Section  4.2.4.

\n
\begin{Lem} \label{prop6}  For $a>0$ and $h>0$  put  $C=\{\z: 0\leq \z\cdot {\bf e} \leq h$\,and\, $|{\rm pr}_{\bf e}\,\z| \leq a\}$ (a circular cylinder of height $h$, radius $a$ and axis $\bf e$). Then,  
 uniformly for  $x>t\vee h$, $a>0$ and $h>0$,  as $t\to\infty$
 \v2
 {\rm  (i)} \, $E_\0[ \,{\rm vol}_d(S_{C}(t))\, | \, B_t= x{\bf e}]  = c_{d-1}a^{d-1}x + O\Big(a^{d-2} \sqrt{(h\vee1)tx}\, +a^{d-1}t \, \Big)$; and
\v2  
{\rm  (ii)} \, $\Big(E_\0[ (\,{\rm vol}_d(S_{C}(t)))^2\, | \, B_t= x{\bf e}] \Big)^{1/2}  =c_{d-1}a^{d-1}x + 
O\Big(a^{d-2} \sqrt{(h\vee1)tx}\, +a^{d-1}t \, \Big)$, 
  \v2\n
  where  $c_n$ denotes  the volume of  $n$-dimensional unit ball.
\end{Lem}
\v2\n
\pf\, We give a proof only for $d=2$, the higher dimensional case being essentially the same.  Let  $\be =\be_0=(1,0)$ and write $B'_t$ and  $B''_t$ for the first and the second component of $B_t$ so that $B_t=(B'_t, B''_t)$.  For the first assertion (i) the lower bound is obvious: the expectation in it is bounded below by $c_{d-1}a^{d-1}x$. To see the upper bound  let $\e>0$ and put $N=\lceil t/\e\rceil$ (for $t>\e$) and $t_k=kt/N$, so that $\e(1-\e/t) < t_{k+1}-t_{k} \leq \e$. Let $v=x/t$. Noting that the conditional law of  $(B_s)$ given $B_t=x\be_0$ equals that of  $(B_s+vs\be_0)$ given $B_t=\0$, we bring in the variables  
$$ \xi^+_{k}= \sup_{t_{k}\leq s\leq t_{k+1}} (B'_s  + v s),\quad\quad \xi^-_k =\inf_{t_{k}\leq s\leq t_{k+1}} (B'_s  + v s); \, \mbox{and} $$
$$  \eta_k=  \sup_{t_{k}\leq s\leq t_{k+1}} |B''_s-B''_{t_k}|.\quad  $$ 
Write $\De_{k} S$ for $S_C(t_{k+1})\setminus S_C(t_{k})$  and $Y_k^\pm= B''_{t_k}\pm a$.  
We decompose $\De_k S$ into two parts  $\De_k S  \setminus (\R\times [Y^-_k,Y_k^+])$ and $\De_k S  \cap (\R\times [Y^-_k,Y_k^+])$. For the first part we have
\beq
 \De_k S  \setminus (\R\times [Y^-_k,Y_k^+])
\subset \Big( [\xi^-_{k},\xi^+_k +h]\times  ( [Y_k^+, Y_k^++\eta_k] \cup [Y_k^- -\eta_k, Y_k^-]) \Big) 
\eeq
of which the volume of the right-hand side is  at most
\beqn\label{est_vol}
2 (h+\xi^+_k -\xi^-_k)\eta_k. 
\eeqn
As for the second part  note that $\xi_k^- \leq \xi_{k-1}^+$ and 
 $|Y_k^+-Y_{k-1}^+|\vee |Y_k^- -Y_{k-1}^-|\leq \eta_{k-1} $
   and  that $S_C(t_k)$ includes the rectangle
 $$    [\xi_{k-1}^-, \xi_{k-1}^++h]\times
 [Y_{k}^- +2\eta_{k-1}  , Y_{k}^+ - 2\eta_{k-1}]$$
($[s,t]$ is understood to be the empty set  if  $s>t$), and  we then  deduce   that 
 \beq 
 &&\De_k S  \cap (\R\times [Y^-_k,Y_k^+]) \\ &&\subset ([\xi^+_{k-1} +h, \xi_k^+ +h]\times [Y_k^-, Y_k^+])\bigcup ([\xi^-_{k-1}\wedge \xi_{k}^-, \xi_{k-1}^-]\times [Y_k^-,Y_k^+])\\
&&\qquad  \bigcup  \Big([\xi^-_{k-1}, \xi_{k-1}^+ +h] \times \Big([Y_{k}^+ - 2\eta_{k-1}  , Y_{k}^+]  \cup [Y_{k}^- , Y_{k}^- +2\eta_{k-1} ] \Big)\Big) 
\eeq
($k=1,\ldots, N-1$)  
 and $\De_0S  \cap (\R\times [Y^-_0,Y_0^+])  \subset [\xi_0^-,\xi^+_0]\times [-a,a]) $.
Combined  with the bound (\ref{est_vol})  this shows
\begin{eqnarray}\label{000}
{\rm area}(S_C(t)) &\leq&  \sum_{k=0}^{N-1}2a(\xi^+_k-\xi^+_{k-1}) + 6\sum_{k=0}^{N-1}  (h+\xi^+_k -\xi^-_k)\eta_k  \\
&& \quad + \, 2a\sum_{k=0}^{N-1} (\xi^-_{k} - \xi_{k+1}^-)\vee 0,  \nonumber 
\end{eqnarray}
where  $\xi_{-1}^+$ is understood to be  $\xi_0^-$.
The expectations under the conditional measure  given  $B_t=\0$ of the variables $\xi_k^+ -\xi_k^-$ and $\eta_k$    are  at most constant multiples of   $v\e +\sqrt\e$ and of $\sqrt \e $, respectively, for every  $k$: indeed,  if   $E_\0^0$ denotes the conditional expectation given  $B_0=B_t=\0$,
$$E_\0^\0[\eta_k^2 ]\leq 2E_\0\Big[ \sup_{t_{k}\leq s\leq t_{k+1}}|B''_s-B''_{t_k}|^2 + (\e t^{-1}|B''_t|)^2\Big] < C  \e + C'\e^2/ t\leq C'' \e $$
 and similarly for $\xi_k^+ -\xi_k^-$. Observe $\xi_k^- -\xi_{k+1}^- \leq -v\e -\inf_{t_{k+1} \leq s\leq t_{k+2}}(B_s'-B'_{t_k})$,  and  by  employing the Schwarz inequality and the expression of  the pinned Brownian  by means of a free Brownian motion  we infer that  
 $$E_\0^\0[(\xi_k^- -\xi_{k+1}^-)\vee 0] \leq C \sqrt \e \Big(P_0\Big[\sup_{t_k<s<t_{k+2}}B'_s-B'_{t_k}+ \e t^{-1}|B_t| >v\e\Big]\Big)^{1/2} \leq C'\sqrt \e e^{-\e v^2/16};$$
also  $$E_\0^\0[\xi_{N-1}^+ - \xi_0^-] \leq  x+  2E_\0^\0[\sup_{0 \leq s \leq \e} |B_s'|]  \leq x + C\sqrt \e.$$  Now   take  $\e= (h\vee 1)/v$, which entails  $\e<t$ for $x>h$.  Putting these bounds together, an elementary computation then  leads to 
$$\Big|E_\0 [{\rm area}(S_C(t))\,|\, B_t=x\be_0] - 2ax \Big| \leq C'''[at +  \sqrt {(h\vee 1)xt} \,]. $$
Thus    (i) is verified. It is easy to check that   $E_\0^\0[(\xi_k^- -\xi_{k+1}^-)\vee 0)^2] \leq C \e e^{-v/16}$  and  the computations carried out  above also deduce  the assertion (ii)  from (\ref{000}).   \qed
\v2

From the proof of Lemma \ref{prop6} its first assertion (i) may be slightly generalized. For  $0\leq \a<\b$,   put 
\beqn\label{saus3}
S_K[\a,\b) = \{\z\in \R^d: \z- B_s  \in  K \,\,\,\mbox{for some} \,\,\, s\in [\a,\b) \}.
\eeqn
Then  for  any  $\de\in (0,1]$,  if $[\a,\b)\subset [0,t)$ with  $\b-\a =\de t$, then
 \beqn\label{extension}
 \Big|E_\0[ \,{\rm vol}_d(S_{C}[\a,\b))\, | \, B_t= x{\bf e}] - c_{d-1}a^{d-1}\de x \Big|  \leq 
   C[at \de+  \sqrt {(h\vee 1)xt}\, \de + \sqrt{(h\vee 1)/v}\,]
 \eeqn
 (provided $x>t\vee h, h>0$). 
Since 
$\{\z\in \Om_K:  \sigma_{- K+\z} \in [\a,\b)\}\subset S_K[\a,\b)  \cap \Om_K$, by Lemma \ref{lem4.1}  we have
 \beqn\label{uppbd_int}
  \int_{\z  \in \Om_K}|d\z| \int_\a^\b ds\int_{\xi\in \partial K}\frac{p_{t-s}(\z-\x-\xi) H_K(\z,s;d\xi)}{p_t(x)} \leq E_\0[{\rm vol}_d(S_K[\a,\b)) \,|\, B_t=x\be ]. 
  \eeqn
We shall need an upper bound of the expectation on the right,  which may be stated as follows.

\begin{Cor}\label{cor4.3}  \,  For each  $\de\in (0,1]$ and $\be\in \partial U(1)$, uniformly for  $0\leq \a<\b\leq t$ with   $\b-\a=\de t$ and for $x\geq t$,  as \,$t\to\infty$,
$$ E_\0[{\rm vol}_d(S_{K}[\a,\b)) \,|\, B_t= x\be]  \leq  {\rm vol}_{d-1}({\rm pr}_{\be} K)\de x + c_{K}\de x/\sqrt{v} +o(x). $$ 
Here $o(\de x)$ may depend on $K$ and  $c_K$ is a constant depending only on $K$.
\end{Cor}
\v2\n
\pf \, Suppose   $\de=1$ so that $[\a, \b)=[0,t)$. Given $\e>0$, we can find  a finite number of   balls $ b_n \subset \De_\be, n=1,\ldots, N$ such that ${\rm pr}_\be K \subset \cup b_n$ and $\sum {\rm vol}_{d-1}(b_n) < {\rm vol}_{d-1}({\rm pr}_{\be} K) +\e$.  By   Lemma \ref{prop6} we have that  $ E_\0[{\rm vol}_d(S_K(t)) \,|\, B_t= x\be] /x \leq \sum_{n=1}^N\Big( {\rm vol}_{d-1}(b_n) +  O(1/\sqrt v)\Big)$, where  $O(1/\sqrt v)$ is independent of the choice of  $(b_n)$, entailing the asserted inequality of the lemma. For the case  $\de<1$ we have only to use (\ref{extension}) instead of Lemma \ref{prop6}.  \qed

\begin{Cor}\label{cor4.2}  \,  Let $K$ be a compact set of  $\R^d$  such  that $K$ lies on a plane perpendicular to  $\be $. Then  as $x/t \to\infty$ and $t\to\infty$,
$E_\0[{\rm vol}_d(S_K(t)) \,|\, B_t=x\be ]  = {\rm vol}_{d-1}(K)x +o(x). $
\end{Cor}
\n
\pf\,  The lower bound holds path-wise due to the hypothesis on $K$, whereas the upper bound follows from Corollary \ref{cor4.3}. \qed
\v2

\v2
{\bf 4.2.3.} {\sc  Proof of Proposition \ref{prop5}}.\, Let $K$ be   a compact set of $\R^d$ ($d\geq 2$). In terms of $H_K$
Proposition \ref{prop5}  may be stated as follows: as $x/\sqrt t\to\infty$ and $ t\to\infty$
\beqn\label{eqn1}
E_\0[ \,{\rm vol}_d\,(S_K(t))\, | \, B_t= \x]  = \int_{\xi\in \partial K} e^{-\x\cdot\xi/t} H_K(\x,t;d\xi) \frac{t}{p_t(x) }+ o(\zeta_d(x,t)).
\eeqn
The proof is performed by showing   Lemmas \ref{lem0a} through \ref{lem4.51} given below. Let  $h=4R_A^2$ and    for $2h<t_1\leq t$ define
\beqn\label{J1}
F_K^{\x,t}(t_1)=  \frac1{p_t(x) }\int_{\z\in \Om_K}|d\z|\int_h^{t_1 -h} p_{t-s}(\x-\z) ds \int_{ \partial K} e^{-\x\cdot\xi /t} H_K(\z,s; d\xi ).
\eeqn
\v2\n
\begin{Lem}\label{lem0a}  For each $\de\in (0,1]$,  as $x/\sqrt t\to\infty$ and $ t\to\infty$
  \beqn\label{J0}
  E_\0[ \,{\rm vol}_d\,(S_K(\de t))\, | \, B_t= \x] = F_K^{\x,t}(\de t)(1 +   o(1)) +o(x).
  \eeqn
\end{Lem}
\v2\n
\pf \,  The proof is given only for the case $\de=1$, the case $\de<1$ being dealt with in the same way.   Note that $F_K^{\x,t}(t)$ is substantially  the same as $F^*_K(t,\x)$ defined in Section 4.1 except for   the contribution to the latter integral from the intervals $[0,h]$ and $t-h,t]$  and the assertion of the lemma   follows from Lemma \ref{lem3.50}  when  $x<Mt$ (with any  $M>0$) (see the remark given at (\ref{remark})).  Thus  we may and do suppose $x/t>1$. By Lemma \ref{prop6}  $E_\0[ \,{\rm vol}_d\,(S_K(t))\, | \, B_t= \x] $ is then bounded by a positive multiple of $x$; by Corollary \ref{cor4.-2} this expectation is  expressed as 
\beqn\label{J}
  \frac1{p_t(x) }\int_{\z\in \Om_K}|d\z|\int_0^t ds\int_{\xi\in \partial K} p_{t-s}(\z-\x -\xi)H_K(\z,s; d\xi)  +{\rm vol}_d(K).
\eeqn

 The rest of the proof is carried out in five steps, where we let $R_A=1$ for simplicity.
 
\v2 
{\bf Step 1.}\,  We first consider  the above  integral  (w.r.t. $|d\z|ds$)  restricted to
the region
\beqn\label{region}
D_t=\Big\{(\z,s): \Big|\z-\frac{s}{t}\x\Big|<\sqrt{8T\lg T}, \,\,   4 \leq s<t-  4, \, \z \in \Om_K\Big\}, 
\eeqn
where  $T=s(t-s)/t$. Observing that
\beqn\label{region5}
\frac{\x-\z}{t-s} =\frac{\x}{t} - \frac{\y}{t-s}\quad \mbox {where}\quad \y:= \z -\frac{s}{t}\x,
\eeqn
 we  obtain  within  $D_t$
\beqn\label{lem4.4_0}
\frac{p_{t-s}(\z-\x -\xi)}{p_{t-s}(\x-\z)} = \exp\bigg[\frac{(\z-\x)\cdot\xi -\frac12|\xi|^2}{t-s} \bigg] = \exp\bigg[- \frac{\x\cdot\xi}{t} + O\bigg(\frac{\sqrt{T\lg T}}{t-s}\bigg)\bigg].
\eeqn


The ratio   $\sqrt{T\lg T}/ (t-s) = (s/t)\sqrt{T^{-1}\lg T}$ is bounded for $ 4< s<t- 4$ and approaches zero as $t-s\to\infty$. For $(\z,s)\in D_t$,   the ratio of two integrands 
 involved  in (4.23) and (4.25), being equal to $\frac{p_{t-s}(\z-\x -\xi)}{p_{t-s}(\x-\z)}e^{\x\cdot\xi/t}$
, therefore is bounded away from zero  and infinity 
   and   approaches unity   uniformly  for   $4<s < s_0$ as $t- s_0\to \infty$.   In view of (\ref{uppbd_int}) and Corollary \ref{cor4.3} we  then  infer  that if $4\leq \a<\b\leq t-4$ and $\b-\a=\de t$,
the contributions from the  interval $\a\leq s\leq \b$  to both of the integrals   in  (\ref{J}) and (\ref{J1})  restricted to $D_t$  are  $[O(\de x) + o(x)] \times p_t(x)$.
  It  also  follows that  the ratio of  the contributions from $s<(1-\de)t$ to  these
 integrals  restricted to $D_t$ approaches unity for each $\de>0$. 
 Taking these into account we  conclude that  the difference of   two integrals restricted to $D_t$ is $o(x p_t(x))$. 
 
 In the sequel
of this proof  we     denote  the right-hand side of  (\ref{J1})  with the outer double integral restricted to a region $D\subset \Om_K\times [0,t]$ by ${\bf F}\{D\}= {\bf F}_{K}^{\x,t}\{D\}$:
$${\bf F}\{D\} :=  \frac1{p_t(x) }\int_{(\z, s)\in D\cap (\Om_K\times [4,t-4])} p_{t-s}(\x-\z) |d\z|ds \int_{ \partial K} e^{-\x\cdot\xi /t} H_K(\z,s; d\xi );$$
similarly denote  by   ${\bf S}\{D\}$  the corresponding integral for  (\ref{J}):
  $${\bf S}\{D\} := \frac1{p_t(x) }\int_{(\z, s) \in D}|d\z|ds \int_{\xi\in \partial K} p_{t-s}(\x-\z +\xi)H_K(\z,s; d\xi).$$
Then  $F^{\x,t}_K(t) = {\bf F}\{\Om_K\times [0, t]\}$ and  $ E_\0[ \,{\rm vol}_d\,(S_K(t))\, | \, B_t= \x] = 
{\bf S}\{\Om_K\times [0, t]\} + O(1)$;  
(\ref{uppbd_int})  is written as
\beqn\label{S[a,b]}
{\bf S}\{\Om_K\times [\a,\b)\} \leq E_{\0}[ {\rm vol}_d(S_K([\a,\b))\,|\, B_t =\x];
\eeqn
and
 by what we have observed right after (\ref{lem4.4_0}) it follows that
\beqn\label{S0} 
{\bf F}\{D_t\cap
(\Om_K\times [\a,\b])\}   \leq C \frac{\b-\a}{t} x + o(x) \qquad (4\leq \a< \b\leq t-4)
\eeqn
and that the difference $ {\bf S}\{D_t\} - {\bf F}\{D_t\} $ is negligible in the sense that
\beqn\label{S1} 
 {\bf S}\{D_t\} - {\bf F}\{D_t\}  = o(x).
\eeqn
Owing to (\ref{S1}) the formula of the lemma   follows   if  the contribution from  the complement of $D_t$  to each of  these two  integrals is   evaluated to be $o(x)$   so as to yield
\beqn\label{S2} 
 {\bf S}\{\Om_K\times [0, t]\} =  {\bf S}\{D_t\}  + o(x),
\eeqn
and
\beqn\label{S3} 
 F^{\x,t}_K(t)  = {\bf F}\{D_t\} + o(x).
\eeqn
We shall verify (\ref{S2})    in Step 2 and (\ref{S3})   in Steps 3 through 5.  In view of what is mentioned at the beginning of this proof these relations  hold if $x<Mt$ for each $M>0$, and hence for the proof we may suppose $v\to \infty$, and we shall do so in Steps 3 through 5.

\v2
{\bf Step 2.}\,
  Put $\hat D_t=(\Om_K\times [0,t)) \setminus D_t$, so that (\ref{S2})
 is written as ${\bf S}\{\hat D_t\} = o(x)$.
 According to Corollary \ref{cor4.-1} 
\beqn\label{426}
{\bf S}\{\hat D_t\}= E_\0[ \,{\rm vol}_d\,(\{\z:  \,    (\z,\sigma_{- K +\z} )\in \hat D_t\} )\,|\,\, B_t=\x].
\eeqn
Plainly  $ \{\z\in \Om_K:  (\z,\sigma_{- K +\z})\in \hat D_t\}$ is included in $ \{\z\in \Om_K: \sigma_{- K +\z} \leq t\} \subset  S_K(t)$ and not empty only under the occurrence of   the event
$${\cal E}:= \{ \exists s\leq t,  (B_s +\xi, s) \in \hat D_t \, \,\,\,\mbox{for some}\, \,\,\,\xi\in K \}.$$
Hence  the conditional expectation in (\ref{426}) is dominated by
$$
E_\0[ \,{\rm vol}_d\,(S_K(t)); {\cal E}\, | \, B_t= \x].
$$
Combined with  (ii) of Lemma \ref{prop6} (applied with $a=1$)  
 an application of Schwarz inequality yields that for each  $\e>0$ the last  conditional   expectation  is at most  ${\rm vol}_{d-1}({\rm pr}_{\be} K)$ times a constant multiple of 
$$ \Big(P_\0\Big[\exists s\in [\e t,(1-\e)t], \;  (B_s +\xi, s) \in \hat D_t \,\,\,\mbox{for some}\,\,\, \xi\in  K \, \Big| \, B_t= \x \Big]\Big)^{1/2} x + \e x.
$$
Here  we have also applied (\ref{S[a,b]}) as well as Corollary \ref{cor4.3}---or rather (\ref{extension})---to obtain  the upper bound $\e x$ for the contribution from the intervals $[0,\e t]$ and $[(1-\e) t,t]$.

We claim that for each $\e>0$ the conditional probability above  tends to zero. For the proof 
we may disregard the dependence on $\xi$ of the event under  $P_\0$.  By scaling property of Brownian motion  we   infer that
 \beq
&& P_\0[\exists s\in [\e t,(1-\e)t], \, (B_s,s)\in \hat D_t \, | \, B_t= \x] \\
 &&= P_\0[ \exists u\in [\e,1-\e], \, |B_u- u\x/\sqrt t\,|\geq \sqrt{8(u(1-u) \lg T} \, | \, B_1 = \x/\sqrt t\, ] \\
 &&= P_\0[ \exists u\in [\e,1-\e],\,  |B_u|\geq \sqrt{8(u(1-u) \lg [tu(1-u)]} \, | \,  B_1 = 0] \\
 &&\to 0
 \eeq
as $t\to\infty$, verifying the claim, which   entails that the conditional expectation in (\ref{426}) is $o(x)$. Thus
(\ref{S2}) is verified as required.

\v2
{\bf Step 3.}\,
We must   verify  (\ref{S3}), which we may write down explicitly as  
\begin{eqnarray}\label{int_J0}
&&F^{\x,t}_K(t) - {\bf F}\{ D_t\} \nonumber \\
&&\quad=\frac1{p_t(x) }\int_ 4^{t- 4}ds\int_{|\z-(s/t)\x|\geq \sqrt{8 T \lg T},\, \z\in \Om_K}  p_{t-s}(\x-\z) |d\z|\int_{ \partial K} e^{-\x\cdot\xi /t} H_K(\z,s; d\xi ) \qquad \nonumber \\
&&\quad= o(x).
\end{eqnarray}
 For simplicity we  suppose that $v\to\infty$ in the rest of the proof as mentioned at the end of Step 1.
For a constant $R \geq 3$, we  take $\sqrt{8T(R  v\vee \lg T)}$ in place of $\sqrt{8T\lg T}$ in the definition of $D_t$ and denote the resulting  region by $D'_x$:
\beq
D'_x=\Big\{(\z,s): \Big|\z-\frac{s}{t}\x\Big|<\sqrt{8T(Rv \vee \lg T)}, \,\,   4 <s<t-  4, \,\z \in \Om_K\Big\},
\eeq
 and its complement in $\Om_K\times [4,t-4]$ by $\hat D'_x$. Plainly  $ D_t \subset  D_x'$, so that  $F^{\x,t}_K(t) - {\bf F}\{D_t\} ={\bf F}\{ \hat D'_x\} + {\bf F}\{D'_x\setminus D_t\}$.

 First we evaluate  ${\bf F}\{ \hat D'_x\}$. 
Using  Propositions \ref{prop6.3} and \ref{lem6.9} and Theorem \ref{thm01}, we see  that $\int_{ \partial K} e^{-\x\cdot\xi /t} H_K(\z,s; d\xi ) \leq e^{ R_Kv}q_K(\z, s) \leq C\gamma_K e^{Rv}p_s(\z)$, and  by using  the identity (\ref{id}) and changing the variable according to  $\y = \z- \frac{s}{t}\x$ (as in (\ref{region5}))  the integral to  be evaluated is at  most 
$$C\ga_Ke^{Rv}p_t(x)\int_ 4^{t- 4}ds\int_{|\y|> \sqrt{8TRv}} p_{T}(\y)|d\y|\leq C'\ga_Kp_t(x)v^{\nu}e^{-3Rv} t,
$$
where $\gamma_K =1$ if $d=2$; $= \cp(K)$ if $d \geq 3$ as in Lemma \ref{lem3.50}.  Thus 
\beqn\label{S4}
{\bf F}\{ \hat D'_x\}  = o(x).
\eeqn
The rest of this proof is devoted to showing ${\bf F}\{D'_x\setminus D_t\}=o(x)$, in which the relation (4.32) (already verified in Step 2) will be used in a
 significant way. Here we write down it   in the following reduced form:
\beqn\label{Stp2}
{\bf S}\{D'_t\setminus D_t\} = o(x).
\eeqn

If $v < \lg (t\vee 2)$  [the choice of $\lg(t\vee 2)$ is rather arbitrary; it may be, e.g., $t^{1/4}$], then the integral defining  ${\bf F}\{D_t'\setminus D_t\}$ restricted to $s\leq t/2$ is negligible owing to  (\ref{Stp2}), since then the $O$ term in (\ref{lem4.4_0}), which is now $O(\sqrt{T(Rv\vee \lg T)}/(t-s))$,  is bounded for $s\leq t/2$ so that we have in effect the same integrand as that of the integral defining ${\bf S}\{D\}$. 
  Therefore we have only to evaluate the integral  in (\ref{int_J0}) over $D_x' \setminus D_t$  in  the situation when
  \beqn\label{AB}
 \mbox{either}\,\, {\rm (i)} \,\, \quad  v\geq \lg (t\vee 2)  \quad\,\, \mbox{or} \,\,\,  {\rm (ii)}   \quad    s> t/2 \,\, \mbox{and}\, \, v < \lg (t\vee 2). 
  \eeqn

Putting ${\bf v} =\x/t$ and
$$\fa(\xi,u) = \int_{\partial K} e^{- {\bf v}\cdot\xi_1 }H_K(\xi, u; d\xi_1),$$
we write the inner integral  in (\ref{int_J0}) in the form
\beqn\label{M-fa}
\int_{ \partial K} e^{-\x\cdot\xi_1 /t} H_K(\z,s; d\xi_1 ) =\int_0^sdu\int_{\partial U(R)} \fa(\xi,u)H_{U(R)}(\z, s-u;d\xi).
\eeqn
Denote  the last repeated integral restricted to 
$[\a,\b]\times \partial U(R)$ by
$$I_{[\a,\b]} = I_{[\a,\b]}(s,\z, {\bf v}; R) \qquad  (0\leq \a<\b\leq s).$$ 
 The contribution   coming from the small interval $[0, M/v]$ with a constant $M>2R$  is dominant and problematic; we make  its evaluation  in the succeeding  two steps.  In the rest of the present step we  ascertain that the other part is negligible.  

 Let  $(\z,s) \in D_x'\setminus D_t$.  Then under (\ref{AB})
 \beqn\label{restr1}
\Big|\frac{\z}{s} - \frac{\x}{t}\Big|\leq \sqrt{8(Rv\vee \lg T) (t-s)/st} \leq C[\sqrt {v/s} \vee \sqrt{(\lg t)/t}\,] =O(\sqrt v\,),
 \eeqn
hence $z/s\sim v$ (where $z=|\z|$) and    we deduce that for $v$ large enough,  
$$
I_{[ M/v, s]}\leq e^{v}\int_{ M/v}^s q(z,s-u; R)du 
\leq C e^{2Rv}v^{-\eta}  p_{s- M/v}(z)
$$
with $\eta = \nu +\frac32>0$. Here we have employed Theorem \ref{thm2.2}  to evaluate  the integral in the middle member, which equals $ \int_0^{s-  M/v} q(z,u; R)du$. Observing
 $z^2/(s- M/v) -z^2/s >  M z^2/s^2v >\frac12 v M$,
  we then find that  for $t$ large enough, 
$$
I_{[ M/v, s]}
\leq C' p_s(z) e^{2Rv}e^{- \frac14 v M},
$$
from which we  infer that if $ M \geq 9 R$, then 
\beqn\label{int_J}
\frac1{p_t(x) }\int_{D_x'\setminus D_t}
  p_{t-s}(\x-\z)\, I_{[ M/v, s]}(s,\z,{\bf v}; R) |d\z|ds = o(x). 
\eeqn
\v2
{\bf Step 4.}\, It remains to verify that  
 for some suitably chosen constant  $R$ 
\beqn\label{S4_1}
\frac1{p_t(x) }\int_{D_x'\setminus D_t}
  p_{t-s}(\x-\z) \,I_{[0,M/v]}(s,\z,{\bf v};R) |d\z|ds = o(x).
\eeqn
 To this end  we are going to apply  Theorem \ref{cor4.1}. We take $M=9R$ for definiteness. 
 
   Look at  the  repeated integral  on the right-hand side of  (\ref{M-fa}) and let  $\th=\th_{\xi,\x}\in [0,\pi]$ be the colatitude  of  $\xi\in \partial U(R)$ relative to   $\x$ taken as a north-pole so that  $\cos \th =\x\cdot \xi/xR$.
 Given a small positive number $\e$, we then   break  $I_{[0, M/v]}$   into two parts by splitting   $\partial U(R)$, the range  of the  inner integral of the corresponding  repeated integral  for  it, along the   parallel of colatitude   $\th =\cos^{-1}\e$ and denote the part  for   $\cos \th >\e$ by  $I_{[0, M/v]}^{ >\e}$ and the other by $I_{[0, M/v]}^{ \leq\e}$: 
 $$I_{[0, M/v]}^{ >\e} = \int_0^{M/v} du\int_{\partial U(R)} \fa(\xi,u)\1(\cos \th >\e)H_{U(R)}(\z, s-u;d\xi)$$
 and similarly for  $I_{[0, M/v]}^{ \leq\e}$. In the rest of this step we prove that the part $I_{[0, M/v]}^{ \leq\e}$ makes only a negligible contribution to (\ref{S4_1}) in comparison with   the other part $I_{[0, M/v]}^{ >\e}$---the latter will be treated in the next step. For the proof,  on putting
$$Z(\xi) = \int_0^{M/v} p_{s-u}(z)\fa(\xi,u)du,$$
 it suffices to show that uniformly for ${(\z,s)\in D_x'}$ subject to  the condition (\ref{AB}), 
 \beqn\label{S65}
  \frac{\int_{\cos \th\leq \e} Z(\xi) H_{U(R)}(\z,s; d\xi)}{\int_{\partial U(R)}  Z(\xi) H_{U(R)}(\z,s; d\xi)} = O( e^{-v}).
 \eeqn

As in (\ref{region5}) we let $\y=\z-\frac{s}{t}\x$ so that   $\z = s{\bf v} + \y$.  Let $(\z,s)\in D_x'$ and (\ref{AB}) be satisfied. On multiplying both sides of   (\ref{restr1}) by  $1/v$,  we infer that  
\beqn\label{St4-0}
|\z/z -\x/x| =O(1/\sqrt v\,),
\eeqn
thus $|\z\cdot \xi/zR -\cos \th| =O(1/\sqrt v)$ ($z=|\z|$). 
   In view of  (\ref{restr1})  $z/s \sim v$ (as $v\to\infty$), hence
 if  $0\leq u \leq  M/v$ in addition,  then  $|\z/(s-u) - \z/s| =  O(1/s)$ and we find
\beqn\label{z/x}
\frac{\z}{s-u}  = \frac{\x}{t} +\frac{\y}{s} + O(1/s); 
\eeqn 
Note that (\ref{z/x}) is valid in $D_t$  (without assuming (\ref{AB})), for   in $D_t$, $Rv$ on the left-hand side of  (\ref{restr1})   may be dropped so that we have   the  better bound   $O(s^{-1}\lg s)$  instead of $O(\sqrt v)$, although this  do not improve the bound $O(1/s)$ in (\ref{z/x}).

With the help of (\ref{St4-0})  we apply   Theorem \ref{cor4.1} (ii) (or rather (\ref{xx})) to see that for each $\e>0$,
\beqn\label{St4-1}
\frac{H_{U(R)}(\z,s; \{\cos \th <\e\})}{H_{U(R)}(\z,s; \partial U(R))} \leq Ce^{-(1-\e)Rv +O(\sqrt v)}.
 \eeqn
 where the  factor $v^{\nu-\frac12}$ that  arises by the application   is absorbed into $e^{O(\sqrt v\,)}$.
Obviously $Z(\xi) \geq e^{-v} \int_0^{M/v} p_{s-u}(z) q_K(\xi,u)du$, of which  we perform integration by parts for the integral on the right. For $v$ large enough   we have
$-\partial_u p_{s-u}(z) \geq \frac13 v^2p_{s-u}(z)$, so that  
\[
Z(\xi) 
\geq  e^{-v}p_{s- M/v}(z)P_\xi[\sigma_K< M/v] 
+ \frac13 v^2e^{-v} \int_0^{M/v} p_{s-u}(z) P_\xi[\sigma_K< u]du. 
 \]
Plainly $p_{s-u}(z) = p_s(z)e^{-v^2u(1+o(1))/2}\{1+o(1)\}$.  On the other hand, from  Proposition \ref{propA6} (applied with $(u, R)$ in place of $(t,x)$) we  infer that if $R \geq 10$ and  $u<1/10$,
$$P_\xi[\sigma_K< u] \geq \k_d\la(K)p_u(R)e^{-2R/u}\qquad (\xi\in \partial U(R))$$
($\la(A)$ appears in Proposition \ref{propA6} and is given at the beginning of Section 6.4).
By the relation 
$\int_0^{(1+h)\eta/v} e^{-\eta^2/2u - v^2u/2}u^{-\nu-1} du \sim 2(v/\eta)^\nu K_p(v\eta)\sim (v/\eta)^\nu\sqrt{{2\pi}/{v\eta} } \, e^{-v\eta}$  as   $v \to\infty$  
valid for each real numbers  $p$, $\eta>0$ and $h>0$ (\cite{E}, p.146) an easy computation shows that for any $\de>0$ small enough,
$$Z(\xi) \geq \k_d'\la(K) R^{-2\nu}e^{-v}p_s(z)\exp\{- vR\sqrt{1+ 4R^{-1}}(1+\de)\}$$
for  $R>1$ and all sufficiently large $v$. Similarly,  using Lemma \ref{lemA60}  with the help of  Lemma \ref{thm5},
we deduce that $P_\xi[\sigma_K< u] \leq \k_d\la(K)u^{(d-3)/2}p_u(R)e^{3R/2u} \leq \k'_d\la(K)p_u(R)e^{2R/u}$, and then, as above, that 
$$Z(\xi) \leq \k_d'' R^{-2\nu} e^{v}p_s(z)
\exp\{- vR\sqrt{1- 4R^{-1}} (1-\de)\} \quad (R>5). $$

Now let $\e=1/10$ in (\ref{St4-1}) and $\de = 1/20$. Then taking  a large $R$ ($R=20$ suffices) we see that  $$Z(\xi')/Z(\xi) \leq
\k_d''' e^{2v+vR(\sqrt{1+4R^{-1}}-\sqrt{1- 4R^{-1}}\,) +2\de vR} \leq \k_d''' e^{10v + \frac1{10}Rv} \quad\mbox{   for all $\xi, \xi' \in \partial U(R)$,}
$$
and that combining this with  (\ref{St4-1})  leads to (\ref{S65}).

 In below we fix  the constant $R$ as chosen above. 

\v2
{\bf Step 5.}\, In order to finish the proof it now suffices to show that 
\beqn\label{int_J6}
\frac1{p_t(x) }\int_{D_x'\setminus D_t}
  p_{t-s}(\x-\z) I^{> 1/10}_{[0, M/v]}(s,\z;{\bf v};R) |d\z|ds =  o(x).
\eeqn
Recall that this integral is obtained by reducing the range of  $(\z,s)$ from that appearing in  ${\bf F}(D_x'\setminus D_t)$, and
we may further reduce the range   by  (\ref{AB})---($\z,s)$ is related to $v$ within  $D_x'\setminus D_t$.

 Let $(\z,s)\in D_x'$, $u<M/v$ and (\ref{AB}) be satisfied. Then (\ref{z/x}) is in force,  by which    we have $\z\cdot \xi/(s-u) = Rv\cos \th + \y\cdot\xi/s  +O(1/s)$,
hence   Theorem \ref{cor4.1} (i) entails  that 
\beqn\label{H4.0}
H_{U(R)}(\z, s-u;d\xi) = 
e^{\y\cdot\xi/s }p_{s-u}(z)  V(d\xi) \Big[1+  O(v^{-1/2})+O(s^{-1})\Big]
\quad \mbox{if} \;\; \cos\th > 1/10,
\eeqn
where  $V(d\xi) = [\om_{d-1}R^{2\nu +1}v] e^{R v\cos \th}\cos \th \,  m_R(d\xi)$,
so that   the integral on the left-hand side of (\ref{int_J6})
   is  dominated by a constant multiple of
\beqn\label{S4_2}
\int_{D_x'\setminus D_t}
  p_{t-s}(\x-\z) d\z ds\int_0^{ M/v}p_{s-u}(z)du\int_{\partial U(R) } e^{\y\cdot\xi/s }\fa(\xi,u)  \1(\cos \th >  {\textstyle \frac1{10}} ) V(d\xi).
\eeqn
  If we replace the range $D_x' \setminus D_t$  by $D_t$ in this expression, the contribution to the resulting integral from the error term $O(s^{-1})$  on the right-hand side of (\ref{H4.0}) is at most $o(xp_t(x))$ since
the contribution from the interval $4<s< \de t$ is $O(\de xp_t(x))$ for each $\de$, so  that this integral  agrees with that in (\ref{int_J6})
with  $D_x' \setminus D_t$ replaced by $D_t$ apart from an error  term  of magnitude   
$o(xp_t(x))$.

What is important for the argument made below is the fact that $V$ is independent of $(\z, s-u)$ since the following computation concerns only the ratio of 
the integral of   
  $$  p_{t-s}(\x-\z)  p_{s-u}(z)e^{\y\cdot\xi/s }   $$
over $D_x'\setminus D_t$ to that over $D_t$  (for each  $4<s<t-4$ and $u<M/v$).
It is convenient to  take  $\y= \z-\frac{s}{t}\x$ rather than  $\z$ as the variable of integration, and   the  ranges of integration accordingly become
  $\sqrt{8T\lg T}< y=|\y|\leq \sqrt{8RTv}$ 
and   $y\leq\sqrt{8T\lg T}$, respectively. 
Put  $T' = (s-u)(t-s)/(t-u)$ so that 
$$p_{t-s}(\x-\z) p_{s-u}(\z) = p_{t-u}(\x) p_{T'}\Big(\z -\frac{s-u}{t-u}\x\Big).$$
Then, observing that  $T' \sim T$,  $s> T'$ ($s> 4, u< M/v$)  and
$$\z -\frac{s-u}{t-u}\x = \y + \frac{(t-s)u}{(t-u)t}\x = \y +  b(u)$$
with $|b(u)| \leq M$ for $u\leq M/v$, we infer that the ratio of the integral 
$$\int_{\sqrt{8T\lg T}< y \leq \sqrt{8RTv}} e^{\xi\cdot \y/s}p_{T'}(\y+ b(u))|d\y|$$
to the same integral but  over  $y\leq \sqrt{8T\lg T}$  is bounded and  tends to zero as   $s\wedge (t-s) \to\infty$ (so that $T\to\infty$). Now we take account of   the repeated integral   (\ref{S4_2}) as well as  the corresponding one for the integral with $D_t$ replacing $D_x'\setminus D_t$, of which the latter admits the bound  (\ref{S0}).
By making comparison between them  
 we  infer that for each $\de>0$, 
the integral in (\ref{int_J6})   restricted to $\Om_K\times (\de t,(1-\de) t)$ is $p_t(x)\times o(x)$  on the one hand  and that restricted to $\Om_K\times ([4,\de t] \cup [(1-\de)t,t-4])$ is dominated by a constant multiple of $\de x\times p_t(x)$  on the other hand. Since $\de$ may be made arbitrarily small this verifies  (\ref{int_J6}).
  The proof of Lemma \ref{lem0a} is complete. 
\qed

\begin{Lem}\label{lem4.52} \, For each $\e \in (0,1)$  there exists a constant  $M$ such that if   $x^2/t >M$ and $t>M$, 
 \beqn\label{lem-eq4.6}
  \int_{W_s} q_K(\z,s)p^{(d)}_{t-s}(\z-\x)|d\z| \leq  C\gamma_K^*(t) p^{(d)}_t(x)e^{- \e x^2/5t}  \quad \mbox{for}  \quad s\in [\e t,t],  
 \eeqn 
 where $W_s := \{\z\in \Om_K: | \z\cdot \x|/x < \frac13 sx/t\}$, $\gamma_K^*(t)= 1/\lg t$ or $\cp(K)$ according as $d=2$ or $\geq 3$ 
and $C$ is a universal constant.
\end{Lem}
\v2\n 
\pf \, We apply  Propositions \ref{prop6.3} ($d\geq 3$) and \ref{lem6.9} ($d=2$)  of Section 6  along with Theorem \ref{thm01} for $d=2$   to see that if $R=2 R_A$,
\beqn\label{H-q-p}
 q_K(\z,s) \leq C\gamma_KR^{2\nu} [(Rz/s)\vee 1]^{\frac12-\nu}p^{(d)}_s(\z)e^{Rz/s}
\eeqn
($z=|\z|$).    Writing $z_1 =\z\cdot\be$ and  $y= |\z-z_1\be|$  ($\be= \x/x$) and using the identity (\ref{id})  we obtain
$$p^{(d)}_s(\z)p^{(d)}_{t-s}(\z-\x) =p^{(d)}_t(x)p_T^{(d-1)}(y)p_T^{(1)}\Big(z_1- \frac{s}{t}x\Big), \quad T=s(t-s)/t.$$
For simplicity let $R_K =1$.  
On substituting the obvious bounds  
$${\textstyle |z_1- \frac{s}{t}x|\geq \frac23  sx/t\quad\mbox{ and}\quad  z\leq y + \frac13 sx/t}$$
valid  on $W_s$,  using $\int_\eta^\infty p_T^{(1)}(y)dy \leq \eta p_T^{(1)}(\eta)$ if $\eta \geq \sqrt T$  and letting $e^{z/s}$ dominate $[(2z/s)\vee 1]^{\frac12-\nu}$,   the right-hand side of (\ref{lem-eq4.6}) is at most
 \beqn\label{C}
 C''\gamma_Ke^{x/t} p^{(d)}_t(x)\Big[\frac{2sx}{3t}p_T^{(1)}\Big(\frac{2sx}{3t}\Big)\Big]\int_{\R^{d-1}} p_T^{(d-1)}(\y)e^{3|\y|/s}|d\y|,
 \eeqn 
 provided that  $x^2/t >M$ and $M$ is large enough.
 An elementary computation shows that  the integral above is at most  $Ce^{9T/2s^2}[1+(T/s^2)^\nu]$, which is bounded by a universal constant for $s>1$. On the other hand, since $\a e^{-\a^2/2}$ is decreasing for $\a>1$ and $s/t\sqrt{T} \geq  \sqrt{s}/t$,  the quantity in the square brackets   is less than   
 $${\textstyle \sqrt{\frac{x^2 s}{t(t-s)}}\exp\Big\{-\frac{2x^2s}{9t(t-s)}\Big\} \leq c\sqrt{\e x^2/ t}\, e^{-2\e x^2/9t}}\quad\mbox{for}\quad \e t\leq s\leq t,$$
if $ 2\sqrt{\e}  x/3\sqrt t>1$.  It is easy to see that if $d\geq 3$,   (\ref{C})  is  dominated by
 the right-hand side of (\ref{lem-eq4.6}), provided  $M$ is
taken so large  that $x>20/\e$. This shows the lemma for  $d\geq 3$.  

For  $d=2$ the factor $\ga_K$ needs to be replaced by $1/\lg t$, which may be ascertained  by slightly modifying the proof above with  simple remarks. Indeed,  if  $x\geq t^{3/4}$, then  (\ref{lem-eq4.6}) is obtained in the  argument made above, since the factor $1/\lg t$ is  blotted out by an exponential factor $e^{-\e x^2/5t}$ of which the number  $5$  could be replaced by a little larger one.  
   In the  case when  $x< t^{3/4}$,  if $y = |\z-z_1 \be|\leq t^{3/4}$ (and $|z_1|< sx/3t$, $s>\e t$) so that
   $z= O(s^{3/4})$,  we have $q_K(\z,s) \leq Cp_s(z)/\lg t$ instead of (\ref{H-q-p}),  hence conclude the bound (\ref{lem-eq4.6}) since the contribution from  $y>t^{3/4}$ is negligible.
\qed

\v2\n
\begin{Lem}\label{lem0b} \,  {\rm (i)}  For    any Borel  set  $E\subset \partial K$ and $0<s<t$,
\beqn\label{eqn3}
H_K(\x,t;E) \leq \int_{\z\in \Om_K} H_K(\z,s;  E)p^{(d)}_{t-s}(\z-\x)|d\z|.
\eeqn

 {\rm (ii)}  For each $\e \in (0,1)$  there exists a constant  $M$ such that if   $x^2/t >M$ and $t>M$, then  for all  $s\in [\e t,t]$  and $E \subset \partial K$,
\begin{eqnarray}\label{eq4.5}
H_K(\x,t; E) &\geq&  (1-\e) \int_{\z\in \Om_K} H_K(\z,s; E)p^{(d)}_{t-s}(\z-\x)|d\z| \nonumber \\
&&\,\,  - \, C\gamma_K^*(t) p^{(d)}_t(x)e^{- \e x^2/5t}, 
\end{eqnarray}
 where $\gamma_K^*(t) $ is the same function as defined in Lemma \ref{lem4.52} and $C$ is a universal constant.
\end{Lem}
\v2\n 
\pf \,  Let $0<s<t$.   The upper bound (\ref{eqn3}) follows from the identity
$$H_K(\x, t;E) = \int_{\z\in \Om_K}H_K(\z, s; E)P_\x[B_{t-s}\in d\z, \sigma_K >t-s].$$

For the proof of the lower bound  (\ref{eq4.5}), let $s\geq\e t$ and put $W_s := \{\z: | \z\cdot \be| < \frac13 sx/t\}$, $\be=\x/x$ as in  Lemma \ref{lem4.52}.  First we observe that   the second term on the right-hand side of (\ref{eq4.5}) is an upper bound of the contribution from  $W_s$  to the integral  in the first term.   Indeed, this follows by simply  substituting  $\partial K$ for $E$   and then applying Lemmas \ref{lem4.52}. 
  
To complete the proof it  suffices to show that as $t\to\infty$ and $x/\sqrt {t} \to\infty$
\beqn\label{ratio1}
\frac{\int_{\z\notin W_s} H_K(\z,s; E)P_\x[B_{t-s}\in d\z, \sigma_K <t-s]}{\int_{\z\in \Om_K} H_K(\z, s; E )p_{t-s}(\z-\x)|d\z| }\,\longrightarrow \, 0
\eeqn
uniformly for  $E\subset \partial K$ with $H_K^\infty(E)>0$ as well as for $s\in [\e t,t]$.  Since
$$P_\x[B_{t-s}\in d\z, \sigma_K <t-s]/|d\z| = P_\x[\sigma_K <t-s | B_{t-s} =\z]p_{t-s}(\z-\x),$$
  we   have only to show that uniformly for $\z\notin W_{s}$ and $s\in [\e t,t]$,
\beqn\label{frac2}
P_\x[\sigma_K <t-s | B_{t-s} =\z] \to 0.
\eeqn
On expressing the Brownian bridge by means of a free Brownian motion and  reversing the time 
 the probability in (\ref{frac2}) is expressed as  $P_{\z}[B_u + \frac{u}{t-s}(\x -B_{t-s})\in  K\,$ for some $u\in(0,t-s) ]$, and hence,   on scaling  the space and time variables,    dominated by
 $$\sup_{s\in [\e t,t]} \,\sup_{\y: \y\cdot\be >sx/3t\sqrt{t-s}} P_\y\bigg[ B_u -uB_1 + \frac{u\x}{\sqrt{t-s} }\in \frac1{\sqrt{t-s} }K\quad\,\mbox{for some}  \,\,\, u\in (0,1)\bigg]. $$
Note that  $B_u -uB_1 $  under  $P_\y$ has the same law as  $B_u-uB_1 +(1-u)\y$ under  $P_\0$.  Now suppose $x/\sqrt t\to\infty$ and $t\to\infty$ and first consider the case when $s< t-1$. Then  
since  $sx/3t\sqrt{t-s} \to \infty$ so that $ | (1-u)\y + u\x/\sqrt{t-s}\,| \to \infty$,  the supremum above obviously converges to zero. The case when $s\geq t-1$ is  also  easy to  dispose of. \qed

\begin{Lem}\label{lem4.51} \, There exists a  universal constant $C$ such that for $ R_K x < t$,
$$ \int_{\z\in \Om_K}q_K(\z,s)p^{(d)}_{t-s}(\x-\z) |d\z| \leq  \left\{ \begin{array} {ll} C \ga_K p^{(d)}_t(x) \quad &\mbox{for}  \quad R_K^2\leq s \leq t/2 \quad (d\geq 2), \\[2mm]
{\displaystyle  \frac{Cp^{(d)}_t(x)}{\lg(t/R_K x)}}   \quad &\mbox{for}  \quad R_K\sqrt t \leq  s \leq  t/2 \quad \mbox{if} \quad  d=2.
\end{array} \right.
$$
Here  $\gamma_K =1$ if $d=2$ and $ \cp(K)$ if $d \geq 3$ (as before).
\end{Lem}

\n
\pf\,  The proof rests on Propositions \ref{prop6.3} and \ref{lem6.9} as in that of Lemma  \ref{lem4.52}.   Let  $z=|\z|$ and $\be'=\z/z$,  
and   we make substitution from  $q_K(\z,s)\leq \ga_Kp_s(z) e^{2R_Kz/s}$. Let $R=2R_K$.   Recalling the identity (\ref{id})  we put $T=s(t-s)/t$. Then, 
observing the identity
$$p_T(\z- {\textstyle \frac{s}{t}}\x)e^{Rz/T} = p_T(\z- {\textstyle \frac{s}{t}}\x -  R\be')e^{R\x\cdot \be'/(t-s)}e^{R^2/T}$$
and  the inequality $1/s<1/T$, we deduce that if $R_K^2 <s <t/2$ and $R_K x<t$, 
\beqn\label{461}
p_{t-s}(\x-\z)q_K(\z,s) \leq C\ga_K p_t(x) p_T(\z- {\textstyle \frac{s}{t}}\x - R\be'),
\eeqn
which shows the first inequality of the lemma.

For the proof of the second one,  let $d=2$ and $R_K=1$,  and put $W_s :=\{\z \in \Om_K: |\z-\frac{s}{t} \x| \leq \sqrt{4T\lg T}\}$. From (\ref{461}) it follows that the integral over the outside of  $W_s$ is negligible. For the evaluation of the integral inside $W_s$ we consider the following  two subcases of   $\sqrt t\leq s\leq t/2$:
$${\rm (i)}\quad {\textstyle \frac{s}{t}} x  < \sqrt{s} \,\lg t\quad\mbox{or}\quad{\rm (ii)}\quad  {\textstyle \frac{s}{t}} x \geq \sqrt{s} \,\lg t$$
(or, equivalently,   $s$  is less than  $[(t/x)\lg t]^2$ or not).
Note that $s/2\leq  T<s$ and   if  (i) is the case, then for $\z\in W_s$, $z = O(\sqrt s \lg s)$  so that $q_K(\z,s) \leq Cp_s(z)/\lg t$ since  $s\geq \sqrt t$. On the other hand,  for $\z \in W_s$,  we have  $z/s = x/t + O(\sqrt{s^{-1}\lg t})$ and,   in the case (ii),   $z/s \sim x/t$, so that $q_K(\z,s) \leq Cp_s(z)/\lg (t/x)$. Hence, in  either case,  $q_K(\z,s) \leq Cp_s(z)/\lg (t/x)$  in $W_s$, so that  for $\z\in W_s$, we may multiply  
 the right-hand side of (\ref{461})  by $1/\lg (t/x)$ and thus find  the second inequality of the lemma.
\qed

\v2\n
{\it Proof  of  Proposition \ref{prop5}.}\, Let $R_A=1$.  In view of Lemma \ref{lem0a} it suffices to prove the required relation (\ref{eqn1}) with $F_K^{\x,t}(t)$ replacing $E_\0[ {\rm vol}_d(S_{K}(t))|B_t=\x]$.
 Picking $0<\e< 1/2$,  we perform the integration by $\z$ first in  the triple integral defining $F_K^{\x,t}(t)- F_K^{\x,t}(\e t+4)$ and then use  the  inequality (\ref{eq4.5}) to   find  that 
\beq
F_K^{\x,t}(t) - F_K^{\x,t}(\e t +4)&=&
\frac1{p_t(x)}\int_{\z\in \Om_K}|d\z|\int_{\e t}^{t-4} p_{t-s}(\x-\z) ds\int_{ \partial K} e^{-\xi\cdot\x /t} H_K(\z,s; d\xi )\\
&\leq& \frac{t}{p_t(x)} \int_{\xi\in \partial K} e^{-\xi\cdot\x/t} H_K(\x,t;d\xi)  
+ \,  Ct \gamma^*_K (t) e^{-\e x^2/5t} \sup_{\xi\in K}  e^{-\x\cdot\xi/t}.
\eeq
In the right most member the second term is  $o(\zeta(x,t))$ as  $x^2/t \to \infty$ and $t\to\infty$.
On the other hand it holds that
$$F_K^{\x,t}(\e t +4)\leq c_K \zeta(x,t)\e,$$ 
as is  deduced immediately from  
 Lemma \ref{lem0a} along with (\ref{extension})  if $R_Ax >t$ and by integrating the inequality of Lemma \ref{lem4.51} over $4< s <\de t$ if $R_Ax \leq t$.  These together yield the upper bound  of $F_K^{\x,t}(t)$ for the asserted equality since $\e$ may be chosen arbitrarily small. The lower bound is verified immediately by   Lemma \ref{lem0b} (i). The proof   is complete.  \qed

\v2
{{\bf 4.2.4.} \sc Proofs of Theorems \ref{thm3.4} and  \ref{thm4.3}.} 
Remember  the measures  $m_{K,\be}(d\xi)$ given in (\ref{mm}) and 
$$\qquad \mu_{t,\x}^K(d\xi)= \frac{e^{-v\be\cdot\xi} H_K(\x, t;  d\xi)}{vp_t(x)}\qquad (\be=\x/x, v=x/t)
$$
(see Section {\bf 3.2}). According to Proposition \ref{prop5},  for any compact set $K$, as $x/t\to\infty$ and $t\to\infty$
\beqn\label{saus/hd}
E_\0[ {\rm vol}_d(S_{K}(t))|B_t=\x]/x = \mu^K_{t,\x}(\partial K) +o(1).
\eeqn

\begin{Lem} \label{lem4.6}   If  $E\subset \partial K$ is compact,  then $\mu_{t,\x}^K ( E) $  is asymptotically dominated by  ${\rm vol}_{d-1}({\rm pr}_\be E)$ in the sense that 
\beqn\label{501}
\limsup_ {x/t\to \infty,\, t\to\infty} \mu_{t,\x}^K ( E)  \leq {\rm vol}_{d-1}({\rm pr}_\be E).
\eeqn 
In particular,
\beqn\label{502}
\limsup_ {x/t\to \infty,\, t\to\infty} \mu_{t,\x}^K (\partial K)  \leq m_{K,\be}(\lan K\ran_\be).
\eeqn 
\end{Lem}
\n
\pf\, 
From the inclusion $E\subset \partial K$ it follows that 
$H_{K}(\x,t;d\xi) \leq H_{E}(\x,t;  d\xi)$ for $d\xi \subset E$ since  Brownian paths that
hits $\partial K\setminus {E}$ before $E$  may contribute to $H_{E}$ but never to $H_{\partial K}|_{E}$.  Hence
$$ \mu_{t,\x}^K(E) \leq \mu^{E}_{t, \x}( E).
$$
If  $E$ is compact, then  by (\ref{saus/hd})  $\mu^{E}_{t, \x}( E)$ is asymptotically dominated by 
$E[ {\rm vol}_d(S_{E}(t))|B_t=\x]/x, $
which in turn  is asymptotically dominated by ${\rm vol}_{d-1}({\rm pr}_\be E)$ owing to Corollary  \ref{cor4.3}. 
 \qed

\begin{Lem}\label{lem4.7}   
For any Borel set  $E\subset \lan K\ran_\be$,  if  ${\rm vol}_{d-1}(\partial_{d-1}({\rm pr}_\be E))=0$, then
\beqn\label{503}
\limsup_ {x/t\to \infty,\, t\to\infty} \mu_{t,\x}^K ( E)  \leq m_{K,\be}(E).
\eeqn
\end{Lem}
\n
\pf\, Let $\overline{ E}$ denote the closure of $E\subset \lan K\ran_\be$ in $\R^d$.
 Then, 
 noting that  ${\rm pr}_\be\overline{E}$ is compact, we have 
 $$ {\rm pr}_\be\overline{ E} = \overline{ {\rm pr}_\be E}.$$
 Hence, from the assumed condition  on $E$  in the lemma it follows that 
  $ m_{K,\be}(E) = {\rm vol}_{d-1}({\rm pr}_\be E) =  {\rm vol}_{d-1}({\rm pr}_\be \overline{E}),$ and  applying Lemma \ref{lem4.6} with $\overline{E}$ in place of $E$ we obtain (\ref{503}). \qed
 \v2
  
  The next lemma  can be shown either directly  or as in  the proof of Lemma \ref{lem4.6}  (with the help of  Proposition \ref{prop5} and Corollary \ref{cor4.2}).
\begin{Lem}\label{cor4.4}  \,Suppose that   $\lan K\ran_{\be} $ lies on a plane perpendicular to  $\be$. Then  as $x/t \to\infty$ and $t\to\infty$,
$\mu_{t,\x}^K  \Longrightarrow m_{K,\be},$ in which the two sides of $\lan K\ran_{\be}$ are distinguished  and the measure $m_{K,\be}$ is  considered to be concentrated in the $+\be$ side of $\lan K\ran_{\be}$.
\end{Lem}

\v2\n
{\it Proof of Theorem \ref{thm3.4}.} \,  For the proof of (\ref{eq3.5})  it  suffices to show  that  
\beqn\label{W5}
\liminf_ {x/t\to \infty,\, t\to\infty} \mu_{t,\x}^K(\lan K\ran_\be)  \geq m_{K,\be}(\lan K\ran_\be).
\eeqn
Indeed,  this together with (\ref{502}) entails    $\lim \mu_{t,\x}^K(\lan K\ran_\be) = m_{K,\be}(\lan K\ran_\be)$.
Hence the first condition for (\ref{eq3.5}), i.e.,  $\lim \mu_{t,\x}^K(\partial K \setminus \lan K\ran_\be) =0$, follows from   Proposition \ref{prop5} (see (\ref{saus/hd})) and Corollary \ref{cor4.3} combined,  and the other condition  (\ref{K6})  from   Lemma  \ref{lem4.7} and (\ref{W5}) since  
$$\partial_{d-1}[ {\rm pr}_\be (\lan K\ran_\be \setminus E)] \subset \partial_{d-1}( {\rm pr}_\be E) \cup \partial_{d-1}( {\rm pr}_\be \lan K\ran_\be) \quad\mbox{for}\quad  E\subset \lan K\ran_\be,$$
 hence (\ref{503}) is valid for $\lan K\ran_\be \setminus E$ in place of $E$
owing to the assumption (\ref{cond1}).
 
The proof of (\ref{W5}) will be carried out  by  employing   Lemma \ref{prop6} and  to this end we  bring in an auxiliary set, $W$ say,  that contains  $K$.  
 
Let   $V_\b =\{ \z\in \R^d: \z\cdot \be \geq -\b\}$ for  $\b>0$  chosen  so that  $K\subset  V_\b $,  and    put   
$$L_\xi = \{\xi- s\be:  s\geq 0\} \cap  V_\b  \quad\mbox{for}\quad \xi\in \lan K\ran_{\be}   \quad\mbox{and}\quad  W =\cup_{\xi\in \lan K\ran_{\be}} L_\xi. $$ 
Note that  $K\subset W$, $\lan W\ran_\be =\lan K\ran_\be$ and $W$ is  compact.  Let 
$$D= W\cap \partial V_\b \quad \mbox{ and}\quad C=\partial W\setminus (\lan K\ran_{\be} \cup D).$$ Then 
$C,D$ and  $\lan K\ran_{\be}$ together make up the decomposition of $\partial W$. 

Now consider the limit procedure as $t\to\infty, x/t\to \infty$.   Due to the inclusion   $K\subset W$ we have the inequality $\mu_{t,\x}^W(d\xi)  \leq \mu_{t,\x}^K(d\xi)$ for $d\xi \subset \lan K\ran_\be$. 
For the proof of (\ref{W5})  it therefore suffices to show that 
\beqn\label{W}
 \liminf \mu_{t,\x}^W(\lan K\ran_\be)  \geq m_{W,\be}(\lan K\ran_\be).
\eeqn
 Corollary \ref{cor4.2} entails 
$
\liminf E_\0[ S_W(t) \,|\, B_t=\x]/x \geq m_{W,\be}(\lan K\ran_\be),
$
and writing  the expectation on the left  as in  (\ref{saus/hd})  (according to Proposition \ref{prop5}) we   find 
\beqn\label{W2}
\liminf \mu_{t,\x}^W(\partial W)   \geq  m_{W,\be}(\lan K\ran_\be).
\eeqn
Since
$$\pr_\be C = \mbox{ dis-ct}_\be (K),$$
the assumption  (\ref{cond1})  implies  ${\rm vol}_{d-1}(\overline{{\rm pr}_{\be}C}\,)=0$; hence $E_\0[ S_C(t)| B_t=\x] =o(x)$ owing to
Lemma \ref{prop6}, which in turn shows that $\mu_{t,\x}^W(C) \leq \mu_{t,\x}^C(C) \to 0$.   
 By Lemma \ref{cor4.4} $\mu_{t,\x}^W(D) \to 0$. Combined with (\ref{W2}) these together show (\ref{W}), as desired.  
 \qed

\v2\n
{\it Proof of  Theorem \ref{thm4.3}.}  This is now clear from  Theorem \ref{thm3.4} and Proposition \ref{prop5}. \qed

\section{ Brownian motion with a constant drift }

The law of a Brownian motion with a  constant drift is absolutely continuous relative to the law of a standard Brownian motion with a  Radon-Nikodym derivative    of a simple form and the results obtained so far is  translated to those for the Brownian motion with a  constant drift   as exhibited below. 
 The Brownian bridge  $P_{{\bf 0}}[\,\cdot \, |B_t=  x{\bf e}] $ with $v:=x/t$ kept away from zero being comparable or similar to  the process $B_s - sv{\bf e},$ $0\leq s\leq t$, we in particular derive   from the results  for the bridge in (vi), (vii) and (viii) of Section 1 the corresponding ones,  (vi$'$), (vii$'$) and (viii$'$)  say, for the latter.  

We fix a unit vector $\be\in \partial U(1)$ ($d\geq 2$).  Given $v>0$, we put  $\bv = v\be$ and  label with the superscript  $\,^{(\bv)}$ the objects defined by means of   $B_t-  {\bf v}t$  in place of $B_t$;  in particular, $B^{(\bv)}_t :=B_t- {\bf v}t$. 
Let $A$ be a  bounded Borel set and consider 
 $$H^{(\bv)}_{A}(\x, t;d\xi)  = P_\x\Big[ B^{(\bv)}_{\sigma^{(\bv)}(A)} \in d\xi,  \sigma^{(\bv)}_{A} \in dt\Big]/dt.$$
We put   $\ga(\cdot)= -\bv$ (constant function)    and  
$Z(t) = e^{\int_0^t \ga(B_u)d B_u -\frac12 \int_0^t \ga^2(B_u)du}.$
Then, for a bounded continuous function $\psi$  on  $[0,\infty)\times \partial A$,
\begin{eqnarray*}
\int_{[0,\infty)\times \partial A} \psi(u,\xi)H^{(\bv)}_{A}(\x, u;d\xi) du
=E_\x[Z(\sigma_A) \psi(\sigma_A, B_{\sigma(A)})]
\end{eqnarray*}
owing  to  the Cameron-Martin formula. Since with $P_\x$-probability one
$$Z(\sigma(A)) = \exp\Big\{-\bv\cdot B_{\sigma(A)}+\bv\cdot\x -\tst12 v^2 \sigma_A\Big\},$$
  we obtain
\beqn\label{5.2}
H^{(\bv)}_{A}(\x, t;d\xi)  =e^{\bv\cdot\x -\frac12 v^2t} e^{-\bv\cdot\xi} H_A(\x,t; d\xi). 
\eeqn

In the regime $x/t\to 0$, we have a simple asymptotic expression of $H_A(\x,t; d\xi)$ (see Section 3.1.3);   in particular, on  taking   a ball ($d\geq 3$)
 or a disc ($d=2)$  for $A$,   this leads to a quite explicit expression of  $P_\x[\sigma_{U(a)}^{(\bv)} \in dt]/dt$ as well as $H^{(\bv)}_{U(a)}(\x, t;d\xi)$ as is exhibited in  \cite[Section 7]{Ucal_b}.

\v2
 In below we restrict ourselves to the case when $\x$ is located not far from  $\bv t$.
On writing  $\x= \bv t+ \y$,  (\ref{5.2}) together with the first formula in  (vi) of Section 1 yields
 \v2
(vi$'$) $ \quad\quad  H^{(\bv)}_{A}(\x, t;d\xi)  \sim \, R_A^{2\nu}\La_\nu(R_Av)(2\pi t)^{-d/2} e^{\bv\cdot\y-\bv\cdot \xi } \la_A(\x/t;d\xi)\quad (d\xi\subset \partial A)$
\v2\n
{\it (as $t\to\infty$) with $\la_A(\bv;d\xi)$ defined in (\ref{lamda}). }
 \v2\n
 
Let $K$ be a compact set of $\R^d$. From the fact that  the law of  Brownian bridge does not depends on the strength of drift, we have
\beqn\label{Wsaus}
E_{{\bf 0}} \Big[{\rm vol}_d( S^{({\bv})}_K(t))\Big] = \int_{\R^d} E_{{\bf 0}} [{\rm vol}_d( S_K(t)) \,|\,B_t=\z ]  p_t^{(d)}(\x- {\bf v} \,t) d\z.
\eeqn
Noting that an overwhelming    contribution to the integral comes from a relatively small range ${\bf v}t + U(t^{\a})$  with any  $\a$  such that  $1/2<\a<1$  we infer from the second formula of (vi) that {\it as  $t\to\infty$, uniformly for  $v$ in a  bounded interval, }
 \v2
 (vi$''$) ${\displaystyle\quad\quad  E_{{\bf 0}} \Big[{\rm vol}_d( S^{({\bf v})}_K(t)) \Big] \,\sim\, \bigg(R_K^{2\nu} \La_\nu(R_Kv) \int_{\partial K} e^{-{\bf v} \cdot\xi}\la_K({\bf v}, d\xi)\bigg)t.}$\\
 \v2\n
By the continuity of $\la_K(\bv,d\xi)$ in $\bv$ at $\bv={\bf 0}$ and the expression of $\la_K({\bf 0},d\xi)$ given  in (\ref{lamda5}) the coefficient of  $t$  on the right-hand side above is asymptotic,  as $v \downarrow 0$,   to 
$$\pi/\lg (1/v) \quad\mbox{ or}\quad    \cp(K)  \quad \mbox{according as  \quad $d=2 \,\,$ or $\,\, d\geq 3$.}$$
{\it  For  the case of   a disc/ball  it follows from (\ref{EQ24}) that 
\[
E_{{\bf 0}} \, [ {\rm vol}_d(S^{(\bv)}_{U(a)}(t))] \,\sim\, \bigg(a^{2\nu}\La_\nu(av) \int_{\partial U(1)} e^{-a\bv\cdot\xi} g_{av}(\th_{\xi,\bv})m_1(d\xi)\bigg) t
\]
with $g_\a(\th)$ given  in (\ref{h_v}) of Section 2. Here $\th_{\xi,\bv}$ denotes the angle that  $\xi$  forms with $\bv$. }

\v2
 Consider the case when  $v\to \infty$ as well as  $t\to\infty$.  Let $K$ satisfy the condition (\ref{cond1}), i.e., ${\rm vol}_{d-1}\Big(\overline {{\mbox{dis-ct}}_{\be}(K)}\Big)=0$. Suppose also that   ${\rm vol}_{d-1}({\rm pr}_{{\bf e}}K)>0$ for simplicity.
Then,   substituting   $ t+\tau$  (with $|\tau|<< t$) for  $t$  in  (\ref{5.2}), and comparing the resulting formula with
 $$p_{t+\tau}(vt)=  p_t(vt)  \exp \Big\{\frac12 v^2\tau \Big(1-  \frac{\tau }{t} + \frac{\tau^2 }{t^2} -\cdots\Big) \Big\}, $$
we infer from (vii)  that  {\it  uniformly for $\tau$ in a finite interval, as $v\to\infty$ with $v/t\to 0$ (so that the ratio $\bv\cdot \xi/t$ as well as the third term in the exponent    tends to zero) } 
\v2
(vii$'$) $ \quad\quad\quad H_{K}^{(\bv)}(\x, t+\tau;d\xi)  \approx \, vp_t(v\tau) m_{K,\x/x}(d\xi) \quad\quad(  d\xi\subset \partial K)$
\v2\n
{\it  for  $\x=\bv t$,  which may be broken into the  two relations
\beqn\label{59}
q^{(\bv)}_{K}(\x,t+\tau)\sim\, vp_t(v\tau) {\rm vol}_{d-1}(\pr_{\x/x}K)
\eeqn
and}
\beqn\label{590}
 P_\x[ B^{(\bv)}_t\in d\xi\,|\, \sigma^{(\bv)}_K=t+\tau] \, \Longrightarrow \, \frac{m_{K,\x/x}(d\xi)}{{\rm vol}_{d-1}(\pr_{\x/x}K)}.
\eeqn
 In these formulae, as is readily ascertained, we may take $\x = \bv t + \y$ in place of $\x=\bv t$, provided  that  $|\y|^2/t \to 0$. 

 From (viii) and (\ref{Wsaus}) it also follows that  {\it as $v\to \infty$ and $t\to\infty$}
\v2
(viii$'$) $ \quad\quad\quad E_{{\bf 0}} \Big[{\rm vol}_d( S^{({\bv})}_K(t))\Big] \,\sim\, 
\Big[{\rm vol}_{d-1}({\rm pr}_{{\bf e}}K)\Big]vt\quad\quad  ({\bf v}= v\be).$
\v2
  Formula (\ref{590}) (at least with $\tau=0$) (as well as (viii$'$))  is intuitively comprehensible, but it
   (as well as (\ref{59}))  is not true (for  $\tau\neq 0$)  if  $v/ t$ is bounded away from zero when   the ratio
   $P_\x[\sigma^{(\bv)}_{\De_1}\in t+d\tau] /P_\x[\sigma^{(\bv)}_{\De_2}\in t+ d\tau] $ does not approaches unity, where $\De_1$ and $\De_2$ are any  (but distinct) two fixed  hyperplane perpendicular to $\be$.

\v2\n

\section{Miscellaneous estimates concerning $\sigma_A$}

The arguments presented  in this section  are made independently of those of preceding sections other than Section 2. Throughout this section $A$ denotes a bounded and non-polar Borel set of $\R^d$ and put $R_A=\sup \{|\y|: \y\in A^r\}$.  

\subsection{Uniform estimates for  $e_A(\x)/P_\x[\sigma_{\partial U(r)} < \sigma_A] $ }
Here we discuss  a  part of classical potential theory in two dimensions  related to  the  function $e_A(\x)$  and thereby prove Lemma \ref{lem3.01}. Most of what are presented below  are  known but some as given in Proposition \ref{prop6.1}   do not seem   to be found  in the existing  literature.

 Let  $d=2$.   Hunt \cite[Section 5.5]{Ht}  defines $e_A(\x)$ by
 \beqn\label{6.1}
 e_A(\x) = \pi g_{\Om_A}(\x,\y) +\lg|\x-\y| - E_\x[\, \lg| B_{\sigma(A)}-\y|]
 \eeqn
 ($H$ is written for $e_A/\pi$ in \cite{Ht}). 
Here  $g_{\Om_A}$ stands for  the Green function for the set  $\Om_A$:   it may be expressed  as 
$$g_{\Om_A}(\x,\y) = \int_0^\infty \frac{ P_\x[B_t\in d\y, \sigma _A >t]}{|d\y|}dt \quad\quad (\x,\y \in \Om_A),$$
where $|\cdot|$ designates the Lebesgue measure on $\R^2$; for  $\x\neq \y$, $g_{\Om_A}(\x,\y)$  is symmetric  and jointly continuous in $\x$ and $\y$ in the interior of $\Om_A$ and tends to zero as $\x$ (or $\y$) approaches $A^r$ (readers may refer to \cite{B}: Section 2.3 for  a precise 
definition, and Section 2.4 for  the above expression of $g_{\Om_A}$ and  the properties of it).
The function  
$e_A(\x)$ is defined for all  $\x\in \Om_A$ and 
   independent of  $\y\in \Om_A$ (this fact is seen from the arguments developed for (\ref{*}) below).  It follows that  \beqn\label{g_e}
   e_A(\x)= \pi \lim_{|\y|\to\infty} g_{\Om_A}(\x,\y);
   \eeqn 
  in particular $e_A$ is harmonic in the interior of $\Om_A$.
  Lemma  \ref{lem3.01}   follows from the following  proposition by using $(1+x)^{-1} > 1-x, x>0$.
   \v2

   \begin{Prop}\label{prop6.1} \, For $r >R_A$ and  $\x\in \Om_A\cap U(r)$,
    \beqn\label{sigma-s}
 P_\x[\sigma_{\partial U(r)} < \sigma_A] =\frac{e_A(\x)}{\lg(r/R_A)}\bigg(1+ \frac{m_{R_A}(e_A)}{\lg(r/R_A)}(1+\de)\bigg)^{-1}
\eeqn   
with $-2(R^{-1}_Ar +1)^{-1}\leq \de \leq 2(R^{-1}_Ar  -1)^{-1}$. Here $m_R(e_A) = \int_{\partial U(R)} e_A(\xi)m_R(d\xi)$.
 \end{Prop}

 In below we shall derive from  (\ref{6.1})  several formulae that relate   $P_\x[\sigma_{\partial U(r)}< \sigma_A]$ and  $e_A$;  Proposition \ref{prop6.1} will be among  them  (see Corollary \ref{cor6.1}).
    They are based on and refine  the  fact  
 that for each $R>R_A$, uniformly for  $\x\in U(R)\cap \Om_A$  
  \beqn\label{*}
(\lg r)P_\x[\sigma_{\partial U(r)} < \sigma_A]   \to  e_A(\x) \quad
 \mbox{ as} \quad  r\to \infty
    \eeqn  
(cf.  \cite[Theorem 11.2.14]{St}, where $A$ may be unbounded).  Put $\Om_r=  U(r) \setminus A^r$ and let $\tau(\Om_r)$ denote the first exit time from $\Om_r$. We give a proof  of (\ref{*})  
   resting on the  well-known identity
 \beqn\label{6.2}
 0=  \pi g_{\Om_r}(\x,\y) + \lg|\x-\y| - E_\x[\, \lg| B_{\tau(\Om_r)}-\y|\,] \quad\quad(\x,\y\in \Om_r, \x\neq \y).
 \eeqn
We  break the last expectations in (\ref{6.1}) and (\ref{6.2})  into two parts according as  $\sigma_{\partial U(r)}$ is larger or smaller than  $\sigma_A$.  Noting  that   
$ B_{\sigma_A}$ agrees with $ B_{\tau(\Om_r)}$ a.s.  on the event $\sigma_A < \sigma_{\partial U(r)} $ we observe
\begin{eqnarray} \label{6.3}
&&E_\x[ \lg| B_{\tau(\Om_r)}-\y|] -E_\x[ \lg| B_{\sigma(A)}-\y|] - (\lg r)P_\x[\sigma_{\partial U(r)} < \sigma_A]  \\
&&= E_\x[ \lg|B_{\tau(\Om_r)}-\y| -\lg r\, ;\, \sigma_{\partial U(r)} <\sigma_A]  - E_\x[ \lg| B_{\sigma(A)}-\y|\,; \sigma_{\partial U(r)} <\sigma_A], \nonumber
\end{eqnarray}
of which  each of the expectations on the right-hand side approaches zero as $r\to\infty$.
Using this equality we find that
\begin{eqnarray*}\label{616} 
 e_A(\x)-  (\lg r) P_\x[\sigma_{\partial U(r)} < \sigma_A] &=& \mbox{ RHS of (\ref{6.1})} - (\lg r) P_\x[\sigma_{\partial U(r)} < \sigma_A]  - \mbox{ RHS of (\ref{6.2})} \\
&=& \pi [g_{\Om_A}(\x,\y)-g_{\Om_r}(\x,\y)]  + \mbox{ RHS of (\ref{6.3})} \nonumber 
\end{eqnarray*}
and readily   conclude (\ref{*}) 
 and incidentally  that the right side of (6.1) is independent of $y$.

We bring in the function
$$\bar e_{A,R}(\x) = E_\x[e_A(B_{\sigma(U(R))})] \quad\quad (x\geq R\geq R_A).$$
  The following lemma is essentially a corollary of    (\ref{*}).
\begin{Lem}\label{lem6.1.0}\, Whenever $ x\geq R\geq R_A$,
\beqn\label{(a)}   
e_A(\x) = \lg (x/R)  + \bar e_{A,R}(\x). 
\eeqn
\end{Lem}
\v2\n
\pf\, For $\x$ with  $R<x<r$,
\begin{eqnarray}\label{e_A}
 P_{\x}[\sigma_{\partial U(r)} < \sigma_A]  
&=& P_{\x}[\sigma_{\partial U(r)} < \sigma_{\partial U(R)}]   \\
&&\, +\,
\int_{\partial U(R)}  P_{\x}[\sigma_{\partial U(R)} < \sigma_{\partial U(r)}, B_{\sigma(U(R))} \in d\xi] 
P_\xi[\sigma_{\partial U(r)} < \sigma_A], \nonumber
\end{eqnarray}
and, noting that the first term on the right-hand side  equals $\lg (x/R) /\lg (r/R)$, we multiply  $\lg r$, let  $r\to\infty$ and apply    (\ref{*})
to  obtain   (\ref{(a)}).  

\begin{Lem}\label{lem6.6} \,  Let $r > R\geq R_A$ and   $\x \in \Om_A\cap U(r)$.  Then
\begin{eqnarray} \label{e_log}
\qquad \quad  e_A(\x)= P_\x[\sigma_{\partial U(r)} < \sigma_A] \Big(\lg (r/R) + E_\x[\bar e_{A,R}(B_{\partial U(r)}) \,|\, \sigma_{\partial U(r)} < \sigma_A]\Big).
\end{eqnarray}
\end{Lem}

\v2\n
\pf \, 
Take $r^*>r$ and apply the strong Markov property to see  that 
\beqn\label{mu0}
P_\x[\sigma_{\partial U(r^*)} < \sigma_A]  =P_\x[\sigma_{\partial U(r)} < \sigma_A] \int_{\partial U(r)} P_\xi[\sigma_{\partial U(r^*)} < \sigma_A] \mu_{r,\x}(d\xi),
\eeqn
where 
$$
\mu_{r,\x}(d\xi) = P_\x[ B_{\sigma(\partial U(r))} \in d\xi \,|\,\sigma_{\partial U(r)} < \sigma_A].
$$
Then multiply the both sides by $\lg r^*$, let  $r^*\to \infty$  and  apply first   the formula (\ref{*}) with  $r^*$ in place of $r$, and then  (\ref{(a)})  to $e_A(\xi)$ that comes up  on  the right-hand side  under the integral sign, and one  then  finds the  identity of the lemma.
\qed
\v2

By using   an explicit form of the Poisson kernel for $\Om_{U(R)}=    \{\z\in \R^d:|\z|>R\}$, we have for $\y\in \partial U(r)$
$$\frac{ r - R}{r+R}   \leq \frac{P_\y[ B_{\sigma(U(R))}\in d\xi]}{m_R(d\xi)} \leq \frac{r+R}{r-R}  \qquad (r>R),$$
entailing $m_R(e_A)\frac{ r - R}{r+R}   \leq \bar e_{A,R}(\y) \leq \frac{r+R}{r-R}\, m_R(e_A)$,
which combined with Lemma \ref{lem6.6} yields
\begin{Cor}\label{cor6.1} \, For  $r> R_A$ and $\x \in \Om_A\cap U(r)$, 
\beqn\label{e3}
 m_{R_A}(e_A)\frac{r-R_A}{r+R_A} \leq \frac{e_A(\x)} {P_\x[\sigma_{\partial U(r)} < \sigma_A] } - \lg\Big(\frac{r}{R_A}\Big) \leq   m_{R_A}(e_A)\frac{r+R_A}{r-R_A}.
 \eeqn
\end{Cor}
\v2
Rearranging  the inequalities (\ref{e3})     leads to Proposition \ref{prop6.1} (write down it as the lower and the upper bounds of the ratio $e_A(\x)/P_\x[\sigma_{\partial U(r)} < \sigma_A]$ and take the reciprocal). It is noted that taking $\x$ from $\partial U(r)$ in it (or rather directly from  (\ref{(a)})) we have  for $x>R_A$
 \beqn\label{*0*}
m_{R_A}(e_A)\frac{x-R_A}{x+R_A} \leq e_A(\x)- \lg\Big(\frac{x}{R_A}\Big) \leq   m_{R_A}(e_A)\frac{x+R_A}{x-R_A}.
\eeqn
The second inequality of (\ref{e3}) entails  $e_A(\x)\leq \a +m_{R_A}(e_A)(1+\a^{-1})$ ($\a>0$, $x< (1+\a)R_A$), from which we deduce that 
\beqn\label{00*}
e_A(\x) \leq 2\sqrt{2m_{R_A}(e_A)}+m_{R_A}(e_A)  \quad \mbox{for}\quad \x\in \partial U(R_A).
\eeqn
Combining this with (\ref{*0*}) leads to the bound (\ref{e2}). 

\v2
The identity (\ref{(a)}) 
shows
\beqn\label{lim_thm}
\lim_{x\to \infty} (e_A(\x)-\lg x) = -\lg R  + m_R(e_A).
\eeqn
Since the left side is independent of $R$, the right-hand side must be a constant depending only on  $A$, which, known as  {\it Robin's constant} associated with $A$,  we denote by $V(A)$:
\beqn\label{RC}
V(A) = -\lg R + m_R(e_A).
\eeqn
  Since  $g_{\Om_A}$ is symmetric,  letting $x\to \infty$ in (\ref{6.1}) we find another (rather classical)  representation
\beqn\label{Rcst}
V(A) = e_A(\y) -\int_{\partial A} \lg |\xi -\y|  H_A^\infty(d\xi), \quad \y \in \Om_A.
\eeqn

Using    formula (\ref{Rcst})   instead of (\ref{(a)}) in the very  last step of the proof  of Lemma \ref{lem6.6}, we deduce from (\ref{mu0}) that
$$e_A(\x) = P_\x[\sigma_{\partial U(r)} < \sigma_A] \bigg[ V(A) + \int_{\partial U(r)}\mu_{r,\x}(d\xi') \int_{\partial A} \lg |\xi -\xi'|  H_A^\infty(d\xi)\bigg].$$
The repeated integral in the large square brackets may be written as $\lg r + \de(r)$ with $|\de(r)| \leq  |\lg( 1- R_A/r)|$.  
Thus 
$$|e_A(\x)/P_\x[\sigma_{\partial U(r)} < \sigma_A] - V(A)-\lg r |\leq  -\lg(1-R_A/r),$$
which  in terms of the {\it logarithmic capacity} defined by
\beqn\label{lcap}
{\rm lcap}(A) = e^{-V(A)}
\eeqn
 (normalized so that ${\rm lcap}(U(a))=a$) may be expressed as  in the following   
 \begin{Prop}\label{prop6.2} \, For  $r > R_A$ and for  $\x\in \Om_A\cap U(r)$  
  \beqn\label{**}
 \bigg| \frac{ e_A(\x)} {P_\x[\sigma_{\partial U(r)} < \sigma_A]}  -  \lg \Big(\frac{r}{ {\rm lcap}(A)}\Big) \bigg| \leq -\lg\Big(1-\frac{R_A}{r}\Big).    \eeqn
\end{Prop}
\v2

By virtue of (\ref{RC}) the twin  inequalities  of Corollary \ref{cor6.1} may be written as 
$$
  -\frac{2m_{R_A}(e_A)}{r/R_A+1} \leq \frac{e_A(\x)} {P_\x[\sigma_{\partial U(r)} < \sigma_A] } -\lg \Big(\frac{r}{ {\rm lcap}(A)} \Big)\leq  \frac{2 m_{R_A}(e_A)}{r/R_A-1}, 
$$
   which combined with Proposition \ref{prop6.2} yields  
\begin{Cor}\label{cor6.21}  \, If $2R_A< x\leq r$, then for some universal constant $C$, 
$$
 \bigg| \frac{ e_A(\x)} {P_\x[\sigma_{\partial U(r)} < \sigma_A]}  -  \lg \Big(\frac{r}{ {\rm lcap}(A)}\Big) \bigg| \leq C[1\wedge m_{R_A}(e_A)]\frac{R_A}{r}.  
 $$  
\end{Cor}

 \subsection{An  upper bound of $q_A$ $(d\geq 3)$}

 Let $d\geq 3 $ and   $\k_d, \k_d', \k_d''$ etc.  designate  constants that depends only on $d$, whose precise values are not important to the present purpose and may vary from  line to line.   In this subsection we prove
 
 \begin{Prop}\label{prop6.3}  \,  For $t\geq R_A^2$ and $\x\in \Om_A$,
 $$q_A(\x,t) \leq \k_{d}\cp(A)\La_\nu(2R_Ax/t)p_t^{(d)}(x).
 $$
 $(\cp(A)$ is the Newtonian capacity defined in Section 3.1.1.$)$
 \end{Prop}
 
 The factor  $\La_\nu(2R_Ax/t)$  in  Proposition \ref{prop6.3} may be replaced by  $\La_\nu(Rx/t)$ if $R>R_A$  (but with $\kappa_d$ depending on $R$) but  not  by $\La_\nu(R_Ax/t)$, for  if replaced,  the factor $\cp(A)$ is  possibly too small for the inequality to be valid. This   is caused by the concentration of the measure  $H_{U(R_A)}(\x,t; d\xi)$ at  $R_A\x/x$ as $x/t\to \infty$. 
 (Compare   with the result for $d=2$ given in Proposition \ref{prop6.4}.)
 
 The proof consists of several lemmas.

\begin{Lem}\label{lem6.1}  Put $\eta_\x=\eta_{A,\x} = {\rm dist}(\x,A^r).$  Then
for all $\x \in \Om_A$ and $t\geq 0$,
 $$P_\x[t<\sigma_A<\infty] \leq  \frac{\k_d}{(t\vee  \eta^2_\x)^{\nu}} \,\cp(A). $$
 \end{Lem}
 \v2\n\pf \,
 Put $\mu_{t,\x}(d\xi)= P_\x[B_{\sigma(A)}\in d\xi, t<\sigma_A <\infty]$. If $\fa$ is a  positive Borel function, 
 \beq
 &&\int \mu_{t,\x}(d\xi)\int \fa(\z)G^{(d)}(|\xi-\z|)|d\z| \\
 &&\quad =  E_\x\bigg[\int_0^\infty E_{B_{\sigma(A)}} [\fa(B_s)] ds; t<\sigma_A <\infty\bigg]\\                                            
  &&\quad =     E_\x\bigg[ \int_{\sigma(A)}^\infty \fa(B_s)ds; t<\sigma_A <\infty\bigg] \\
 &&\quad \leq  E_\x\bigg[\int_t^\infty \fa(B_s)ds\bigg] = \int_t^\infty ds\int p^{(d)}_s(\x -\z)\fa(\z)|d\z|.
 \eeq
Taking into account the fact that  the potential of $\mu_{t,\x}$ is lower semi-continuous and  maximized on the set $\overline{ A^r}$,   from the inequality above we infer  that
  \beq
 \int G^{(d)}(|\xi-\z|)\mu_{t,\x}(d\xi) &\leq &  \sup_{\z\in  \overline{ A^r}}\int_t^\infty p^{(d)}_s(\z-\x)ds\\
 &\leq&  \int_t^\infty p^{(d)}_s(\eta_\x)ds = \frac1{(2\pi)^{d/2}\eta_\x^{d-2}}\int_0^{\eta_\x/\sqrt t}u^{d-3}e^{-u^2/2}du.
 \eeq
 We integrate the left-most member w.r.t. $\z$ with  the equilibrium measure  of  $A$, denoted by $\mu_A$.
 The integration   results in  $\int P_\xi[\sigma_A<\infty]\mu_{t,\x}(d\xi)$, which in turn  equals the total charge of  $\mu_{t,\x}$ since $\mu_{t,\x}$ is concentrated on $A^r$.  On the other hand    the right-most member that is independent of $\z$  is evaluated to be at most  $\kappa_d (\sqrt t\vee \eta_\x)^{-(d-2)} =\k_d[t\vee  \eta^2_\x]^{-\nu}$, which  multiplied by $\mu_A(\bar A)=\cp(A)$  thus dominates $\mu_{t,\x}(A^r) = P_\x[t<\sigma_A<\infty]$,
  yielding   the  inequality of the lemma. \qed

 \begin{Lem}\label{lem6.2} \, Let $\eta_\x= {\rm dist}(\x,A^r)$ as above. Then,  {\rm (i)}  for all $\x\in \Om_A$ and $t\geq 1$,
$$q_A(\x,t)\leq \k_d (t\vee \eta_\x^2)^{-\nu}\,\cp(A); $$
and  {\rm (ii)}   for all  $t>0$ and  $\x\in \Om_A$ with  $\eta_\x \geq R_A$, $q_A(\x,t)\leq \k_d \cp(A)/R_A^d$.
 \end{Lem}
   \v2\n
   \pf \,  Let  $t\geq 1$. Since then $\inf_{\y\in U(1)} (t\vee\eta_{\x+\y}^2) \geq \frac14(t\vee \eta_\x^2)$,   by the preceding lemma  we have 
   $$\int_{t/3}^{2t/3} q_A(\x +\y,s)ds \leq \k_d' (t\vee \eta_\x^2)^{-\nu }\, \cp(A)\quad \mbox{ for}\quad  \y\in U(1),$$
    so that there exists $s^*\in [\frac13 t,\frac23 t]$ such that
 \beqn\label{6.01}
\int_{U(1)} q_A(\x +\y, s^*)|d\y|\leq  3\k_d' t^{-1} (t\vee \eta_\x^2)^{-\nu}\,\cp(A).
 \eeqn
 Here   we set  $q_A(\z, s)=0$ for $\z \notin \Om_A$, $s>0$.
 Denote by $\tau =\tau_\x$ the first exit time from the ball  $\x+ U(1)$. Then by strong Markov property 
 \begin{eqnarray}\label{6.009}
 q_A(\x,t) &=&\int_0^{t-s^*}  \int_{\partial U(1)} P_\x[\sigma_A > s, \tau \in ds, B_{\tau} -\x\in d\xi ]q_A(\x+\xi,t-s) \nonumber\\
 &&+  \int_{ U(1)} P_\x[B_{t-s^*} -\x \in d\y, \tau \wedge \sigma_A >t- s^*]q_A(\x+\y,s^*).
 \end{eqnarray}
By rotational symmetry we have
 \beqn\label{density_bdd}
 P_\x[\sigma_A > s, \tau \in ds, B_{\tau} -\x\in d\xi ] \leq P_\x[\tau\in ds, B_\tau -\x\in d\xi]\leq \k_d''  m_1(d\xi)ds,
 \eeqn
and using   this  as well as  Lemma \ref{lem6.1} we infer  that the first term (i.e., the  repeated integral) is dominated by   $4\k_d(t\vee \eta_\x^2)^{-\nu}\,\cp(A)$.
As for the second term we see   that  for  some universal constant $\la>0$, $P_\x[B_{t-s^*} -\x \in d\y, \tau \wedge \sigma_A >t- s^*]/|d\y|\leq C e^{-\la t}$, hence  the bound (\ref{6.01}) yields an estimate enough  for the one asserted in  (i). 

Let $R_A=1$ and  $\eta_\x\geq 1$. Putting  $f(s) = P_\x[\tau_\x \in ds]/ds (= P_0[\tau_{U(1)}\in ds]/ds$) with  the same $\tau_\x$ as above, we see that   for $t<1$,  $f(1-s)\geq Cf(t-s) $ for some $C>0$,  and 
$$q_A(\x,1) \geq \int_{0}^t f(1-s)ds\int q_A(\x+\xi,s) m_1(d\xi)\geq Cq_A(\x,t).$$
On taking the scaling relations given in  (\ref{scp}) into account  (ii) follows from (i). \qed
\v2

 The next lemma  provides   Harnack type estimates for  $q_A$. 
  \begin{Lem}\label{lem6.30} 
 If   $ 3R_A\leq x\leq |\y| \leq t/R_A$, then 
 \beqn\label{6.05} 
 q_A(\y, t) \leq \k_d \Big[q_A(\x, t) + \cp(A) R_A^{-d}e^{-y/2R_A}\Big]\quad\quad  (y=|\y|).
 \eeqn 
  \end{Lem}
 \v2\n
 \pf \, Put $R=\frac32 R_A$ and let   $ 2R\leq  |\y| \leq t/R_A$.  We use the representation
 \beqn\label{6.051}
  q_A(\y, t)  =\int_0^t ds \int_{\partial U(R)} H_{U(R)}(\y, s;d\xi) q_A(\xi,t-s).
  \eeqn
  In view of  Theorem \ref{thm2.3} there exist  positive constants $\k'_d$, $\kappa''_d$ such that  for $s> R_Ay$ $(y=|\y|)$,
 \beqn\label{6.06}
 \k'_d \,q(y,s;R) \leq H_{U(R)}(\y, s;d\xi)/m_R(d\xi) < \k''_d \,q(y,s;R).
 \eeqn
Now let $2R\leq x <y$. Then  $q(y,s;R) \leq Cq(x,s;R)$  ($s>0$).  Hence, comparing the integral (\ref{6.051}) restricted  on 
  the interval $[R_Ay,t]$ with the same integral but  with $\x$ replacing  $\y$ we have  
  $$ \int_{R_Ay}^t ds \int_{\partial U(R)} H_{U(R)}(\y, s;d\xi) q_A(\xi,t-s) \leq \k''_d\, q_A(\x,t).$$
   On the other hand 
\beq
\int_0^{R_Ay} ds \int_{\partial U(R)} H_{U(R)}(\y, s;d\xi) q_A(\xi,t-s) &\leq& \k_d\,\cp(A)\int_0^{R_Ay} q(y,s;R)ds \\
&\leq& \k'''_d \,\cp(A) p_{R_Ay}(y), 
\eeq
where Lemma \ref{lem6.2} (ii) and Theorem \ref{thm2.2} (or rather (\ref{6.9_d_Int})) are applied for the first and second  inequalities, respectively. Since $p_{R_Ay}(y)\leq R_A^{-d}e^{-y/2R_A}$, we obtain (\ref{6.05}) as desired. \qed

   \begin{Lem}\label{lem6.20} \,  
  If  $x /R_A>  (t/R_A^2)^{1/d} >1$,   then 
$q_A(\x,t)\leq  \k_d \cp(A)t^{-d/2}.$
  \v2
 \end{Lem}
  \v2\n\pf \,  We may suppose $R_A=1/3$ and $x>1$. 
     The proof parallels  that of  Lemma \ref{lem6.2} but  this time we consider the hitting of the ball $U(1)$  (instead of the exiting from $\x+U(1)$),  choose $s^*\in [\frac13 t,\frac23 t]$ so that
  \beqn\label{6.012}
  q_A(\xi_0, s^*)\leq   \k_d  t^{-1}(t\vee \eta_{\x_0}^2)^{-\nu}\cp(A)
 \eeqn
with any fixed $\xi_0\in \partial U(1)$ (which is possible owing to Lemma  \ref{lem6.1}) and  look at the decomposition
 \begin{eqnarray}\label{6.001}
 q_A(\x,t) &=&\int_0^{t-s^*}  ds\int_{\partial U(1)} H_{U(1)}(\x,s;  d\xi)q_A(\xi,t-s) \nonumber\\
 &&+  \int_{ \y\notin U(1)} P_\x\Big[B_{t-s^*} \in d\y, \sigma_{U(1)}   >t- s^*\Big]q_A(\y,s^*).
 \end{eqnarray} 
 On applying Lemma \ref{lem6.30} (with $\xi_0, s^*$ in place of $\x,t$)  and  (\ref{6.012})  in turn we deduce that
 $$
 q_A(\y,s^*) \leq  \k_d [q_A(\xi_0, s^*) + \cp(A)e^{-y/2R_A}] \leq \k_d' \cp(A)(t^{-d/2} + e^{-y}).
 $$
Since $E_\x[e^{-|B_s|}] \leq p_{s}(0)\int e^{-|\y|}|d\y| \leq Cs^{-d/2}$,  the second term on the right-hand side of (\ref{6.001})  is dominated by a constant multiple of $t^{-d/2}\cp(A)$.

With the help of $ p^{(d)}_s(x) \leq  p^{(d)}_{x^2/d}(x) = (d/2\pi e)^{d/2}x^{-d}$, we infer from Corollary \ref{upper_bd}  that 
  \beqn\label{ineqS}
  q(x,s;1) \leq \k_d x^{-d}\quad \mbox{if} \quad s\wedge x \geq1.
  \eeqn
 Together with  (\ref{6.06}) this shows that the inner integral of the repeated integral of the  first term is at most $\k_d x^{-d}\int q_A(\xi,t-s)m_1(d\xi)$ if $s> x$, whereas it is at most  $\k_d \cp(A)e^{-(x-1)^2/3s}$ for all $s>0, x>2$ in view of  Lemma \ref{lem6.2} (ii).  Hence,
  on employing  Lemmas  \ref{lem6.1} for the integral over $s\in [x, t-s^*]$ (if $x <t-s^*$) as well as the assumption of the  present lemma  the  first term  is dominated by $\cp(A)$ times 
  $\k_d' t^{-\nu} x^{-d} \leq \k_d' t^{-d/2}$
 as desired. \qed

  \begin{Lem}\label{lem6.3} 
If $t>R_A^2$, then for all  $\x\in \Om_A$,
 $$q_A(\x,t)\leq \k_d \cp(A) t^{-d/2}.$$
 \end{Lem}
 \v2\n
 \pf \,  Suppose  $R_A =1/2$ for simplicity.  Owing to the preceding lemma  it suffices to show  the inequality for  $x\leq t^{1/d}$. 
    For $r > 1$/2 we make decomposition
   \begin{eqnarray}\label{eq6.6}
   q_A(\x, t)  &=& \int_0^{t/2}\!\! \!\int_{\partial U(r)}  P_\x\Big[B_{s}\in d\xi,\sigma_{\partial U(r)}\in ds,  s<\sigma_A \Big]q_A(\xi,t- s) \nonumber\\
 && + \,\int_{U(r)\cap \Om_A}  P_\x\Big[B_{t/2}\in d\y, \tst 12 t< \sigma_{\partial U(r)}\wedge\sigma_A \Big]q_A(\y,\tst12 t).
\end{eqnarray}

Taking $r=1$ in it, we observe that the range of $x$ may be   further restricted to  $x> 2R_A=1$. Indeed,    for $x\leq 1$,  estimation of the first term on the right   is reduced to that in the case $x= 1$ and   the second term   is at most  $\k_d \cp(A) e^{-\la t}$  owing to Lemma \ref{lem6.2} (ii).

 Now  let $1\leq x \leq t^{1/d}$
 and put  $r =2 t^{1/d}$ in (\ref{eq6.6}).  Then in the right-hand side of the decomposition 
 the second term is at most   $ \k_d e^{-\la t/ r^2}\cp(A)$ for some  constant  $\la>0$ owing to  the bound  $q_A(\cdot\,, t) \leq  \k_d \cp(A)$ ($t\geq 1$), whereas Lemma \ref{lem6.20} shows that the repeated integral is dominated by $\k_d t^{-d/2}\cp(A)$. Thus Lemma \ref{lem6.3} has been proved. \qed

 \v2\n
 
 Combined  with   Corollary \ref{upper_bd}   as well as with the last lemma
the following one, virtually a corollary of Lemma \ref{lem6.20},     concludes the proof of proposition \ref{prop6.3}.
\begin{Lem}\label{lem6.71} 
There exists a constant  $\k_{d}$ such that  
  if   \,$x\geq  \sqrt t \geq 2R_A$, then 
  $$q_A(\x,t)\leq \k_d \frac{\cp(A)}{R_A^{2\nu}} q(x,t; 2R_A).$$
 \end{Lem}
 \v2\n
 \pf \, Let $R_A=1/2$. Recalling (\ref{3.6a}), we apply (ii) of  Lemma \ref{lem6.2} and  Corollary \ref{upper_bd2} in turn we  deduce that  for $x>t$,
\beqn\label{upper_bd3}
 q_A(\x,t)\leq \k_d'\cp(A)\int_0^t q(x,s;1)ds\leq \k_d\cp(A)\Lambda_\nu(2x/t)p_t^{(d)}(x),
\eeqn
hence the upper bound of the lemma in view of Theorem \ref{thm01} and the scaling relation.

For the case $x\leq t$ we split the outer integral  in the right-hand side of (\ref{3.6a}), and write  $I_{[0,t/2]}$ and $I_{[t/2,t]}$ for the corresponding parts  (as in Section 3.1). Then applying Lemma \ref{lem6.20} and Corollary \ref{upper_bd2}  
we have
$$I_{[t/2,t]} \leq  \k_d\cp(A)t^{-d/2} \int_0^{t/2} q(x,s;1)ds \leq \k_d'\cp(A)q(x,t;1).$$
 On the other hand,
   using  Theorem \ref{thm2.3}   and   the inequality  $q(x,t-s;1)\leq \k_d' q(x,t;1)$ ($0<s<t/2$)  we obtain  that for $1\leq x \leq t$,
$$I_{[0,t/2]} \leq  C q(x,t;1)\int_0^{t/2} ds\int_{\partial U(1)} q_A(\xi,s)m_1(d\xi)    \leq \k_d'q(x,t;1) \int_{\partial U(1)} P_{\xi}[\sigma_A<\infty] m_1(d\xi),$$
Owing to   (\ref{cap1}) the last integral is equals $\cp(A)/\cp(U(1))$. Thus we obtain the required bound in the case $x\leq t$.  \qed

 \subsection{Some upper bounds of $q_A$ $(d = 2)$}

The statements  corresponding to  Proposition \ref{prop6.3}  for $d=2$ are given   by the following  one.

   \begin{Prop}\label{prop6.4}  There exists a universal constant $C$ such that for   $t>R_A^2$,
    \[
  q_A(\x,t)\leq \left\{\begin{array} {ll}  {\displaystyle C\frac{e_A(\x)}{t[\lg (t/R_A ^2)]^2}}  \quad &\mbox{if} \quad   2R_A\leq x< \sqrt t,\\[4mm]
C \La_0 (R_Ax/t)p^{(2)}_t(x) \quad &\mbox{if} \quad  x \geq \sqrt t,
 \end{array}\right.
     \]
and 
$$q_A(\x,t) \leq \frac{C \b_A}{t[\lg (t/R_A ^2)]^2} \quad \quad \mbox{if} \quad  \x\in \Om_A\cap U(2R_A),
$$
where $\b_A = m_{2R_A}(e_A)$.
\end{Prop}

  One may follow  the proof  for  $d\geq 3$;      in place of Lemma \ref{lem6.1} we can  derive  the following bound 
   \[
  P_\x[\sigma_A>t] \leq  Ce_A(\x)/\lg (t/R_A^2) \quad  (x>2R_A,   t> R_A^2)
  \]
  by an  argument analogous to that of the second half of the proof of Lemma \ref{lem6.3} (applied to $P_\x[\sigma_A>t]$ in place of $q_A(\x,t)$)  with the help of Proposition \ref{prop6.1}.  However we adopt a somewhat different method in which  the bound above though sharp for itself  (if $\lg x<\!\!<\lg t$) is not so useful. 
  The arguments  set forth  in below   rests  on   the trivial bound
   $\int_{t_1}^{t_2}  q_A(\x,s)ds \leq 1$;  on using it  in place of Lemma \ref{lem6.1}  the same proof of Lemma \ref{lem6.2} leads to 
  \begin{Lem}\label{lem0} \,  If $  t\vee {\rm dist}(\x, A^r) \geq 1, \x\in \Om_A$, then 
   $\, q_A(\x, t) \leq C_0\,$ for some  universal constant $C_0$.
  \end{Lem}
  \v2\n 
   
 \begin{Prop} \label{lem6.9}  There exists a universal constant $C$ such that  for $t > R^2_A$ and 
     $\x\in \Om_A$,
 \[
  q_A(\x,t)\leq \left\{\begin{array} {ll}  {C}/{t\lg (t/R^2_A)} \quad &\mbox{if} \quad x<\sqrt t,\\ [2mm]
C q(x,t; 2R_A) \quad\,\, &\mbox{if} \quad  x\geq \sqrt t.  
 \end{array}\right.
     \]
  \end{Prop}
\v2\n
  \pf\,  
   Suppose $R_A=1/2$  for simplicity and    let  $R=2R_A=1$.
  We split the outer integral in (\ref{3.6a}) at $t/2$ and make the decomposition 
 \beqn \label{B9}
q_A(\x,t) = I +II\quadd (x > 1),
\eeqn
where
$$
I = \int_0^{t/2} \!ds \int_{ U(1)} H_{ U(1)} (\x, t-s; d\xi) q_A(\xi, s)
$$
and
$$
II = \, \int_0^{t/2} ds \int_{ U(1)}H_{U(1)} (\x, s; d\xi) q_A(\xi, t-s).$$
 We shall  apply the following bounds  
 \begin{eqnarray}\label{6.9}
 P_\x[\sigma_{U(1)}< t]& \asymp&  \frac{1}{1\vee  \lg t}\Big(1\vee\lg \frac{t}{x^2}\Big) \quad \mbox{if}\quad t>x^2>1,  \\
 &\asymp& \frac{t^2}{x^2} \La_0\bigg(\frac{x}{t}\bigg) p_t^{(2)}(x)\quad \mbox{if}\quad x^2\geq t>1/4,
 \label{6.10}
 \end{eqnarray}
 which are  deduced from Theorems \ref{thm2.2} and \ref{thm2.20} 
 (cf. Corollary 14 of \cite{Ubes}).

First of all we point out two facts that are easy to verify. Firstly, using   (\ref{6.10})  as well as Lemma \ref{lem0}  we deduce as in (\ref{upper_bd3})    that
\beqn\label{large_x}
q_A(x,t)\leq Cq(x,t;1)\quad \mbox{for}\quad x>t> 1/4,
\eeqn
namely  the bound of the lemma holds true for $x\geq t$. 
    Secondly,   use  Theorem \ref{thm2.3}   and   the inequality  $q(x,t-s;1)\leq C' q(x,t;1)$ ($0<s<t/2$) to obtain  that for $1\leq x\leq t$, $I \leq C_1\int_0^{t/2}ds\int_{\partial U(1)}m_1(d\xi)p_A(\xi,s)  q(x,t;1) \leq C_1q(x,t;1)$,  which combined with (\ref{large_x}) yields
\beqn \label{I_1}
I \leq C q(x,t;1) \quad \mbox{for}\quad x>1, t>1/4.
\eeqn

For estimation of $II$ we start with   the  crude  bound  
\beqn\label{Bbd}
q_A(\x,t) \leq C/(x^2\vee 1)  \quadd (t>1/4, \x \in \Om_A).
\eeqn
For verification  we apply  the bound  $q(x,s;1) \leq C'/x^2$ valid for all $s>0,\,  x>1$, which  is   derived  from   Theorem \ref{thm01} and Lemma \ref{thm5}  in view  of $x^2p_s^{(2)}(x)\leq 1$. In the decomposition (\ref{B9}) we take $x$ in place of $t/2$ as the splitting point of the integral in (\ref{3.6a}) and, by   using Theorem \ref{thm2.3}, observe that for   $x<t$ 
$$q_A(\x,t) \leq C' \int_x^t  ds\int_{ U(1)}  m_1(d\xi)  q(x, s;1)q_A(\xi, t- s)+ C_0\int_0^x q(x,s;1)ds,$$
of which the first term on the right-hand side is dominated by $C''/(x^2\vee1)$ and the second one by $C''e^{-x/2}$ in view of (\ref{6.10}). Combined with  (\ref{large_x}) this shows (\ref{Bbd}). 
 
 On using the expression
\begin{eqnarray} \label{0C}
q_A(\x,t)  
= \int_{\Om_A} P_\x[  \sigma_{A}>  \tst 12 t, B_{t/2}\in d\y] q_A(\y,\tst 12 t) 
\end{eqnarray}
the bound (\ref{Bbd})  entails
 $$q_A(\x,t)  \leq \int_{\R^2}p^{(2)}_{t/2}(\y-\x)  \frac{C}{|\y|^2\vee1} |d\y| \quadd (t>1/4, x > 1).$$
 The integral on the right-hand side  attains the maximum at $\x=0$ and an easy computation shows that it is $O(t^{-1}\lg t)$. By continuity the bound is valid for $x =1$. Thus
\begin{eqnarray} \label{E}
q_A(\xi,t) \leq  C'\frac{1\vee \lg t}{t} \quadd \xi \in \partial U(1).
\end{eqnarray}
 Substitution into  the  expression that defines $II$   yields 
 \beqn\label{6.90}
 II\leq C' P_\x[\sigma_{U(1)}<\tst 12 t]\frac{1\vee \lg t}{t}.
 \eeqn
Using  (\ref{6.9})  and (\ref{6.10}),   we infer that
\beqn\label{oi}
II\leq \frac{C}{t} \Big(1\vee\lg \frac{t}{x^2}\Big)e^{-(x-1)^2/t}   \quadd (x > 1, t> 1/4).
\eeqn
 This together with  the bound  (\ref{I_1})  of  $I$ gives an upper bound of  $q_A(\x,t) $. 
Noting  that the upper bound  of  $II$ above is also an upper bound of $q(x,t;1)$ on $x< \sqrt t$, we   substitute  for   $q_A(\y,t/2)$ in   (\ref{0C})  the  bound of  it  just obtained and   observe that  for $\xi\in \partial U(1)$,
\begin{eqnarray} \label{6.92}
q_A(\xi,t)& \leq& C_1\int_{|\y|<\sqrt t} p^{(2)}_{t/2}(|\y-\xi|)\frac{1+\lg (t/|\y|^2)}{t}|d\y|  + \frac{C_1}{t}\int_{|\y|\geq\sqrt t} p^{(2)}_{t/2}(|\y-\xi|) |d\y| \nonumber \\
&\leq& \frac{C_2}{t}.
\end{eqnarray}

We have  derived (\ref{oi}) from (\ref{E})  and then (\ref{6.92}) from (\ref{oi}). Repeating the same procedure once more but with the bound (\ref{6.92}) in place of (\ref{E}), we plainly gain the  factor of $1/\lg t$ for the right-hand sides of (\ref{oi}) and (\ref{6.92}). The resulting bound of $II$  being also an upper bound of    $q(x,t;1)$ for $x< \sqrt t$,  we can repeat it further once, which results in  
 \beqn\label{6.91}
II \leq \frac{C}{t(1\vee \lg t)^2} \Big(1\vee \lg \frac{t}{x^2}\Big)e^{-(x-1)^2/t}  \quadd (x \geq 1, t> 1/4).
 \eeqn 
Combined with (\ref{I_1}) again this shows  the bound of the lemma for $x>1$.  The case   $x \leq  1$ can be  easily reduced to the case $x = 1$ owing to the bound $q_A(\x,t)< C_0$ ($\x\in \Om_A, t>1$) (see the second paragraph of the proof of Lemma  \ref{lem6.3}).
  The proof  of Proposition \ref{lem6.9} is complete.    \qed

 \v2\n
{\it Proof of Proposition \ref{prop6.4}. } 
Since  $e_A(\x) \geq \lg (x/R_A)$ for $x> R_A$ owing to   (\ref{(a)}), the first bound of the proposition  follows from Proposition \ref{lem6.9} in  the case $x>t^{1/4}$  (note that $q(x,t;1)\asymp 1/t\lg t$ uniformly for  $t^\de< x<\sqrt t$ with any $\de>0$). On the other hand,  
on employing   Proposition \ref{prop6.1} as well as  Proposition \ref{lem6.9}, we make the same argument as in the second  half  of the proof of Lemma \ref{lem6.3} but with $r= 2 t^{1/4}$ instead of $r=2t^{1/d}$  to obtain  the asserted bound   for  $2R_A\leq x\leq t^{1/4}$.  

The last bound is deduced from what we have just proved. For, if we let $R_A=1$, then using  the decomposition (\ref{-51}) with $T=\sqrt t$ and $r=2$ (and with $d=2$) leads to
$$q_A(\x,t) \leq \sup_{\xi\in \partial U(2), s<\sqrt t} q_A(\xi, t-s) + Ce^{-\la\sqrt t},$$
by which the stated deduction is immediate.

It remains  to show that
$$q_A(\x,t) \leq \k_d q(x,t; R_A) \quad \mbox{for} \quad x \geq t/R_A.$$
(Note that here we are concerned with  $\sigma_{U(R_A)}$ instead of $\sigma_{U(2R_A)}$.)
For the proof we make use of  Theorem  2.4,   which entails that there exists a family of Borel  measures on $\partial U(a)$, say $(\mu_{v}(d\xi),{\bf v}\in \R^2)$,   such that $\mu_{{\bf v}}(\partial U(a)) <  C $   and 
\beqn \label{Eq_pp}
 H_{U(a)}(\x,s;d\xi) \leq  q(x,t; a)\mu_{a\x/t}(d\xi) \quad \mbox{if}\quad t/2 <s <t.
 \eeqn
Let $R_A=1/2$ as in the proof of Lemma \ref{lem6.10} and  define $I$ and $II$ as in  (\ref{B9}) but with $U(\frac12)$ in place of $U(1)$. Since $\int_0^{t/2}q_A(\xi,s)ds \leq 1$ for all  $\xi\in \partial U(\frac12)$ (where  $q_A(\xi,\cdot)$ is regarded as the Dirac delta  function if $\xi\in A^r$),  (\ref{Eq_pp})  immediately gives the required bound  for $I$ in (\ref{B9}), whereas  (\ref{6.90}) (together with  (\ref{6.9}) with $s=t/2$) provides a  bound of $II$ sufficient for the present purpose. 
\qed
\v2

We remark that the method used in the proof  of Proposition \ref{lem6.9} can be applied for the case $d\geq 3$ although equally involved as a whole.

It is also  remarked that  the procedure  mentioned  in the last step of the proof of Proposition \ref{lem6.9} may be applied  once more.    This   no longer  yields a  better bound of $q_A(\xi, t), \xi\in \partial U(1)$,  for  on the interval $\sqrt t/e <x <M\sqrt t$ with any $M>1$,  $q(x,t;1)$ is larger than a positive multiple of the bound of  $II$ in (\ref{6.91}) and     gives   rise to a term which is  the same order as $ 1/t\lg t$. If we apply Proposition \ref{prop6.4} however we can  obtain a better bound of $II$ but with the extra factor $\b_A$, which result we state as a lemma.

 \begin{Lem}\label{lem6.10} ~ With some universal constant $C$ and  $\b_A = m_{2R_A}(e_A)$
 $$ \int_0^{t/2} \int_{ \partial  U(R)}H_{U(R)} (\x, s; d\xi) q_A(\xi, t-s)ds\leq \frac{C\b_A}{t(\lg t)^3} \Big(1\vee\lg \frac{t}{x^2}\Big)e^{-x^2/t} \quad\quad (x> R_A, t>R_A^2).
$$
\end{Lem}

 \subsection{A lower bound of  $P_\x[\sigma_A<t]$ in case   $x/t>1$}
 
 Let
 $$\la(A)  = \left\{ \begin{array} {ll} R_A^{-1}\cp(A\times [-R_A,R_A])  &\quad \mbox{if} \quad d=2, \\
R_A^{-d+2}  \cp(A)  &\quad \mbox{if} \quad d\geq 3.
  \end{array}\right.
  $$
  ($\cp$ in  case $d=2$ stands for the three-dimensional capacity.)
  Note that $\la(A/R)= \la(A)$.  The following result is used in Step 4 of the proof of Lemma \ref{lem0a}.
  \v2
  \begin{Prop}\label{propA6} \, For  $\e\in (0,1/6d]$, $t>0$ and  $x\geq (t/R_A)\vee R_A$,
$$P_\x[\sigma_A < t]  \geq [\e \k_d R_A^{d} ] \,\la(A)  p^{(d)}_t(x)\exp\{-(R_Ax/t)[ 1+ (R_Ax/t)^{-\frac12} + R_Ax^{-1} +  6\e]\}.$$
Here   $\k_d$ designates  a constant that depends only on $d$ as in Section 6.2.
\end{Prop}
\v2
  Before proceeding to the proof of Proposition \ref{propA6}  
  we present an upper bound.
\v2
\begin{Lem}\label{lemA60}
For any $\e>0$ there exists a constant $C$ such that  for $x> (1+\e)R_A$ and for $t>0$ if $d\geq 3$ and  $0<t< R_A^2$ if $d=2$,
 $$P_\x[\sigma_A < t] \leq C  \la(A) P_\x[\sigma_{U((1+\e)R_A)} <t].$$
  \end{Lem}
  \v2\n
  \pf\,
Let $R_A=1$. If $d\geq 3$, then on writing  the probability $P_\x[\sigma_A < t]$ as the double integral 
$\int_0^t ds \int_{\partial A}P_\xi[\sigma_A <t-s] H_{ U(1+\e)}(\x,s;d\xi)$ 
the inequality follows from  the Harnack inequality that shows $P_\xi[\sigma_A <\infty] \leq C_{\e,d} \cp(A)$, $\xi\in \partial U(1+\e)$.  As for the case  $d=2$ we have
$$P_\x[\sigma_A < t] = P^{BM(3)}_{(\x,0)}[ \sigma_{A\times \R}<t] < C'P^{BM(3)}_{(\x,0)}[ \sigma_{A\times [-1,1]}<t]$$
for $t < 1$ with $C'=1/ P^{BM(1)}_0[\sigma_{[\R\setminus [-1,1]} \geq 1]$, where $P_{\y}^{BM(d)}$ denotes the law of $d$-dimensional Brownian motion started at $\y$,  and the result follows from that for $d=3$. \qed

\v2
 For the proof of Proposition \ref{propA6} we show two preliminary lemmas. 
 \v2
\begin{Lem} \label{lemA61}\, Let $d\geq 3$. Then
\quad ${\displaystyle   \int_{U(R_A)} P_\y[\sigma_A < R_A^2]|d\y| \geq \k_d  R_A^2 \cp(A).}$
\end{Lem}
\v2\n
\pf\,  Let  $R_A=1$ and put  $\fa(\y) = P_\y[\sigma_A\leq 1/2]$. 
Then on the one hand
\beqn\label{A6-1}
\int_{U(1)} P_\y[\sigma_A < 1]|d\y| \geq \int_{U(1)} |d\y| \int_0^{1/2}ds\int_{\partial U(1)} H_{\partial U(1)}(\y,s; d\xi) \fa(\xi) = C m_1(\fa),
\eeqn
where $m_r(\fa) = \int \fa dm_r$ and $C= \int_{U(1)} P_\y[\sigma_{\partial U(1)} <\frac12]|d\y| $.
On the other hand,  making  decomposition  
$\cp(A)/\cp(U(1)) = P_{m_1}[\sigma_A<\infty] =m_1(\fa) + \sum_{k=1}^\infty P_{m_1}[\frac12 k <\sigma_A\leq \frac12(k+1)],$
we deduce
\beq
{\textstyle P_{m_1}[\frac12 k <\sigma_A \leq \frac12(k+1)]} 
&\leq&  E_{m_1}[\fa(B_{k/2})]\\ 
&\leq& \frac1{(\pi k)^{d/2}}\bigg(\int_{U(1)}\fa(\y)|d\y| + \om_{d-1}\int_1^\infty m_r(\fa) r^{d-1}dr \bigg)
\eeq
for $k\geq 1$. Noting $m_r(\fa) \leq [\int_0^{1/2}q(r,s;1)ds ] m_1(\fa)$, we apply Lemma \ref{thm5} to evaluate the second integral in the big parentheses to be bounded above by a constant multiple of $m_1(\fa)$. Now employing
 (\ref{A6-1}) we conclude that 
 $\cp(A) \leq C \int_{U(1)} P_\y[\sigma_A < 1]|d\y|$ as desired. \qed
\v2
\begin{Lem} \label{lemA62}\, Let $d\geq 3$.  Then for $t\leq R_A^2$,
 $${\displaystyle \int_{U((1+\sqrt{dt})R_A)} P_\y[\sigma_A < t]|d\y| \geq \k_d  \cp(A)  t.}$$
\end{Lem}
\v2\n\pf\,  Let  $R_A=1$. Denote the $d$-dimensional cube of side  $r$ and center $\x$ by $Q(\x,r)$
 and 
 for lattice points ${\bf k}\in \Z^d$ put  $A_{\bf k} =  Q({\bf k}, 1)\cap (A^r/\sqrt t)$.
 Note that  $Q({\bf k},1)  \subset U((1+ \sqrt{dt})/\sqrt t)$ if $A_{\bf k}  \neq \emptyset$, 
Then by scaling property of Brownian motion
\beq
\int_{U(1+\sqrt{dt})} P_\y[\sigma_A < t]|d\y| &=& \int_{U(1+\sqrt{dt})} P_{\y/\sqrt t}[\sigma_{A/\sqrt t} < 1]|d\y| \\
& =& t^{d/2}\sum_{{\bf k}\in \Z^d} \int_{U([1+ \sqrt{dt}]/\sqrt t) \cap Q({\bf k},1) } 
P_\z[\sigma_{A/\sqrt t} <1] |d\z|\\
&\geq& t^{d/2}\sum_{{\bf k}\in \Z^d} \int_{ Q({\bf k},1)} P_\z[\sigma_{A_{\bf k}} <1] |d\z|\\
&\geq & t^{d/2}\sum_{{\bf k}\in \Z^d} \k_d \cp(A_{\bf k})\\
&\geq &\k_d   t^{d/2} \cp(A/\sqrt t)\\
&=& \k_d \cp(A)t,
\eeq
where Lemma \ref{lemA61} is used  for  the second inequality and the sub-additivity    of the capacity for the third.   \qed
\v2\n
{\it Proof of Proposition \ref{propA6}.} \, Suppose $d\geq 3$ and write  $v= x/t$.  Obviously 
$$P_\x[\sigma_A<t] \geq \int_{\R^d} p_{t-\e R_A/v}(\y-\x) P_\y[\sigma_A < \e R_A/v]|d\y|.$$
Let $R_A=1$.  Restricting the range of integration of the integral above  to $U(1+\sqrt{d\e/v})$ we are going  to apply 
Lemma \ref{lemA62}. Using $1/(1-r)<1+2r$ ($0< r<\frac12$) and putting  $\a:=\sqrt{d\e/v}$, we see that
if $\y\in U(1+\sqrt{d\e/v}) $ and $\e/x<\frac12$, 
$$\frac{|\y-\x|^2}{2(t-\e/v)} \leq \frac{(x+ 1+\a)^2}{2t}\Big(1+ \frac{2\e}{tv}\Big) =\frac{x^2}{2t} +v\bigg[ 1+\a + \frac{(1+\a)^2}{2x} +\e\Big(1+\frac{1+\a}{x}\Big)^2 \bigg]. $$
Thus if $\sqrt{d\e}<\sqrt 2-1$ (valid if $d\e< 1/6$) and $v> 1$ (so that $1+\a <\sqrt 2$) and if $x>1$, this entails
$$p_{t-\e /v}(\y-\x)  >  p_t(x) e^{-v(1+ (\sqrt 2-1)  v^{-1/2}+ x^{-1} + 6\e)},$$
and by Lemma \ref{lemA62} $P_\x[\sigma_A<t] \geq \e \k'_dp_t(x) e^{-v(1+ v^{-1/2}+ x^{-1} + 6\e)}\cp(A)$, as desired. The case $d=2$ is reduced to the case $d\geq3$ as in the proof of  Lemma \ref{lemA60}. \qed

\section{Appendix}

Let $A$ be a bounded, non-polar,  Borel  set as  before.
\v2\v2\n
{\sc A.1.  Harmonic measure of heat operator}

\v2
Here we give a  brief exposition  of  the well known fact that 
$H_A(\x,t; d\xi)dt$ is the lateral component of the caloric measure for the exterior of the cylinder with base $A$.  
 Given $t>0$, the space-time Brownian motion $Y_s=(B_s,t-s)$ ($0\leq s\leq t$) under the law $P_\x$
is regulated by the heat operator $\frac12 \De - (\partial/\partial s)$. Let  $A^r$ designate the set of  regular  points of $A$ and  $\Om_A$   the unbounded component of $\R^d\setminus A^r$ as in Section 1, put
$$D = \Om_A \times (0,\infty) =\{(\z,t): \z\in \Om_A, t>0\}$$
and consider the following  \lq Dirichlet problem'  for the heat operator:
\begin{eqnarray}
&&\Big(\frac12 \De - \frac{\partial}{\partial t}\Big) u =0 \quad \mbox{on} \quad D \nonumber\\
&& u = \fa \quad \mbox{on} \quad \partial_{\rm reg} D := A^r\times [0,\infty) \cup \Om_A\times \{t=0\},
\label{bry_cond}
\end{eqnarray}
where, $\fa$ is any bounded continuous function on $\partial_{\rm reg} D$  and  the boundary condition  (\ref{bry_cond}) is   interpreted in a reasonable way.
The caloric measure,  $\mu_D(\x,t, d\xi ds)$ say, for  $D$ at a reference point $(\x,t)\in D$ is defined  as such a measure kernel  that    the above boundary value problem may be  solved in the form 
\[
u(\x,t)= \int_{\partial_{\rm reg} D} \fa(\xi, s) \mu_D(\x,t; d\xi ds),
\]
whereas the solution to
the same   problem is  represented by the expectation:
$$ u(\x,t)=  E_\x[\fa(B_{\sigma(A)}, t- \sigma_A); \sigma_A <t] +E_\x[\fa(B_t,0); \sigma_A  >t].$$
The first expectation on the right-hand side  above is expressed as an integral by the measure kernel $Q_A(\x,dtd\xi)$   
given in (\ref{Eq1}).  
The function $u^A(\x,t) := P_\x[\sigma_A \leq t] =  Q(\x, \partial A \times [0,t))$ satisfies the heat equation in the interior  $D^\circ$, hence its partial derivative  $q_A(\x,t) = u_t^A(\x,t)$   is not only well-defined but smooth in $D^\circ$ as is $u^A$; similarly for the derivative $H_A(\x,t; d\xi) = Q_A(\x,dtd\xi)/dt$. This assures that 
  the caloric measure  $\mu_D$ restricted to the lateral part of $\partial_{\rm reg}D$  is written as
$$\mu_D(\x,t; d\xi dt))\Big|_{A^r \times [0,\infty)}= H(\x,t-s;d\xi)ds,
$$
thus providing  the probabilistic expression of $\mu_D$ (cf. Hunt \cite{Ht}),   the one on the initial  boundary   $\Om_A\times \{0\}$ being  of course given by $P_\x[ B_t\in d\xi, \sigma_A>t]$ ($\xi\in \Om_A$).  If $\R^d\setminus A^r$  
 is assumed to be a Lipschitz domain, we know that the measure $Q(\x, \cdot)$ and  surface measure on $\partial A \times (0,\infty)$ are  mutually absolutely continuous \cite{FS}.

\v2\v2\n
{\sc A.2.  Asymptotics of the distribution of  $\sigma_A$}
\v2
 
 In below we give some asymptotic   estimates of  the distribution function  $P_\x[\sigma_A \leq t]$ for large time only in the case  $x/t\to 0$, when   those  of the density $q_A(\x,t)$ are explicit enough. For the other case Propositions \ref{prop6.3} and \ref{prop6.4} would be enough  if it is upper bound what one might need
 (see also Section 6.4).   
 \v2\v2\n
{\bf Proposition A.1.} {\it   Let $d\geq 3$. Then, as $t\to \infty$ and  $x/t\to 0$   
 \beqn\label{6.5.1}
 P_{\x}[t<\sigma_A <\infty] = \cp(A)P_\x[\sigma_A =\infty] \int_t^\infty p^{(d)}_s(x)ds(1+o(1)) 
 \eeqn
 uniformly for  $\x\in \Om_A$;  and as $x \to \infty$ and  $x/t\to 0$ }
  \beqn \label{6.5.2}
   P_{\x}[\sigma_A \leq t] = \cp(A)\int_0^t  p^{(d)}_s(x)ds(1+o(1)). 
\eeqn

 \v2
The proposition above  follows from Theorem \ref{thm3.1}:  (\ref{6.5.1}) is immediate, whereas  for the proof of  (\ref{6.5.2}),   if $\lim  x^2/t  = \infty$, one may apply  Proposition \ref{prop6.3}   to see that  $P_\x[\sigma_A<t/2]$ is negligible; if otherwise, use $P_{\x}[\sigma_A <\infty] = {\rm Cap}(A)G^{(d)}(x)(1+o(1))$ together with (\ref{6.5.1}).   It is noted that 
 $$\int_0^t p^{(d)}_s(x) ds =  \frac{G^{(d)}(x)}{\Gamma(\nu)}\int_{x^2/2t}^\infty  e^{-y}y^{\nu-1}dy$$
  and similarly for   $\int_t^\infty p^{(d)}_s(x) ds$.  
\v2\v2\n   
{\bf Proposition A.2.} {\it  Let $d=2$ and $\e>0$. Then, uniformly for $\x\notin  {\rm nbd}_\e(A^r)$, as  $t\to \infty$ and $x/t \to 0$}
\beq
 {\rm (i)} \quad  &&P_{\x}[\sigma_A >t] = \frac{2e_A(\x)} {\lg t}\bigg(1+ O\Big(\frac{1}{\lg t}\Big)\bigg)\quad\mbox{ for} \quad x \leq \sqrt t ,\\
 {\rm (ii)} \quad &&P_{\x}[\sigma_A \leq t] = \frac{1}{2 \lg (t/x)} \int_{x^2/2t}^\infty  e^{-y}y^{-1}dy(1+o(1)) \quad \mbox{for}  \quad  x \geq  \sqrt{t/\lg t}.
 \eeq
 
 \v2\n
 \pf \, For (i), use $\int_t^\infty  (1- e^{-x^2/2s})[s (\lg s)^2]^{-1} ds \leq Cx^2 /t(\lg t)^{2}$. See \cite{Ubes} (the proof of Theorem 15) for (ii).   \qed
\v2\v2

{\bf Acknowledgments.}\, 
 I wish to thank the anonymous  referee for his/her   pointing out many mistakes as well as  inadequate   arguments involved in the original 
  manuscript. His/Her referee reports  
 were very helpful in revising it.

\end{document}